\documentclass[12pt,a4paper,reqno]{article}
\usepackage{latexsym,amsfonts,amssymb,amsthm,upref,amscd,fontawesome5}
\usepackage[T1]{fontenc}
\usepackage[reqno]{amsmath}
 
\usepackage[mathlines]{lineno}
\newcommand*\patchAmsMathEnvironmentForLineno[1]{
	\expandafter\let\csname old#1\expandafter\endcsname\csname #1\endcsname
	\expandafter\let\csname oldend#1\expandafter\endcsname\csname end#1\endcsname
	\renewenvironment{#1}{\linenomath\csname old#1\endcsname}
	{\csname oldend#1\endcsname\endlinenomath}}
\newcommand*\patchBothAmsMathEnvironmentsForLineno[1]{\patchAmsMathEnvironmentForLineno{#1}
	\patchAmsMathEnvironmentForLineno{#1*}}
\AtBeginDocument{\patchBothAmsMathEnvironmentsForLineno{equation}
	\patchBothAmsMathEnvironmentsForLineno{align}
	\patchBothAmsMathEnvironmentsForLineno{flalign}
	\patchBothAmsMathEnvironmentsForLineno{alignat}
	\patchBothAmsMathEnvironmentsForLineno{gather}
	\patchBothAmsMathEnvironmentsForLineno{multline}}

\usepackage{times,mathrsfs,constants,booktabs,verbatim,array,xcolor}
\definecolor{softblue}{RGB}{25,25,112}   
\definecolor{softred}{RGB}{165,42,42}    
 
\usepackage{pstricks,ifthen,calc}
\usepackage[labelformat=simple,labelsep=colon]{caption}
 
\usepackage[breaklinks,colorlinks=true,linkcolor=softblue,linktoc=all,
anchorcolor=blue,
bookmarks=true,
bookmarksopen=true,
citecolor=softred,
urlcolor=magenta]{hyperref}
\usepackage{xcolor}

\usepackage[misc]{ifsym}
\usepackage{fancyhdr}
\usepackage{geometry}
\usepackage{graphicx}
\usepackage{subfigure}

 \usepackage[abbrev,lite,nobysame]{amsrefs}
\usepackage{authblk}

\newcommand{\R}{\mathbb R}
\newcommand{\RN}{\mathbb R^N}

\newcommand{\dd}{\,\mathrm d}
\newcommand{\1}{\mathbf 1}
\newcommand{\veps}{\varepsilon}

\newcommand{\cB}{\mathcal{B}}
\newcommand{\cF}{\mathcal{F}}

\newcommand{\cN}{\mathcal{N}}
\newcommand{\cL}{\mathcal{L}}

\newcommand{\xqedhere}[2]{%
	\rlap{\hbox to#1{\hfil\llap{\ensuremath{#2}}}}}

\newtheorem{Theorem}{Theorem}[section]
\newtheorem{definition}[Theorem]{Definition}
\newtheorem{lemma}[Theorem]{Lemma}
\newtheorem{proposition}[Theorem]{Proposition}

\newtheorem{remark}[Theorem]{Remark}
\newtheorem{Open Problem}[Theorem]{Open Problem}

\renewcommand{\theequation}{\arabic{section}.\arabic{equation}}

\makeatletter \@addtoreset{equation}{section} \makeatother
\newcommand{\abs}[1]{\lvert#1\rvert}

\makeatletter
\def\thanks#1{\protected@xdef\@thanks{\@thanks
		\protect\footnotetext{#1}}}
\makeatother
 
\title{\Large
	\textbf{Positive solutions to the nonlinear Choquard equation: symmetry, uniqueness, and nondegeneracy
		}
\\[0.3em]
	 
}

\author[a]{Jiankang Xia\footnote{%
	\noindent
	\hspace{-1.0em}%
	\begin{minipage}{0.6\textwidth}
		\raggedright
		\raisebox{-0.15ex}{\faEnvelope}\hspace{5pt} jiankangxia@nwpu.edu.cn \\
		\phantom{\raisebox{-0.15ex}{\faEnvelope}}\hspace{5pt} zcx@bnu.edu.cn \\
		\phantom{\raisebox{-0.15ex}{\faEnvelope}}\hspace{5pt} zhangxu0725@csu.edu.cn; \; darkblue1121@163.com
	\end{minipage}%
	}}
	
\author[b]{Chengxiang Zhang$^*$\thanks{*Corresponding Author}}
\author[c]{Xu Zhang}

\affil[a]{\footnotesize School of Mathematics and Statistics, Northwestern Polytechnical University, Xi’an 710129, China}

\affil[b]{\footnotesize Laboratory of Mathematics and Complex Systems (Ministry of Education), 
	\authorcr \hspace*{1.5cm} 
	School of Mathematical Sciences, Beijing Normal University, Beijing 100875,  China}

\affil[c]{\footnotesize  School of Mathematics and Statistics, HNP-LAMA, Central South University, Changsha 410083, China}

\date{ }

\begin{document}
	\maketitle
	

		\vskip -0.25 in
	
	\begin{minipage}{15.5cm}
		{\small {\bf Abstract:}
We study the symmetry, uniqueness, and nondegeneracy of positive solutions to the Choquard equation
$-\Delta u+u=(I_\alpha*u^p)u^{p-1}$ in $\mathbb{R}^N$, which  have long remained open problems.
For $N\geq3$, $\alpha\in(0,N)$, and $\frac{N+\alpha}{N}<p<\frac{N+\alpha}{N-2}$,
we prove that every positive $H^1$-solution is radially symmetric and strictly decreasing in the radial variable.
The symmetry proof relies on the moving-plane method, where the sublinear singularity is controlled
by analyzing the far-field behavior of the nonlocal term.
Uniqueness and nondegeneracy are established through a functional framework hinging on a nonnegative joint form and a Schur complement argument.
Radial nondegeneracy and Morse index one are obtained by ruling out the spectral interval
$(0,p-1]$ via a strict cumulative Riesz identity, which leads to uniqueness through a Nehari two-well argument.
Nondegeneracy in the full space $H^1(\mathbb{R}^N)$ follows from the spherical-harmonic decomposition
of the linearized form.
These results hold for the full range of $p$ when $\alpha=2$, and for $p\geq2$ when $\alpha\in(2,N)$.
The cases $\alpha\in(0,2)$ and $\alpha\in(2,N)$ with $p<2$ remain open.
 	}
 	

	\end{minipage}

\tableofcontents

	 \section{Introduction and main results}

In this paper, we investigate the nonlinear Choquard equation with homogeneous nonlinearity:
\begin{equation} \label{eq:choquard}
	-\Delta u + u = (I_\alpha \ast |u|^p) |u|^{p-2}u, \quad x \in \mathbb{R}^N,
\end{equation}
where $N \geq 3$, $p > 1$, and $*$ denotes the standard convolution on $\mathbb{R}^N$. The Riesz potential of order $\alpha \in (0, N)$ is given by $I_\alpha(x) = c_{N,\alpha} |x|^{\alpha-N}$, where the normalization constant $c_{N,\alpha} > 0$ is explicitly defined as
\begin{equation}\label{eq:Riesz-normalization}
	c_{N,\alpha} = \frac{\Gamma\left(\frac{N-\alpha}{2}\right)}{2^\alpha \pi^{\frac{N}{2}} \Gamma\left(\frac{\alpha}{2}\right)}.
\end{equation}
 Here $\Gamma$ denotes the Gamma function.
 Furthermore, for any $\beta, \gamma > 0$ with $\beta + \gamma < N$, the following semigroup identity holds:
\begin{equation}\label{eq:Riesz-semigroup}
	I_\beta \ast I_\gamma = I_{\beta+\gamma}.
\end{equation}
In particular, the factorization $I_\alpha = I_{\alpha/2} \ast I_{\alpha/2}$ is pivotal to the uniqueness and spectral analysis  below.

	The Choquard equation \eqref{eq:choquard} originates from fundamental problems in quantum mechanics and many-body physics. In the  prototypical case where $N = 3$, $\alpha = 2$, $p = 2$, Equation \eqref{eq:choquard} collapses to the celebrated Choquard--Pekar equation (also known as the Schr\"odinger--Newton equation):
	\begin{equation} \label{eqChPekar}
		-\Delta u + u = \Big(\frac{1}{4\pi |x|} \ast |u|^2 \Big) u, \quad x \in \mathbb{R}^3.
	\end{equation}
This nonlocal model was first formulated by Pekar \cite{Pekar1954} to characterize the polaron in condensed matter physics, capturing the self-trapped state of an electron within an ionic crystal lattice. Choquard later encountered the same equation while analyzing the Hartree--Fock approximation for a one-component plasma, employing it to illustrate the mean-field limit in many-body quantum systems (see, e.g., \cite{Lieb1977}). Moreover, Equation \eqref{eq:choquard} occupies a central position in astrophysics and quantum gravity. Often termed the Schr\"odinger--Newton equation, it governs the stationary states of self-gravitating boson stars, thereby constituting a semi-classical framework that unifies quantum mechanics with Newtonian gravitational interactions \cite{MPT98}. Across these physical contexts, the nonlocal term $(I_\alpha \ast |u|^p)$ precisely mediates the long-range Coulomb or Newtonian gravitational interactions among particles. For further physical background, we refer readers to \cite{MVS17}.

Since the seminal contributions of Lieb \cite{Lieb1977} and Lions \cite{Lions1980}, the existence of solutions to the Choquard equation \eqref{eq:choquard} has been extensively explored via variational methods. In particular, Moroz and Van Schaftingen \cite{MVS2013} established both the existence and qualitative properties including radial symmetry, monotonicity, and asymptotic decay of ground states for the equation ~\eqref{eq:choquard}. It turns out that the optimal range of exponents for the existence of finite-energy solutions is given by
\begin{equation}\label{eq:p-range}
	\frac{N+\alpha}{N}<p<\frac{N+\alpha}{N-2} 
\end{equation}
where both endpoints are critical with respect to the Hardy--Littlewood--Sobolev inequality. Although we refrain from an exhaustive bibliography here and focus primarily on positive solutions, it is worth highlighting that the asymptotic behavior of arbitrary positive solutions, particularly the algebraic decay in the regime $p<2$ has been fully characterized in our previous work \cite{XZZ2025}.

Among these properties, the uniqueness and nondegeneracy of positive solutions to the Choquard equation \eqref{eq:choquard} are of particular interest due to the 
 $H^1$-evolution dynamics of the associated Cauchy problem,
	\begin{equation*}
	\label{eq:1.3}
	\begin{aligned}
	\begin{cases}
		i\partial_t \psi + \Delta \psi + (I_\alpha \ast |\psi|^p)|\psi|^{p-2}\psi = 0, & (t, x) \in \mathbb{R} \times \mathbb{R}^N, \\
		\psi(0, x) = u(x). &
	\end{cases}
	\end{aligned}
\end{equation*}
These static properties play a crucial role in determining the orbital stability and instability of standing waves 
 $\psi(x,t) = e^{it}u(x)$. 
For instance, classification results in \cite{GS18} indicate that linearized stability of ground states holds for $p \in (1+\alpha/N, 1+(\alpha+2)/N)$, whereas instability occurs for the remaining interval of \eqref{eq:p-range}.
Central to this analysis is the nondegeneracy of positive solutions. Defining the linearized operator $L_{+,u}$ at a postive solution $u$ by
\begin{equation}\label{linearization}
L_{+,u} \varphi = -\Delta \varphi + \varphi - p (I_\alpha * (u^{p-1}\varphi)) u^{p-1} - (p-1) (I_\alpha * u^p) u^{p-2}\varphi,
\end{equation}
we say that $u$ is non-degenerate if $\ker(L_{+,u})$ is spanned exactly by the trivial elements generated by translational symmetries, i.e., 
\[
\ker(L_{+,u}) = \mathrm{span}\{ \partial_{x_1} u, \dots, \partial_{x_N} u \}.
\]
This property plays an indispensable role in nonlinear analysis: it is a crucial prerequisite for applying the Grillakis--Shatah--Strauss theory to establish orbital stability \cite{GSS87,GSS90,CGNT07}, and it guarantees the invertibility of $L_{+,u}$ on the orthogonal complement, which is necessary for constructing multi-peak solutions via Lyapunov--Schmidt reduction in perturbed problems such as the semiclassical limit.

Despite its significance, establishing symmetry, uniqueness, and nondegeneracy for the Choquard equation \eqref{eq:choquard} remains extremely difficult due to its nonlocal nature. As $\alpha \to 0$, Equation \eqref{eq:choquard} formally reduces to the local nonlinear Schr\"odinger equation $-\Delta u + u = |u|^{2p-2}u$, for which these properties are fully understood, see e.g., \cite{GNN79, MS87, Kwong89, McL93, W85}. The fundamental obstacle in the nonlocal setting is that the convolution term $(I_\alpha \ast |u|^p)$ depends on the global profile of $u$. Consequently, even after reducing to the radial setting, the nonlocal
coupling term persists, and the resulting integro-differential equation
is ill-posed as an initial value problem. This renders classical ODE
techniques and Sturm--Liouville theory inapplicable. Beyond this general obstruction, two further difficulties arise specifically for uniqueness and nondegeneracy. First, the linearized operator contains a nonlocal integral perturbation of indefinite rank, inducing a coupled spectral problem and preventing
direct application of standard Morse-index arguments for Schr\"odinger operators. Second, when $p < 2$, the effective potential term $V=(I_\alpha \ast \abs{u}^p) \abs{u}^{p-2}$ in the linearized operator is no
longer radially monotone  and positive solutions decay only algebraically. The loss of radial monotonicity obstructs strict spectral comparison. The slow algebraic decay further precludes the uniform second-order
expansions at infinity on which variational arguments rely.

	In light of the inherent challenges of the nonlocal Choquard equation, progress on symmetry, uniqueness, and nondegeneracy has been largely confined to specific parameter regimes. The foundational result on symmetry was established by Lieb \cite{Lieb1977}, who proved the radial symmetry of global minimizers for the Choquard--Pekar equation \eqref{eqChPekar} via symmetric rearrangement. Ma and Zhao \cite{MaZhao2010} subsequently classified positive solutions under certain restrictions on $p$, using the method
	of moving planes in integral forms \cite{ChenLiOu}. More recently, de Souza B\"oer, dos Santos, and de Paula Ramos \cite{BoerDosSantosRamos} removed these restrictions, extending the classification to a strongly coupled Hartree system. Nevertheless, the standard moving planes argument relies on $p\ge 2$ to ensure the boundedness of the coefficients in the resulting linear inequality for the reflected difference function; the factor $u^{p-2}$, arising from the mean value theorem, diverges when $p<2$, leaving the symmetry problem open in the sublinear regime.
	
	In this paper, we first establish the radial symmetry of positive solutions for the full range of exponents, including the sublinear regime $p<2$. Unlike Moroz and Van Schaftingen \cite{MVS2013}, who first established radial symmetry of ground states and then deduce algebraic decay, our previous nonradial decay estimate \cite{XZZ2025} demonstrates that for $p<2$, positive solutions exhibit algebraic decay independently of any symmetry assumption. Building
	on this nonradial decay, we establish a precise far-field asymptotic 
	identity that captures the balance between $u^{2-p}$ and the Hartree
	potential $I_\alpha \ast u^p$.
	Specifically, we prove the following without assuming radiality (see \eqref{eq:V-limit} of Lemma \ref{lem:regularity-and-far-field}):
	\[
	(I_\alpha*u^p)(x)u(x)^{p-2}\to 1,
	\quad \text{ as } \quad   |x|\to\infty.
	\]
Consequently, the coefficient $(p-1)(I_\alpha*u^p)u^{p-2}$ in the
reflected-difference inequality converges to $p-1 < 1$ as
$|x|\to\infty$. This strict bound allows the moving-plane argument to proceed in both its starting and continuation steps.

	
	Our first result establishes symmetry throughout the full variational
range.

	\begin{Theorem}\label{thm:symmetry}
	Let \(N\geq3\), \(0<\alpha<N\), and assume \eqref{eq:p-range}.
	If \(u\in H^1(\RN)\) is a positive weak solution of
	\eqref{eq:choquard}, then there exist \(x_0\in\RN\) and
	\(\mathcal{U}\in C^1([0,\infty))\) such that
	\[
	u(x)=\mathcal{U}(|x-x_0|),\qquad\text{for every }x\in\RN,
	\]
	and
	\[
	\mathcal{U}'(r)<0, \qquad\text{for every }r>0.
	\]
\end{Theorem}
	From this theorem, every positive $H^1$-solution is, up to translation,
	radially symmetric and strictly decreasing in the radial variable;
	we may therefore assume that all positive solutions are centered at the origin.
	Throughout the paper, we make use of the standard regularity,
	integrability, and decay properties of positive solutions to \eqref{eq:p-range};
	see~\cite{MVS2013,XZZ2025}.

We now turn to uniqueness and nondegeneracy. While these properties
are well understood when $\alpha = 2$ (the Newtonian case),
particularly for the Hartree case $p = 2$, they remain largely open for general Riesz potentials or general nonlinearities.

In the Newtonian setting, uniqueness was established by Lieb
\cite{Lieb1977} for global minimizers of the Choquard--Pekar equation
\eqref{eqChPekar} and by Choquard, Stubbe, and Vuffray \cite{CSV08}
for positive radial ground states in arbitrary dimensions via a
shooting method. Wang and Yi \cite{WangYi2017} subsequently extended
Lieb's rearrangement technique to dimensions $N=3,4,5$.
Nondegeneracy, significantly more delicate due to the coupled spectral
structure induced by the nonlocal term, was resolved by Lenzmann
\cite{Lenzmann2009} through spherical harmonic decomposition, see also \cite{WW09}, and adapted by Chen \cite{Chen2021} to
dimensions $N=3,4,5$. Together with the classification result of
Ma and Zhao \cite{MaZhao2010}, which reduces all positive solutions
to radial ones up to translation, these works completely settle the
uniqueness and nondegeneracy of positive solutions to the
Choquard--Pekar equation. However, all these approaches rely crucially
on special properties of the Newtonian potential, particularly the
spherical symmetry properties captured by Newton's theorem
\cite[Theorem 9.7]{LiebLoss}, and are restricted to the
Newtonian--Hartree case $\alpha = p = 2$.

Beyond the Newtonian case, results remain scarce and are largely
confined to perturbative regimes near known limits. Uniqueness and
nondegeneracy of ground states have been established in the endpoint
regimes of $\alpha \in (0,N)$
\cite{Seok2019}, as well as in perturbative neighborhoods of the Newtonian or Hartree limits:
$\alpha\to 2$ \cite{LZY2024,Li2022} and $p\to 2$ \cite{Xiang2016,Li2022}.
These results rely fundamentally on variational characterizations of ground states, rearrangement inequalities, or one-dimensional oscillation arguments \cite{FrankLenzmann2013,FrankLenzmannSilvestre2016}. For $\alpha = 2$ and $2 < p < (N+2)/(N-2)$, Georgiev, Tarulli, and Venkov \cite{GTV2024} recently proved local uniqueness of radial positive ground states using a method of analytic extension and showed that the radial kernel of the linearized operator is at most one-dimensional. Full nondegeneracy without symmetry assumptions, however, remains open even in this case.

A complete theory of uniqueness and nondegeneracy for general $\alpha$ and $p$ remains elusive. The fundamental obstruction is that the linearized operator $L_{+,u}$ lacks a tractable spectral structure when $\alpha \neq 2$: standard
tools such as Sturm oscillation theory and Perron--Frobenius arguments do not apply to the nonlocal eigenvalue problem. This difficulty is compounded when $p < 2$, where the effective potential loses radial
monotonicity and solutions decay only algebraically.

We establish uniqueness and nondegeneracy in the following range,
which is narrower than \eqref{eq:p-range}:
\begin{equation}\label{eq:range}
	\begin{cases}
		\displaystyle \frac{N+2}{N}<p<\frac{N+2}{N-2}, & \alpha=2,\\[1.2ex]
		\displaystyle 2\leq p<\frac{N+\alpha}{N-2}, & 2<\alpha<N.
	\end{cases}
\end{equation}

For a positive solution $u$, let $L_{+,u}$ denote the linearized
operator defined in \eqref{linearization}, regarded as a bounded map from \(H^1(\RN)\) to \(H^{-1}(\RN)\). Its associated symmetric, closed, lower-bounded form also defines a
self-adjoint operator on \(L^2(\RN)\); we use \(L_{+,u}\) for both.
The Morse index $i_{-,\mathcal{V}}(L_{+,u})$ denotes the maximal dimension of
a subspace of $\mathcal{V}$ on which the quadratic form is negative definite,
and $i_{0,\mathcal{V}}(L_{+,u})$ the dimension of $\ker(L_{+,u}) \cap \mathcal{V}$;
see Definition~\ref{def:index-nullity}. 

\begin{Theorem}\label{thm:main}
	Let $N\geq3$, and assume $\alpha,p$ satisfy \eqref{eq:range}.
	Then \eqref{eq:choquard} has exactly one positive $H^1$-solution
	up to translations. Moreover, the positive solution $u$ is
	nondegenerate modulo translations. More precisely, after translating
	its centre to the origin, its full-space linearization satisfies
	\begin{equation}\label{eq:full-nondegeneracy}
		i_{-,H^1(\RN)}(L_{+,u})=1,
		\qquad
		\ker L_{+,u}=\operatorname{span}\{\partial_{x_1}u,\ldots,\partial_{x_N}u\}.
	\end{equation}
\end{Theorem}

The Newtonian case $\alpha=2$ is covered throughout the
full variational range \eqref{eq:p-range}, whereas for $\alpha>2$
the condition $p\geq2$ is required. This restriction arises from
successive steps in the spectral proof, not from the symmetry
theorem or the half-potential reduction; the proof-idea discussion
below identifies where each condition enters.
\vskip 0.1in
	
We now outline the proof of Theorem~\ref{thm:main}.
Uniqueness rests on two ingredients: the nondegeneracy and
Morse-index-one property of $L_{+,u}$ in $H^1_r(\mathbb{R}^N)$
(Theorem~\ref{thm:index-one}), and a Nehari two-well argument
(Sec.~\ref{sec:nehari}).  Specifically,
Theorem~\ref{thm:index-one} implies that the quadratic form
$\mathcal{Q}(\phi):=\langle L_{+,u}\phi,\phi\rangle$ is coercive on
the tangent space
$T_u\cN:=\{\phi\in H^1_r(\mathbb R^N)\mid \langle L_{+,u}u,\phi\rangle=0\}$
of the Nehari manifold $\cN$. The coercivity requires splitting the
form into a coercive part and compact perturbations, which holds when
$p\geq2$ or $\alpha\geq2$; for $p<2$, an exact Bregman identity and
a tail estimate supply the quadratic lower bound instead.  Thus every
positive critical point is a strict local Nehari well.  If two distinct
positive radial critical points existed, a path joining them in $\cN$
would force the minimax level above both energies; the critical point
at that level would itself be a strict local well, pushing the level
still higher, yielding a contradiction (see \S ~\ref{subsection82}).

To establish Theorem~\ref{thm:index-one}, we introduce the joint
two-field form $Q$ with self-adjoint operator $\cL_r$ and its scalar
Schur complement $S_\lambda$.  Set $A=I_\alpha*u^p$ and
$W=I_{\alpha/2}*u^p$.  The half-potential factorization
$I_\alpha=I_{\alpha/2}*I_{\alpha/2}$, together with the substitutions
$\phi=uf$, $w=Wg$, rewrites the nonlocal interaction as a coupling
between two fields via Markov operators
\[
R:L^2_r(\dd\mu)\to L^2_r(\dd\nu) \quad\text{and}\quad
R^*:L^2_r(\dd\nu)\to L^2_r(\dd\mu),
\]
where $\dd\mu=Au^p\dd x$ and $\dd\nu=W^2\dd y$.  Completing the
square in the second field produces
\[
S_\lambda=L-\lambda I-\frac{p}{p-\lambda}R^*R,
\]
where $L$ is the unique nonnegative self-adjoint operator on
$L^2_r(\dd\mu)$ associated with
$\langle uf,uh\rangle=\langle Lf,h\rangle$; see
\eqref{eq:L-representation}.  A direct computation shows
${\rm ind}(\cL_r-\lambda)={\rm ind}\, S_\lambda$ and
${\rm nul}(\cL_r-\lambda)={\rm nul}\, S_\lambda$ for all
$\lambda<p$, so the spectral properties of $\cL_r$ transfer to
$S_\lambda$.  The proof then reduces to establishing
$\sigma(\cL_r)\cap(0,p-1]=\emptyset$, carried out in three steps.

\textit{Step 1: Fredholm properties.}
The index, nullity and Fredholm properties transfer between $\cL_r$
and $S_\lambda$; see Proposition~\ref{prop:joint-spectrum}. This reduction holds for all $0<\alpha<N$ in
\eqref{eq:p-range}; compactness and Fredholm conclusions require
$p\geq2$ or $\alpha\geq2$.

\textit{Step 2: Spectral gap via monotone zero mode.}
To rule out spectrum in $(0,p-1]$, suppose for contradiction that
$\lambda$ is the smallest spectral point there.  Radial monotonicity
is preserved by the coupling through the total positivity of the
radial Riesz kernel, which relies on the log-concavity in
Lemma~\ref{lem:log-concavity} and hence requires $\alpha\geq2$.
A one-dimensional semigroup argument for the local operator, combined
with the Trotter product formula, shows that $e^{-sS_\lambda}$
preserves the cone of radial nonincreasing functions.  The zero mode,
obtained as the limit of $e^{-sS_\lambda}h$ for a suitable
nonincreasing $h$ as $s\to\infty$, is therefore a nonzero radial
nonincreasing $f\in\ker S_\lambda$ with
$\int_{\RN}f\,\dd\mu=0$.

\textit{Step 3: Strict-sign contradiction.}
Define $X=-ru'/u$ and $Y=-rA'/A$.  The cumulative Riesz identity,
its dual form, and the scaling identity
\eqref{eq:L-scaling} ($L_{+,u}(u\Lambda)=-2u$) together yield
a strict sign for the pairing of $f$ with $X+Y/(p-\lambda)$.
The Schur-pairing identity then gives
\[
0\leq-(p-1-\lambda)
\int_{\RN}\Bigl(X+\frac{Y}{p-\lambda}\Bigr)f\,\dd\mu
=2\int_{\RN}u^2f\,\dd x<0,
\]
a contradiction.  The mechanism behind the strict sign depends on
the parameter range: for $p\geq2$, it follows from the
monotonicity of $V=Au^{p-2}$ and Chebyshev's inequality; for $p<2$,
this monotonicity fails and the Newtonian cumulative
inequality $-\Delta A=u^p$ is used instead, forcing $\alpha=2$.
This is precisely why the proof does not cover $\alpha>2$, $p<2$.

It follows that $(0,p-1]\cap\sigma(\cL_r)=\emptyset$.  Since $\cL_r$
is nonnegative with $\ker\cL_r=\operatorname{span}\{(1,1)\}$
(see Lemma~\ref{lem:joint-operator}), shifting
by $p-1$ leaves exactly one negative eigenvalue and no zero mode.
The Schur-complement transfer and the ground-state transform then give
\begin{equation}\label{indexoneandzerokernel}
	i_{-,H^1_r}(L_{+,u})=1 \quad\text{and}\quad i_{0,H^1_r}(L_{+,u})=0.
\end{equation}

The three parameter restrictions enter at distinct steps:
$p\geq2$ or $\alpha\geq2$ for Fredholm properties (Step~1),
$\alpha\geq2$ for total positivity (Step~2), and
$\alpha=2$ when $p<2$ for the strict-sign argument (Step~3).
Their intersection gives precisely \eqref{eq:range}.

\medskip 
Nondegeneracy in the full space $H^1(\mathbb{R}^N)$ follows from the
spherical-harmonic decomposition of the linearized form.  Writing
$\phi=\sum_k \phi_k(r)Y_k(\theta)$ decomposes $L_{+,u}$ into
independent angular channels.  The radial channel ($k=0$) is covered
by Theorem~\ref{thm:index-one}, which gives
\eqref{indexoneandzerokernel}, hence
$\ker L_{+,u}\cap H^1_r(\mathbb R^N)=\{0\}$.  The first angular
channel ($k=1$) has kernel spanned by the translation modes
$\partial_{x_i}u$, $i=1,\dots,N$; this follows by differentiating
the Euler--Lagrange equation satisfied by $u$, and these are the only
elements in the $k=1$ kernel.  All higher channels ($k\geq 2$) are
strictly positive: the angular eigenvalue $k(k+N-2)$ provides an
additional positive contribution to the quadratic form, and the
index-one property of the radial channel ensures that no negative
direction survives; see e.g., \eqref{eq:higher-channel-positive}.
Since the only elements of $\ker L_{+,u}$ in $H^1(\mathbb{R}^N)$ are
the translations $\partial_{x_1}u,\dots,\partial_{x_N}u$, the ground
state $u$ is nondegenerate.

\vskip 0.2in

The paper is organized as follows.
Section~\ref{sec:prelim} gathers preliminary analytic results:
key properties of the Bessel kernel, regularity and decay estimates
for positive solutions, and far-field behavior of potentials $A$ and $V$, all of
which hold without symmetry assumptions; the asymptotics of the radial
potential are also analyzed there under radial symmetry.
Section~\ref{sec:symmetry} establishes the radial symmetry and monotonicity
of positive solutions (Theorem~\ref{thm:symmetry}) via the integral
moving-plane method. The uniqueness part of Theorem~\ref{thm:main}
is proved in Sections~\ref{sec:half-potential}--\ref{sec:nehari}:
the half-potential Schur-complement reduction is developed in
Section~\ref{sec:half-potential}, a monotone zero mode is constructed
in Section~\ref{sec:monotone-zero-mode}, the cumulative identities and
radial spectral gap are proved in Section~\ref{sec:spectral-gap}, and
the Nehari two-well argument is completed in Section~\ref{sec:nehari}.
Section~\ref{sec:full-nondegeneracy} removes the radial restriction
and proves the full nondegeneracy.
Appendix~\ref{app:operator-facts} collects the operator-theoretic
facts used in the spectral analysis, and
Appendix~\ref{app:newtonian} contains the proof of the Newtonian
cumulative inequality.
Apart from the radial symmetry and monotonicity supplied by
Section~\ref{sec:symmetry}, the spectral argument is self-contained
and does not depend on the moving-plane proof.

\section{Preliminaries}
\label{sec:prelim}

This section collects, without any symmetry assumption, the analytic background used in the rest of the paper. Every statement below holds for arbitrary positive finite-energy solutions of \eqref{eq:choquard}; these facts are applied both in Section~\ref{sec:symmetry} and, after Theorem~\ref{thm:symmetry}, in the radial analysis of Sections~\ref{sec:half-potential}--\ref{sec:nehari}. The analysis rests on two independent tools: the Bessel potential lemma for $-\Delta+1$ (Lemma~\ref{lem:bessel-potential}) and a resolvent-tail comparison lemma (Lemma~\ref{lem:resolvent-tail}).
Their consequences for regularity, decay and far-field asymptotics are collected in
Lemma~\ref{lem:regularity-and-far-field} without radial symmetry, in addition, the asymptotics of radial potential $A$ and $V$ are analyzed.

Throughout, we adopt the Riesz normalization \eqref{eq:Riesz-normalization},
under which \(I_\alpha(x)=c_{N,\alpha}|x|^{\alpha-N}\), and the
semigroup identity \eqref{eq:Riesz-semigroup}.
For a positive finite-energy solution \(u\) of \eqref{eq:choquard},
we set
\begin{equation}\label{eq:A-def}
	A:=I_\alpha*u^p,
\end{equation}
\begin{equation}\label{eq:V-def}
	V:=Au^{p-2},
\end{equation}
and
\begin{equation}\label{eq:M-def}
M:=\int_{\RN}u^p\,\dd x>0.
\end{equation}

\medskip

Following \cite[Section 2.2]{BoerDosSantosRamos}, we define the
Bessel kernel associated with \(-\Delta+1\) by
\begin{equation}\label{eq:G-definition}
	\mathcal G(x):=\frac{1}{4\pi}\int_0^\infty
	s^{-N/2}
	\exp\left(
	-\frac{\pi|x|^2}{s}-\frac{s}{4\pi}
	\right)\dd s,
	\qquad x\in\R^N\setminus\{0\}.
\end{equation}
Equivalently, after the change of variables \(s=4\pi t\),
\[
\mathcal G(x)=\frac{1}{(4\pi)^{N/2}}
\int_0^\infty
t^{-N/2}
\exp\left(-t-\frac{|x|^2}{4t}\right)\dd t.
\]
The integral in \eqref{eq:G-definition} is finite for every
\(x\neq0\).  It defines a function
$
\mathcal G\in C^\infty(\R^N\setminus\{0\})\cap L^1(\R^N)
$
which satisfies
\begin{equation*}\label{eq:G-fundamental}
	(-\Delta+1)\mathcal G=\delta_0
	\qquad\text{in }\mathscr D'(\R^N).
\end{equation*}
These are the standard Bessel-potential properties; see
\cite[Chapter V, Section 3]{Stein}.

The representation \eqref{eq:G-definition} also shows that \(\mathcal G\)
is positive and radial.  Writing \(\mathcal G(x)=g(|x|)\), for every
\(r>0\) differentiation under the integral sign gives
\[
g'(r)
=-\frac r2\int_0^\infty
s^{-N/2-1}
\exp\left(
-\frac{\pi r^2}{s}-\frac{s}{4\pi}
\right)\dd s<0.
\]
Consequently,
\begin{equation*}\label{eq:G-facts}
	\mathcal G>0,\qquad
	\mathcal G(x)=g(|x|),\qquad
	g'(r)<0\quad(r>0),\qquad
	\int_{\R^N}\mathcal G(x)\dd x=1.
\end{equation*}
Moreover, since \(N\geq3\),
\[
\mathcal G(x)\leq C|x|^{2-N}\qquad(0<|x|\leq1),
\]
and \(\mathcal G\) decays exponentially as \(|x|\to\infty\).  In particular,
\begin{equation}\label{eq:G-Lt}
	\mathcal G\in L^t(\R^N)
	\qquad\text{for every }
	1\leq t<\frac{N}{N-2}.
\end{equation}

We shall use the following standard Bessel-potential theorem; see
\cite[Chapter V, Section 3.3]{Stein}.

\begin{lemma}
	\label{lem:bessel-potential}
	Let \(1<r<\infty\).  There exists a constant \(C=C(N,r)>0\) such that
	\begin{equation*}\label{eq:bessel-W2r}
		\|\mathcal G*h\|_{W^{2,r}(\R^N)}
		\leq C\|h\|_{L^r(\R^N)}
		\qquad\text{for every }h\in L^r(\R^N);
	\end{equation*}
	that is, \(h\mapsto\mathcal G*h\) is a bounded linear operator from
	\(L^r(\R^N)\) to \(W^{2,r}(\R^N)\).  Moreover, for every
	\(h\in L^r(\R^N)\) the function \(v=\mathcal G*h\) is the unique
	strong solution in \(W^{2,r}(\R^N)\) of
	\begin{equation*}\label{eq:linear-bessel}
		-\Delta v+v=h
		\qquad \text{ a.e. in }\R^N.
	\end{equation*}
\end{lemma}

Choose \(q>\max\{N,1/(p-1)\}\).  Then
\(Au^{p-1}\in L^q(\R^N)\).  Once the regularity
\(u\in W^{2,q}(\RN)\) is established
(Lemma~\ref{lem:regularity-and-far-field} below), Lemma
\ref{lem:bessel-potential} shows that \(\mathcal G*(Au^{p-1})\) is the unique
\(W^{2,q}(\R^N)\)-strong solution of
$
-\Delta v+v=Au^{p-1}\ \text{in }\R^N,
$
and consequently
\begin{equation}\label{eq:bessel-rep}
	u(x)
	=\int_{\R^N}\mathcal G(x-y)A(y)u(y)^{p-1}\dd y
	\qquad\text{for every }x\in\R^N.
\end{equation}

\medskip

The second tool is a comparison lemma that extracts the far-field
behavior of solutions to \(-\Delta z+\mu z=\mu H\) from the known
asymptotics of the right-hand side \(H\).

\begin{lemma} \label{lem:resolvent-tail}
	Let \(\mu,\beta,\ell>0\), and suppose
	$
	H(x)=\ell\abs{x}^{-\beta}+o(\abs{x}^{-\beta})
	\ \text{ as }\abs{x}\to\infty.
	$
	Let \(z\in C^2(\R^N\setminus B_R)\cap
	C(\overline{\R^N\setminus B_R})\) have bounded boundary data and
	satisfy
\[
\begin{cases}
	-\Delta z+\mu z=\mu H & \text{in } \R^N\setminus B_R, \\
	z(x)\to0 & \text{as } |x|\to \infty.
\end{cases}
\]
	Then
	\begin{equation}\label{eq:resolvent-tail}
		z(x)=\ell |x|^{-\beta}+o(|x|^{-\beta})\qquad\text{ as }|x|\to\infty.  
	\end{equation}
\end{lemma}

\begin{proof}
	Let \(h(x)=|x|^{-\beta}\).  Direct calculation gives
	\begin{equation*}
		(-\Delta+\mu)h
		=\mu h+\beta(N-2-\beta)|x|^{-\beta-2}
		=\mu h\bigl(1+O(|x|^{-2})\bigr).
	\end{equation*}
	Fix \(0<\veps<1/4\) and enlarge the exterior radius to \(R_\veps\) so
	that
	\begin{equation*}
		(1-\veps)\ell h\leq H\leq(1+\veps)\ell h
		\quad\text{and}\quad
		(1-\veps)\mu h
		\leq(-\Delta+\mu)h
		\leq(1+\veps)\mu h\quad\text{on \(\R^N\setminus B_{R_\veps}\).}
	\end{equation*}
	Let \(\Phi_\veps\) be the positive exterior solution of
\[
\begin{cases}
	\begin{aligned}
	&	(-\Delta+\mu)\Phi_\veps = 0,      && \text{in } \R^N\setminus B_{R_\veps}, \\
	&	\Phi_\veps               = 1,      && \text{on } \partial B_{R_\veps}, \\
	&	\Phi_\veps               \to 0,    && \text{as } |x|\to\infty.
	\end{aligned}
\end{cases}
\]
	Then $\Phi_\veps$ decays exponentially by the comparison principle.  From
	$(1-3\veps)(1+\veps)\leq1-\veps$ and $(1+3\veps)(1-\veps)\geq1+\veps$,
	the functions
	\[
	\underline z=(1-3\veps)\ell h-K\Phi_\veps,
	\qquad
	\overline z=(1+3\veps)\ell h+K\Phi_\veps
	\]
	are respectively a subsolution and a supersolution for the equation
	satisfied by \(z\), once \(K\) is chosen large enough to enforce the
	boundary inequalities.  The maximum principle
	\cite[Chapter 3]{GilbargTrudinger} on expanding annuli, followed
	by passage to the exterior domain, gives
	\[
	\underline z\leq z\leq\overline z.
	\]
	Since \(\Phi_\veps=o(h)\), dividing by \(\ell h\), letting
	\(|x|\to\infty\) and then \(\veps\downarrow0\), we obtain
	\eqref{eq:resolvent-tail}.
\end{proof}

\medskip

We end this section with the regularity, integrability,
decay and far-field asymptotics of positive finite-energy
solutions. Parts~\emph{(1)}--\emph{(3)} are used throughout the paper;
part~\emph{(4)} is needed only for \(p<2\), where it provides the
far-field behavior \(Au^{p-2}\to1\) in place of the monotonicity of
\(V\) available for \(p\ge2\) (see e.g., Remark~\ref{rem:V-not-monotone}).
 \begin{lemma}\label{lem:regularity-and-far-field}
 	Let \(u\) be a positive finite-energy solution of \eqref{eq:choquard}
 	with \(\alpha,p\) in the range \eqref{eq:p-range}.
 	Then
 	\(u\in C^2(\RN)\cap W^{2,r}(\RN)\) for every \(r>1\),
 	\(u^p\in L^1(\RN)\cap L^\infty(\RN)\),
 	\(V\in L^\infty(\RN)\), and
 	\(A\in C^2_{\mathrm{loc}}(\RN)\cap C^1(\RN)\cap L^\infty(\RN)\).
 	Moreover:
 	\begin{enumerate}
 		\item[(1)]
 		If \(p\geq2\), then \(V(x)\to0\) as \(|x|\to\infty\) and for every \(\mu\in(0,1)\),
 		\begin{equation}\label{eq:u-exponential-decay}
 			0<u(x)\leq C_\mu e^{-\mu|x|},
 		\end{equation}
 		and the same bound holds for \(D^\beta u\) with
 		\(|\beta|\leq2\) whenever \(\alpha\geq2\).
 		
 		\item[(2)]
 		If \(p<2\), then with \(\beta_p:=(N-\alpha)/(2-p)\),
 		\begin{equation}\label{decayp}
 			0<u(x)\leq C(1+|x|)^{-\beta_p},
 			\qquad p\beta_p>N+\alpha.
 		\end{equation}
 		
 		\item[(3)] With $M$ defined by \eqref{eq:M-def},
 		as \(|x|\to\infty\),
 		\begin{equation}\label{eq:A-tail}
 			A(x)=c_{N,\alpha}M|x|^{-(N-\alpha)}
 			+o\bigl(|x|^{-(N-\alpha)}\bigr).
 		\end{equation}
 		If \(p<2\), then
 		\begin{equation}\label{eq:V-limit}
 			\lim_{|x|\to\infty}V(x)=1,
 		\end{equation}
 		and, as \(|x|\to\infty\),
 		\begin{equation}\label{eq:balance}
 			u(x)=(c_{N,\alpha}M)^{1/(2-p)}|x|^{-\beta_p}
 			\bigl(1+o(1)\bigr).
 		\end{equation}
 	\end{enumerate}
 	If, in addition, \(u\) is radially symmetric and strictly decreasing, then:
 	\begin{enumerate}
 		\item[(4)]
 		\(A\) is radial with \(A'(r)<0\) for all \(r>0\), and
 		\begin{equation}\label{eq:A-prime-far-field}
 			A'(r)=-(N-\alpha)c_{N,\alpha}M\,r^{\alpha-N-1}
 			\bigl(1+o(1)\bigr)
 			\quad\text{ as } \quad r\to+\infty.
 		\end{equation}
 		
 		\item[(5)]
 		If \(\alpha=2\), then \(-\Delta A=u^p>0\) and
 		\begin{equation}\label{eq:A-flux}
 			-r^{N-1}A'(r)
 			=\int_0^r s^{N-1}u(s)^p\,\dd s>0
 			\quad\text{ for all } \quad r>0.
 		\end{equation}
 	\end{enumerate}
 \end{lemma}
 
 \begin{proof}
 	By \cite[Proposition 4.1]{MVS2013},
 	\(u\in C^2(\RN)\cap L^1(\RN)\cap L^\infty(\RN)\) and
 	\(u\in W^{2,r}(\RN)\) for every \(r>1\); in particular, \(u\) is
 	bounded and vanishes at infinity.
 	Hence \(u^p\in L^1(\RN)\cap L^\infty(\RN)\), and the Riesz-potential
 	mapping property gives
 	\(A\in L^\infty(\RN)\cap C(\RN)\) with \(A(x)\to0\) as
 	\(|x|\to\infty\).

 	Assume \(p\geq2\).
 	Since \(u\) is bounded and \(A(x)\to0\), one has \(V=Au^{p-2}\in L^\infty(\RN)\) and \(V\to0\) as $\abs{x}\to \infty$.
 	Given \(\mu\in(0,1)\), choose \(R>0\) such that
 	\(Au^{p-2}\leq1-\mu^2\) on \(\RN\setminus B_R\).
 	Then \(-\Delta u+\mu^2 u\leq0\) there, while
 	\((-\Delta+\mu^2)e^{-\mu|x|}
 	=(N-1)\mu|x|^{-1}e^{-\mu|x|}>0\).
 	The exterior comparison principle yields
 	\eqref{eq:u-exponential-decay}.
 	When \(\alpha\geq2\), the same exponential bound for
 	\(D^\beta u\), \(|\beta|\leq2\), follows from interior
 	second-derivative estimates applied to
 	\(-\Delta u+u=Au^{p-1}\) with \(A\in L^\infty(\RN)\).
 	
 For  Item (2), since \(p<2\), The decay estimate \eqref{decayp} is due to
 \cite[Theorem 1.1]{XZZ2025}; the strict inequality
 \(\beta_p>N\) holds because \(p>1+\alpha/N\).

 	In both the \(p\geq2\) and \(p<2\) cases, \(u^p\) possesses the
 	integrability and decay required by
 	\cite[Lemma 6.2]{MVS2013}, which yields the far-field asymptotics \eqref{eq:A-tail} of \(A\). 
 We now prove the precise asymptotic of $V$ when \(p<2\). 
 	Let \(\mathcal{G}_{2-p}\) denote the Bessel kernel of
 	\(-\Delta+2-p\), defined by
 	\[
 	\mathcal{G}_{2-p}(x)
 	:=(2-p)^{(N-2)/2}\mathcal{G}(\sqrt{2-p}\,x).
 	\]
 	Since \(A\) and \(A^{1/(2-p)}\) are bounded, vanish at infinity, and
 	belong to \(L^q(\RN)\) for all sufficiently large \(q\),
 	Lemma~\ref{lem:bessel-potential} applies to
 	\[
 	Z:=(2-p)\mathcal{G}_{2-p}*A,
 	\quad\text{and}\quad
 	Y:=(2-p)\mathcal{G}_{2-p}*A^{1/(2-p)}.
 	\]
 	In particular,
 	\[
 	(-\Delta+2-p)Z=(2-p)A,
 	\qquad
 	(-\Delta+2-p)Y=(2-p)A^{1/(2-p)},
 	\qquad\text{in }\RN,
 	\]
 	and \(Z,Y\to0\) as \(|x|\to\infty\).
 	Applying Lemma~\ref{lem:resolvent-tail} with $\mu=2-p$ and $H=A$ or
 $H=A^{1/(2-p)}$, we obtain
 \begin{equation}\label{eq:Z-over-A}
 	\frac{Z(x)}{A(x)}\to 1, \quad
 	\frac{Y(x)}{A(x)^{1/(2-p)}}\to 1,
 	\qquad\text{ as }|x|\to\infty.
 \end{equation}
 By the weighted Young's inequality, we have
 \[
 -\Delta u+u=Au^{p-1}
 \leq (2-p)A^{1/(2-p)}+(p-1)u
 =(-\Delta+2-p)Y+(p-1)u.
 \]
 It then follows that
 \[
 -\Delta u+(2-p)u
 \leq  (-\Delta+2-p)Y.
 \]
 On the other hand, since \(u>0\), the chain rule, together with
 \eqref{eq:choquard}, gives
 \begin{align*}
 	(-\Delta+2-p)u^{2-p}
 	&=(2-p)A+(2-p)(p-1)u^{-p}|\nabla u|^2
 	\geq (2-p)A=(-\Delta+2-p)Z.
 \end{align*}
 Since \(Y(x)-u(x)\to0\) and \(u(x)^{2-p}-Z(x)\to0\) as
 \(|x|\to\infty\), the comparison principle implies
 \[ u\leq Y,\quad \text{ and }\quad
 u^{2-p}\geq Z,  \quad \text{ in }\R^N.
 \]
 Combining this with \eqref{eq:Z-over-A}, we obtain
 \begin{equation*}\label{asyAu}
 	\lim_{\abs{x}\to\infty}\frac{u^{2-p}}{A(x)}=1,
 	\end{equation*}
which gives \eqref{eq:V-limit} since \(V=A/u^{2-p}\),
and yields \eqref{eq:balance} by combining this limit
with \eqref{eq:A-tail}.  Since \(V\) is continuous on \(\RN\)
with  \eqref{eq:V-limit}, one has \(V\in L^\infty(\RN)\).
 	
 	\emph{Item~(4).} 
 	Since \(I_\alpha\) and \(u^p\) are radial, so is \(A\).
 	For the strict monotonicity, we use a reflection argument.
 	Since
 	\(\partial_{x_1}|x-y|^{-(N-\alpha)}
 	=-\partial_{y_1}|x-y|^{-(N-\alpha)}\),
 	an integration by parts in \(y_1\) gives, for \(x=re_1\),
 	\begin{equation}\label{Aprime}
 	A'(r)
 	=\int_{\RN}|re_1-y|^{-(N-\alpha)}
 	\,(u^p)'(|y|)\,\frac{y_1}{|y|}\,\dd y.
 	\end{equation}
 	The integration by parts is justified by the decay of \(u\),
 	and is understood in the distributional sense when
 	\(0<\alpha\le1\), where \(\partial_{x_1}|x-y|^{-(N-\alpha)}\) is not locally integrable.
 	Splitting \(\RN\) into \(\{y_1>0\}\) and \(\{y_1<0\}\), and
 	reflecting \(y\mapsto\bar y:=(-y_1,y_2,\dots,y_N)\) in the latter,
 	we obtain
 	\[
 	A'(r)
 	=\int_{\{y_1>0\}}(u^p)'(|y|)\,\frac{y_1}{|y|}
 	\Bigl[|re_1-y|^{-(N-\alpha)}
 	-|re_1-\bar y|^{-(N-\alpha)}\Bigr]\,\dd y.
 	\]
 	For \(r>0\) and \(y_1>0\),
 	\(|re_1-\bar y|^2-|re_1-y|^2=4ry_1>0\), hence
 	\(|re_1-y|<|re_1-\bar y|\) and the bracket is strictly positive.
 	On \(\{y_1>0\}\) one has \(y_1/|y|>0\), while
 	\((u^p)'=pu^{p-1}u'<0\) by assumption.
 	The integrand is therefore strictly negative, giving
 	\(A'(r)<0\) for all \(r>0\).
 	
  For the far-field expansion of \(A'\), expand the kernel on
  \(B_{r/2}\):
  \[
  |re_1-y|^{-(N-\alpha)}
  =r^{-(N-\alpha)}\Bigl[1+(N-\alpha)\frac{y_1}{r}
  +O\bigl(|y|^2/r^2\bigr)\Bigr].
  \]
  The term independent of \(y_1\) integrates to zero by symmetry,
  while the contribution of \(\RN\setminus B_{r/2}\) is
  \(o(r^{\alpha-N-1})\) since \(|(u^p)'|\) decays faster than
  \(|y|^{-N-1}\) by Items~(1) and~(2). Note that $p\beta_p>N$. By symmetry, we deduce 
  \[
  \int_{\RN}(u^p)'(|y|)\,\frac{y_1^2}{|y|}\,\dd y
  =\frac1N\int_{\RN}(u^p)'(|y|)\,|y|\,\dd y
  =\frac{|\mathbb S^{N-1}|}{N}
  \int_0^\infty (u^p)'(t)\,t^N\,\dd t
  =-M.
  \]
  We then conclude \eqref{eq:A-prime-far-field} from \eqref{Aprime}.

 	\emph{Item~(5).} 
 	Assume \(\alpha=2\).
 	By \eqref{eq:Riesz-normalization}, \(-\Delta I_2=\delta_0\), so
 	\(-\Delta A=u^p>0\).
 	Since \(A\) is radial,
 	\(\Delta A(r)=r^{1-N}\bigl(r^{N-1}A'(r)\bigr)'\), and therefore
 	\(\bigl(r^{N-1}A'(r)\bigr)'<0\).
 	Integrating from \(0\) to \(r\), and using
 	\(\lim_{s\downarrow0}s^{N-1}A'(s)=0\), we obtain
 	\eqref{eq:A-flux}.
 \end{proof}

\section{Radial symmetry and strict radial decrease}
\label{sec:symmetry}
This section is devoted to the proof of Theorem~\ref{thm:symmetry}, which asserts that every positive finite-energy solution $u$, without any a priori assumption of radial symmetry, is radially symmetric about some point and strictly decreasing in the radial variable, for the full range of exponents specified in \eqref{eq:p-range}. The proof employs the integral formulation of the moving-plane method. The reflected differences of $ u $ and $ A $ are expressed as convolutions against the respective kernels of $ -\Delta+1 $ and $ I_\alpha $. A contraction estimate in $ L^2 $ then propagates the inequality $ u\geq u_\mu $ from the initial position $ \mu\to-\infty $ to a limiting position, provided that $ (p-1)Au^{p-2} $ is bounded by a constant strictly less than one. This bound is supplied, for $ p\geq2 $ , by the decay of $ A $ , and for $ p<2 $ , by the far-field asymptotics $ Au^{p-2}\to1 $ established in Lemma~\ref{lem:regularity-and-far-field}. Extending the argument to all directions yields a single center of symmetry common to every direction and confirms the strict monotonicity of the radial profile.

\subsection{Reflection identities and the main estimate}
\label{subsec:reflection}

Without loss of generality, we fix the \(x_1\)-direction.  For each \(\mu\in\R\), let
\[
\Sigma_\mu:=\{x=(x_1,x')\in\R^N:x_1>\mu\}
\]
denote the half-space to the right of the hyperplane \(\{x_1=\mu\}\), and set
\[
x^\mu:=(2\mu-x_1,x'),
\quad
u_\mu(x):=u(x^\mu),
\quad
A_\mu(x):=A(x^\mu).
\]
We also introduce the sets on which the reflected functions exceed the originals:
\begin{equation*}\label{eq:positive-sets}
	\Sigma_\mu^u
	:=\{x\in\Sigma_\mu:u_\mu(x)>u(x)\},
	\quad \text{ and } \quad 
	\Sigma_\mu^A
	:=\{x\in\Sigma_\mu:A_\mu(x)>A(x)\}.
\end{equation*}
\begin{lemma}\label{lem:reflection}
	For every \(x\in\Sigma_\mu\),
	\begin{align}
		u_\mu(x)-u(x)
		&=
		\int_{\Sigma_\mu}
		\bigl(\mathcal G(x-y)-\mathcal G(x^\mu-y)\bigr)
		\bigl(
		A_\mu(y)u_\mu(y)^{p-1}
		-A(y)u(y)^{p-1}
		\bigr)\dd y,
		\label{eq:u-reflection}\\
		A_\mu(x)-A(x)
		&=
		\int_{\Sigma_\mu}
		\bigl(I_\alpha(x-y)-I_\alpha(x^\mu-y)\bigr)
		\bigl(u_\mu(y)^p-u(y)^p\bigr)\dd y.
		\label{eq:A-reflection} 
	\end{align}
	If \(x\neq y\in\operatorname{int}\Sigma_\mu\), then
	\begin{align}
		0&<\mathcal G(x-y)-\mathcal G(x^\mu-y)<\mathcal G(x-y),
		\label{eq:G-difference}\\
		0&<I_\alpha(x-y)-I_\alpha(x^\mu-y)<I_\alpha(x-y).
		\label{eq:I-difference}
	\end{align}
\end{lemma}

\begin{proof}
	Split the integral in \eqref{eq:bessel-rep} over
	\(\Sigma_\mu\) and its complement.  In the complement we
	make the change of variables \(z=y^\mu\).  Using the reflection
	symmetry of \(\mathcal G\), we obtain
	\[
	u(x)=\int_{\Sigma_\mu}
	\bigl(
	\mathcal G(x-y)A(y)u(y)^{p-1}
	+
	\mathcal G(x^\mu-y)A_\mu(y)u_\mu(y)^{p-1}
	\bigr)\dd y.
	\]
	Replacing \(x\) by \(x^\mu\) gives the analogous formula for
	\(u_\mu(x)\);subtracting the two identities gives 
	\eqref{eq:u-reflection}.  Applying the same argument to
	\(A=I_\alpha*u^p\) establishes \eqref{eq:A-reflection}.

Observe that for \(x,y\in\operatorname{int}\Sigma_\mu\),
	\[
	|x^\mu-y|^2-|x-y|^2
	=4(x_1-\mu)(y_1-\mu)>0.
	\]
	The strict radial decrease of \(\mathcal G\) and \(I_\alpha\) proves
	\eqref{eq:G-difference}--\eqref{eq:I-difference}.
\end{proof}

\begin{lemma}\label{lem:pointwise}
	If \(1<p<2\), then, for every \(x\in\Sigma_\mu^u\),
	\begin{align}
		0<u_\mu(x)-u(x)
		\leq{}&
		(p-1)\int_{\Sigma_\mu^u}
		\mathcal G(x-y)A(y)u(y)^{p-2}
		\bigl(u_\mu(y)-u(y)\bigr)\dd y
		\notag\\
		&+
		\int_{\Sigma_\mu^A}
		\mathcal G(x-y)u_\mu(y)^{p-1}
		\bigl(A_\mu(y)-A(y)\bigr)\dd y.
		\label{eq:u-pointwise-sublinear}
	\end{align}
	If \(p\geq2\), then, for every \(x\in\Sigma_\mu^u\),
	\begin{align}
		0<u_\mu(x)-u(x)
		\leq{}&
		(p-1)\int_{\Sigma_\mu^u}
		\mathcal G(x-y)A(y)u_\mu(y)^{p-2}
		\bigl(u_\mu(y)-u(y)\bigr)\dd y
		\notag\\
		&+
		\int_{\Sigma_\mu^A}
		\mathcal G(x-y)u_\mu(y)^{p-1}
		\bigl(A_\mu(y)-A(y)\bigr)\dd y.
		\label{eq:u-pointwise-superlinear}
	\end{align}
	In both cases, for every \(x\in\Sigma_\mu^A\),
	\begin{equation}\label{eq:A-pointwise}
		0<A_\mu(x)-A(x)
		\leq p\int_{\Sigma_\mu^u}
		I_\alpha(x-y)u_\mu(y)^{p-1}
		\bigl(u_\mu(y)-u(y)\bigr)\dd y.
	\end{equation}
\end{lemma}

\begin{proof}
	We begin with the algebraic decomposition
	\[
	A_\mu u_\mu^{p-1}-Au^{p-1}
	=
	A(u_\mu^{p-1}-u^{p-1})
	+u_\mu^{p-1}(A_\mu-A),
	\]
	which will be substituted into \eqref{eq:u-reflection}.
From the definitions of \(\Sigma_\mu^u\) and \(\Sigma_\mu^A\), we see that
		\begin{equation}\label{eq:complement-signs1}
		u_\mu^{p-1}-u^{p-1}\leq0 \quad\text{on } \Sigma_\mu\setminus\Sigma_\mu^u,
		\qquad
		A_\mu-A\leq0 \quad\text{on } \Sigma_\mu\setminus\Sigma_\mu^A.
	\end{equation}
Now let \(x\in\Sigma_\mu^u\).  Combining the decomposition above
with \eqref{eq:u-reflection}, \eqref{eq:A-reflection}, and \eqref{eq:I-difference}, we obtain
	\begin{equation}\label{eq28}
		\begin{aligned}
			0<u_\mu(x)-u(x)
			&=
			\int_{\Sigma_\mu}
			\bigl(\mathcal G(x-y)-\mathcal G(x^\mu-y)\bigr)
			A(y)\bigl(u_\mu(y)^{p-1}-u(y)^{p-1}\bigr)\dd y
			\\
			&\quad +
			\int_{\Sigma_\mu}
			\bigl(\mathcal G(x-y)-\mathcal G(x^\mu-y)\bigr)
			u_\mu(y)^{p-1}
			\bigl(A_\mu(y)-A(y)\bigr)\dd y\\
			&\leq 
			\int_{\Sigma_\mu^u}
			\bigl(\mathcal G(x-y)-\mathcal G(x^\mu-y)\bigr)
			A(y)\bigl(u_\mu(y)^{p-1}-u(y)^{p-1}\bigr)\dd y
			\\
			&\quad +
			\int_{\Sigma_\mu^A}
			\bigl(\mathcal G(x-y)-\mathcal G(x^\mu-y)\bigr)
			u_\mu(y)^{p-1}
			\bigl(A_\mu(y)-A(y)\bigr)\dd y.
		\end{aligned}
	\end{equation}
Here the inequality holds because, by \eqref{eq:G-difference},
the kernel difference \(\mathcal G(x-y)-\mathcal G(x^\mu-y)\) is
nonnegative on \(\Sigma_\mu\), while the integrands over the
complements \(\Sigma_\mu\setminus\Sigma_\mu^u\) and
\(\Sigma_\mu\setminus\Sigma_\mu^A\) are nonpositive by
\eqref{eq:complement-signs1}.
	
	We next bound \(u_\mu^{p-1}-u^{p-1}\) on \(\Sigma_\mu^u\),
	where \(u_\mu>u>0\).  The fundamental theorem of calculus gives
	\[
	u_\mu^{p-1}-u^{p-1}
	=(p-1)\int_u^{u_\mu}t^{p-2}\dd t.
	\]
When \(1<p<2\), the map \(t\mapsto t^{p-2}\) is strictly decreasing, we thus have
\[
u_\mu^{p-1}-u^{p-1}
\leq
(p-1)u^{p-2}(u_\mu-u)
\quad\text{on }\Sigma_\mu^u.
\]
When \(p\geq2\), the map \(t\mapsto t^{p-2}\) is nondecreasing,  we then deduce
\[
u_\mu^{p-1}-u^{p-1}
\leq
(p-1)u_\mu^{p-2}(u_\mu-u)
\quad\text{on }\Sigma_\mu^u.
\]
Substituting these estimates into \eqref{eq28} yields
\eqref{eq:u-pointwise-sublinear} and
\eqref{eq:u-pointwise-superlinear}, respectively.

It remains to establish \eqref{eq:A-pointwise}. 
On \(\Sigma_\mu\setminus\Sigma_\mu^u\) we have
\(u_\mu^p-u^p\leq0\), while on \(\Sigma_\mu^u\),
\[
u_\mu^p-u^p
=p\int_u^{u_\mu}t^{p-1}\dd t
\leq p u_\mu^{p-1}(u_\mu-u).
\]
Inserting these bounds into \eqref{eq:A-reflection} and invoking \eqref{eq:I-difference} yields \eqref{eq:A-pointwise}.
\end{proof}

We define the auxiliary exponent
\begin{equation*}\label{eq:s-exponent}
	\hat{s}:=\frac{2N(p-1)}{\alpha}.
\end{equation*}
The assumption \(p>(N+\alpha)/N\) ensures \(\hat{s}>2\).
Throughout the proof we shall make use of the two
H\"older conjugate relations
\begin{equation}\label{eq:Holder-exponents}
	\frac{N+\alpha}{2N}
	=\frac{p-1}{\hat{s}}+\frac12,
	\qquad
	\frac12
	=\frac{p-1}{\hat{s}}+\frac{N-\alpha}{2N}.
\end{equation}
Since
$
1<\frac{2N}{N+\alpha}<\frac N\alpha,
$
the Hardy--Littlewood--Sobolev inequality
\cite[Theorem 4.3]{LiebLoss} implies that convolution with
\(I_\alpha\) maps $L^{\frac{2N}{N+\alpha}}(\R^N)$ continuously into
$L^{\frac{2N}{N-\alpha}}(\R^N)$.
Moreover, Lemma \ref{lem:regularity-and-far-field} yields
$
u,u_\mu\in L^2(\R^N)\cap L^{\hat{s}}(\R^N),
$
so that all norms appearing below are finite.

We next state a tail contraction property that will
drive the moving-plane argument.
When \(1<p<2\), Lemma \ref{lem:regularity-and-far-field} gives
$
(p-1)A(x)u(x)^{p-2}\to p-1<1
$
as \(|x|\to\infty\).
When \(p\geq2\), we have \(A(x)\to0\) and
\(u_\mu\leq\|u\|_\infty\).
In either case there exist \(R>0\) and \(\vartheta\in(0,1)\) such that
\begin{equation}\label{eq:tail-contraction}
	\begin{cases}
		(p-1)A(x)u(x)^{p-2}\leq\vartheta,
		&1<p<2,\\[1mm]
		(p-1)A(x)\|u\|_\infty^{p-2}\leq\vartheta,
		&p\geq2,
	\end{cases}
	\qquad |x|>R.
\end{equation}

\begin{lemma} \label{lem:main-estimate}
	Let \(R\) and \(\vartheta\) be as in
	\eqref{eq:tail-contraction}.  Then, for every \(p>1\) and
	\(\mu\in\R\),
	\begin{align*}
		&\|u_\mu-u\|_{L^2(\Sigma_\mu^u)}\\
		&\leq\Big[
		\vartheta
		+C\bigl\|
		A\max\{u^{p-2},u_\mu^{p-2}\}
		\bigr\|_{L^N(\Sigma_\mu^u\cap B_R)}
		+
		C\|u_\mu\|_{L^{\hat{s}}(\Sigma_\mu^A)}^{p-1}
		\|u_\mu\|_{L^{\hat{s}}(\Sigma_\mu^u)}^{p-1}
		\Big]
		\|u_\mu-u\|_{L^2(\Sigma_\mu^u)}.
	\end{align*}
	The constant \(C\) is independent of \(\mu\).
\end{lemma}

\begin{proof}
	Equations \eqref{eq:u-pointwise-sublinear} and
\eqref{eq:u-pointwise-superlinear}, together with the triangle
inequality, imply
\begin{align*}
	\|u_\mu-u\|_{L^2(\Sigma_\mu^u)}
	\leq&
	\bigl\|
	\mathcal G*\bigl(
	(p-1)A
	\max\{u^{p-2},u_\mu^{p-2}\}
	(u_\mu-u)
	\1_{\Sigma_\mu^u}
	\bigr)
	\bigr\|_{L^2(\RN)}
	\\
	&+
	\bigl\|
	\mathcal G*\bigl(
	u_\mu^{p-1}(A_\mu-A)
	\1_{\Sigma_\mu^A}
	\bigr)
	\bigr\|_{L^2(\RN)}.
\end{align*}
Here and in the sequel, $\1_\Omega$ denotes the characteristic function of the set $\Omega$.

We first estimate the second convolution. 
	By \eqref{eq:A-pointwise}, the Hardy--Littlewood--Sobolev inequality,
	and the first identity in \eqref{eq:Holder-exponents}, we obtain
	\begin{align*}
		\|A_\mu-A\|_
		{L^{\frac{2N}{N-\alpha}}(\Sigma_\mu^A)}
		&\leq
		C\|
		u_\mu^{p-1}(u_\mu-u)
		\|_
		{L^{\frac{2N}{N+\alpha}}(\Sigma_\mu^u)}
		\leq
		C\|u_\mu\|_
		{L^{\hat{s}}(\Sigma_\mu^u)}^{p-1}
		\|u_\mu-u\|_
		{L^2(\Sigma_\mu^u)}.
	\end{align*}
	Applying the second identity in \eqref{eq:Holder-exponents}
	to the product \(u_\mu^{p-1}(A_\mu-A)\), we deduce
	\begin{align*}
		\|u_\mu^{p-1}(A_\mu-A)\|_
		{L^2(\Sigma_\mu^A)}
		&\leq
		\|u_\mu\|_
		{L^{\hat{s}}(\Sigma_\mu^A)}^{p-1}
		\|A_\mu-A\|_
		{L^{\frac{2N}{N-\alpha}}(\Sigma_\mu^A)}
		\leq
		C\|u_\mu\|_
		{L^{\hat{s}}(\Sigma_\mu^A)}^{p-1}
		\|u_\mu\|_
		{L^{\hat{s}}(\Sigma_\mu^u)}^{p-1}
		\|u_\mu-u\|_
		{L^2(\Sigma_\mu^u)}.
	\end{align*}
	We now invoke Young's convolution inequality
	\cite[Theorem 4.2]{LiebLoss}: if
	\(1\leq q_1,q_2,q\leq\infty\) satisfy
	$
	1+1/q=1/{q_1}+1/{q_2},
	$
	then for \(f\in L^{q_1}(\R^N)\) and
	\(g\in L^{q_2}(\R^N)\),
	\begin{equation}\label{eq:Young-convolution}
		\|f*g\|_{L^q(\R^N)}
		\leq
		\|f\|_{L^{q_1}(\R^N)}
		\|g\|_{L^{q_2}(\R^N)}.
	\end{equation}
	Taking \(q_1=1\) and \(q_2=q=2\) in
	\eqref{eq:Young-convolution} and using \(\|\mathcal G\|_{L^1(\RN)}=1\), we derive
		\begin{equation*}\label{eq:nonlocal-convolution-L2}
		\begin{aligned}
		\bigl\|
			\mathcal G*\bigl(
			u_\mu^{p-1}(A_\mu-A)
			\1_{\Sigma_\mu^A}
			\bigr)
			\bigr\|_2
		&\leq
			\|u_\mu^{p-1}(A_\mu-A)\|_
			{L^2(\Sigma_\mu^A)}\\
		&\leq
			C\|u_\mu\|_{L^{\hat{s}}(\Sigma_\mu^A)}^{p-1}
			\|u_\mu\|_
			{L^{\hat{s}}(\Sigma_\mu^u)}^{p-1}
			\|u_\mu-u\|_
			{L^2(\Sigma_\mu^u)}.
		\end{aligned}
	\end{equation*}
	
	We now estimate the first convolution. To this end, we split the integral into the far-field part
	\(\Sigma_\mu^u\setminus B_R\) and the near-field part
	\(\Sigma_\mu^u\cap B_R\). On the far-field part, the tail contraction
	\eqref{eq:tail-contraction} controls the coefficient.
	When \(1<p<2\), we have
	\begin{align*}
		&\bigl\|
		\mathcal G*\bigl(
		(p-1)Au^{p-2}(u_\mu-u)
		\1_{\Sigma_\mu^u\setminus B_R}
		\bigr)
		\bigr\|_{L^2(\RN)}
		\leq
		\vartheta
		\|u_\mu-u\|_{L^2(\Sigma_\mu^u)}.
	\end{align*}
	When \(p\geq2\), the bound
	\(u_\mu^{p-2}\leq\|u\|_\infty^{p-2}\) replaces
	\(u^{p-2}\) in the integrand, and
	\eqref{eq:tail-contraction} yields the same estimate:
	\begin{align*}
		&\bigl\|
		\mathcal G*\bigl(
		(p-1)Au_\mu^{p-2}(u_\mu-u)
		\1_{\Sigma_\mu^u\setminus B_R}
		\bigr)
		\bigr\|_{L^2(\RN)}
		\leq
		\vartheta
		\|u_\mu-u\|_{L^2(\Sigma_\mu^u)}.
	\end{align*}
	
	On the near-field part \(\Sigma_\mu^u\cap B_R\),
	we use the integrability
	$
	\mathcal G\in L^{\frac{N}{N-1}}(\R^N)
	$
	provided by \eqref{eq:G-Lt}, together with the exponent identities
	\[
	1+\frac12
	=\frac{N-1}{N}+\frac{N+2}{2N},
	\qquad
	\frac{N+2}{2N}
	=\frac1N+\frac12.
	\]
	Young's and H\"older's inequalities then give, for \(1<p<2\),
	\begin{align*}
		&\bigl\|
		\mathcal G*\bigl(
		(p-1)Au^{p-2}(u_\mu-u)
		\1_{\Sigma_\mu^u\cap B_R}
		\bigr)
		\bigr\|_{L^2(\RN)}
		\leq
		C\|(p-1)Au^{p-2}\|_
		{L^N(\Sigma_\mu^u\cap B_R)}
		\|u_\mu-u\|_
		{L^2(\Sigma_\mu^u)},
	\end{align*}
	and for \(p\geq2\),
	\begin{align*}
		&\bigl\|
		\mathcal G*\bigl(
		(p-1)Au_\mu^{p-2}(u_\mu-u)
		\1_{\Sigma_\mu^u\cap B_R}
		\bigr)
		\bigr\|_{L^2(\RN)}
		\leq
		C\|(p-1)Au_\mu^{p-2}\|_
		{L^N(\Sigma_\mu^u\cap B_R)}
		\|u_\mu-u\|_
		{L^2(\Sigma_\mu^u)}.
	\end{align*}
Combining the preceding estimates completes the proof of the lemma.
\end{proof}

\subsection{Starting and the limit position}
\label{subsec:start}

\begin{proposition}\label{prop:start}
	There exists \(\mu_*>0\) such that
	$
	u(x)\geq u_\mu(x)
	\ \text{for every }x\in\Sigma_\mu,
	$
	whenever \(\mu\leq-\mu_*\).
\end{proposition}

\begin{proof}
Fix $R>0$ and $\vartheta\in(0,1)$ as in \eqref{eq:tail-contraction}.
Since $u(x)\to0$ as $|x|\to\infty$, and
\[
\inf_{x\in B_R}|x^\mu|
\to\infty
\quad\text{as }\mu\to-\infty,
\]
it follows that $u_\mu\to 0$ uniformly on $B_R$.
On the other handj, the positivity and continuity of $u$ yield
$\min_{B_R}u>0$.
Consequently,
\begin{equation}\label{eq:start-local-empty}
	\Sigma_\mu^u\cap B_R=\varnothing
\end{equation}
for all $\mu$ sufficiently negative.

Next, the change of variables $y=x^\mu$ gives
\begin{equation}\label{eq:start-tail}
	\|u_\mu\|_{L^s(\Sigma_\mu)}^s
	=
	\int_{\R^N\setminus\Sigma_\mu}u(y)^s\dd y
	\to 0
	\qquad\text{as }\mu\to-\infty,
\end{equation}
where the convergence is a consequence of $u\in L^s(\R^N)$.
Since $\Sigma_\mu^u\subset\Sigma_\mu$ and $\Sigma_\mu^A\subset\Sigma_\mu$,
\eqref{eq:start-tail} implies, as $\mu\to-\infty$,
\[
\|u_\mu\|_{L^s(\Sigma_\mu^u)}
+
\|u_\mu\|_{L^s(\Sigma_\mu^A)}
\to 0.
\]
In view of \eqref{eq:start-local-empty},
\[
\bigl\|
A\max\{u^{p-2},u_\mu^{p-2}\}
\bigr\|_{L^N(\Sigma_\mu^u\cap B_R)}
=0.
\]
We may therefore apply Lemma~\ref{lem:main-estimate} and invoke the
above convergence to obtain
\[
\|u_\mu-u\|_{L^2(\Sigma_\mu^u)}
\leq
\bigl(\vartheta+o(1)\bigr)
\|u_\mu-u\|_{L^2(\Sigma_\mu^u)},
\qquad\text{ as }\mu\to-\infty.
\]
Since $\vartheta<1$, this forces
$\|u_\mu-u\|_{L^2(\Sigma_\mu^u)}=0$,
hence $\Sigma_\mu^u=\varnothing$.
We conclude that
$u\geq u_\mu$ in $\Sigma_\mu$
for all $\mu$ sufficiently negative.
\end{proof}

In view of Proposition~\ref{prop:start}, we may define
\begin{equation*}\label{eq:lambda0}
	\mu_0:=
	\sup\left\{t\in\R:
	u\geq u_\mu\text{ in }\Sigma_\mu
	\text{ for every }\mu\leq t\right\}\in(-\infty,+\infty).
\end{equation*}
Proposition \ref{prop:start} ensures \(\mu_0>-\infty\).
It remains to show that $\mu_0<+\infty$.
Suppose, for a contradiction, that $\mu_0=+\infty$.
Then for every $\mu\in\R$,
\[
\int_{\Sigma_\mu}u^2\dd x
\geq\int_{\Sigma_\mu}u_\mu^2\dd x
=\int_{\R^N\setminus\Sigma_\mu}u^2\dd x.
\]
Letting $\mu\to+\infty$, the left-hand side vanishes while
the right-hand side converges to $\|u\|_{L^2(\RN)}^2>0$,
a contradiction.

\begin{proposition} \label{prop:critical}
	One has
	$
	u=u_{\mu_0}\ \text{in }\Sigma_{\mu_0}.
	$
\end{proposition}

\begin{proof}
	By continuity in \(\mu\) and the definition of \(\mu_0\),
	\begin{equation}\label{eq:critical-inequality}
		u\geq u_{\mu_0}
		\qquad\text{in }\Sigma_{\mu_0}.
	\end{equation}
	We argue by contradiction and suppose that equality does not hold
	identically.
	
	By \eqref{eq:critical-inequality},
	\(u^p-u_{\mu_0}^p\geq0\), and this difference is positive on a
	nonempty open set since \(u\not\equiv u_{\mu_0}\).  Then
	\eqref{eq:A-reflection} and \eqref{eq:I-difference} imply
	$
	A>A_{\mu_0}
	\ \text{in }\operatorname{int}\Sigma_{\mu_0}.
	$
	\[
Au^{p-1}>
A_{\mu_0}u^{p-1}
\geq
A_{\mu_0}u_{\mu_0}^{p-1}
\qquad\text{in }\operatorname{int}\Sigma_{\mu_0}.
\]
Invoking \eqref{eq:u-reflection} and \eqref{eq:G-difference},
we deduce the strict inequalities
	\begin{equation}\label{eq:strict-critical}
		u>u_{\mu_0},\qquad A>A_{\mu_0}
		\quad\text{in }\operatorname{int}\Sigma_{\mu_0}.
	\end{equation}
	From \eqref{eq:strict-critical} and the continuity of \(u\) and
	\(A\), it follows that
	$
	\1_{\Sigma_\mu^u\cap B_R}\to 0,
	\ 
	\1_{\Sigma_\mu^A\cap B_R}\to 0
	\ \text{a.e. as }\mu\downarrow\mu_0.
	$
	Moreover, it holds
	\begin{equation}\label{eq:bad-measure}
		|\Sigma_\mu^u\cap B_R|
		+
		|\Sigma_\mu^A\cap B_R|
		\to 0
		\qquad\text{as }\mu\downarrow\mu_0.
	\end{equation}
	
	Next we show that the function
$A\max\{u^{p-2},u_\mu^{p-2}\}$
is uniformly bounded on $\Sigma_\mu^u\cap B_R$
for $\mu$ near $\mu_0$.
On $\Sigma_\mu^u$ we have $u_\mu>u$, so the maximum
equals $u^{p-2}$ when $1<p<2$,
and $u_\mu^{p-2}$ when $p\geq2$.
In the former case, boundedness follows from
	$
	\min_{B_R}u>0,
	$
	while in the latter, from $u_\mu\leq\|u\|_{L^\infty}$.
	Since \(A\) is bounded, \eqref{eq:bad-measure} gives
	\begin{equation}\label{eq:local-small}
		\|
		A\max\{u^{p-2},u_\mu^{p-2}\}
		\|_{L^N(\Sigma_\mu^u\cap B_R)}
		\to 0
		\qquad\text{as }\mu\downarrow\mu_0.
	\end{equation}
	
	It remains to control the \(L^s\)-norms.  Restrict \(\mu\) to a compact
	interval containing \(\mu_0\), and choose \(\ell>0\) such that
	\(|\mu|\leq \ell\) on this interval.  Since
	$
	x^\mu=(-x_1,x')+(2\mu,0),
	$
	we have
	$\abs{x^\mu}\geq\abs{x}-2\ell$. 
The change of variables $y=x^\mu$ then yields
	\begin{equation*} 
		\int_{\Sigma_\mu\setminus B_R}u_\mu^s
		\leq
		\int_{|y|>R-2\ell}u(y)^s\dd y.
	\end{equation*}
Therefore,
	\begin{equation}\label{eq:norm-control}
		\begin{aligned}
			\|u_\mu\|_{L^s(\Sigma_\mu^u)}
			&\leq
			o_R(1)
			+
			\|u\|_\infty
			|\Sigma_\mu^u\cap B_R|^{1/s}
			=
			o_R(1)+o_\mu(1),
			\\
			\|u_\mu\|_{L^s(\Sigma_\mu^A)}
			&\leq
			o_R(1)
			+
			\|u\|_\infty
			|\Sigma_\mu^A\cap B_R|^{1/s}
			=
			o_R(1)+o_\mu(1),
		\end{aligned}
	\end{equation}
	where
	\[
	o_R(1):=
	\left(
	\int_{|y|>R-2\ell}u(y)^s\dd y
	\right)^{1/s}
	\to0
	\qquad\text{as }R\to\infty,
	\]
	and \(o_\mu(1)\to0\) as
	\(\mu\downarrow\mu_0\), with \(R\) fixed.
	
	Substituting \eqref{eq:local-small} and
	\eqref{eq:norm-control} into Lemma \ref{lem:main-estimate}, we obtain
	\[
	\|u_\mu-u\|_{L^2(\Sigma_\mu^u)}
	\leq
	\bigl[
	\vartheta+o_R(1)+o_\mu(1)
	\bigr]
	\|u_\mu-u\|_{L^2(\Sigma_\mu^u)}.
	\]
	Choosing $R$ large enough and then
	$\mu>\mu_0$ sufficiently close to $\mu_0$ ensures
	$\vartheta+o_R(1)+o_\mu(1)<1$,
	which forces
	$\|u_\mu-u\|_{L^2(\Sigma_\mu^u)}=0$
	and hence $\Sigma_\mu^u=\varnothing$.
	Consequently, there exists $\delta>0$ such that
	$u\geq u_\mu$ in $\Sigma_\mu$
	for every $\mu\in(\mu_0,\mu_0+\delta)$. Combined with the inequality $u\geq u_\mu$ for all
	$\mu\leq\mu_0$, this contradicts the definition
	of $\mu_0$.
	Therefore $u\equiv u_{\mu_0}$ in $\Sigma_{\mu_0}$.
\end{proof}

\subsection{Common center and strict radial decrease}
\label{subsec:center}

The preceding argument extends to an arbitrary direction.
For $e\in\mathbb S^{N-1}$ and $\mu\in\R$, let
\[
\Sigma_{\mu,e}:=\{x\in\R^N:x\cdot e>\mu\},
\qquad
x^{\mu,e}:=x-2(x\cdot e-\mu)e,
\]
and set $u_{\mu,e}(x):=u(x^{\mu,e})$,
$A_{\mu,e}(x):=A(x^{\mu,e})$.
Since every ingredient of the proof, namely the operator,
the integrability assumptions, and the decay of $u$,
is invariant under rigid rotations,
the moving-plane argument applies verbatim in the direction $e$.
For each $e\in\mathbb S^{N-1}$ we thus obtain a number $\mu(e)$
such that $u$ is symmetric with respect to the hyperplane
\[
H_e:=\{x\in\R^N:x\cdot e=\mu(e)\}.
\]

\begin{proposition}\label{prop:center}
	There exist \(x_0\in\R^N\) and a function
	\(\mathcal{U}:[0,\infty)\to(0,\infty)\) such that
	\[
	u(x)=\mathcal{U}(|x-x_0|)
	\qquad\text{for every }x\in\R^N,
	\]
	and
	\[
	\mathcal{U}'(r)<0
	\qquad\text{for every }r>0.
	\]
\end{proposition}

\begin{proof}
Fix $e\in\mathbb S^{N-1}$ and $\mu<\mu(e)$.
By the moving-plane construction,
\[
u\geq u_{\mu,e}
\qquad\text{in }\Sigma_{\mu,e}.
\]
We first rule out the identity $u\equiv u_{\mu,e}$.
If equality held, $u$ would be symmetric
with respect to both $\{x\cdot e=\mu\}$ and $H_e$;
composing these two reflections would produce a nontrivial
translation in the $e$-direction,
under which $u$ would be periodic, contradicting \(u(x)\to0\) as \(|x|\to\infty\) and \(u>0\).  
Hence $u\not\equiv u_{\mu,e}$, so $u-u_{\mu,e}\geq0$ and, by continuity, is strictly positive on a nonempty open subset of $\Sigma_{\mu,e}$.Applying \eqref{eq:A-reflection} and \eqref{eq:I-difference}
in the $e$-direction yields
$A>A_{\mu,e}$ in $\operatorname{int}\Sigma_{\mu,e}$.
Since $u\geq u_{\mu,e}>0$, we deduce
	\[
	Au^{p-1}>A_{\mu,e}u^{p-1}\ge
	A_{\mu,e}u_{\mu,e}^{p-1}
	\qquad\text{in }\operatorname{int}\Sigma_{\mu,e}.
	\]
Invoking \eqref{eq:u-reflection} and
\eqref{eq:G-difference} in the $e$-direction, we obtain
	$
	u>u_{\mu,e}
	\ \text{in }\operatorname{int}\Sigma_{\mu,e}.
	$
	
	Since \(u>0\), \(u\) is continuous, and \(u(x)\to0\) as
	\(|x|\to\infty\), it attains a positive global maximum at some
	\(x_0\in\R^N\). For any $e\in\mathbb S^{N-1}$, the strict comparison above,
	combined with the symmetry of $u$ about $H_e$,
	implies that $u$ is strictly increasing along every line
	parallel to $e$ before reaching $H_e$, and strictly
	decreasing thereafter. Consequently $x_0\in H_e$.
As $e$ is arbitrary, every symmetry hyperplane $H_e$ passes
through $x_0$.
It follows that $u$ is symmetric about every hyperplane
through $x_0$, and therefore
\[
u(x)=\mathcal{U}(|x-x_0|)
\]
for some strictly decreasing function $\mathcal{U}\colon[0,\infty)\to\R$.

To show that $\mathcal{U}$ is strictly decreasing, let $r>0$
and $e\in\mathbb S^{N-1}$, and set $\mu=x_0\cdot e-r$.
The plane $\{x\cdot e=\mu\}$ lies strictly before $H_e$.
Define
\[
w=u-u_{\mu,e}
\qquad\text{in }\Sigma_{\mu,e}.
\]
By the strict comparison established above,
\[
w>0,\qquad
-\Delta w+w
=Au^{p-1}-A_{\mu,e}u_{\mu,e}^{p-1}>0
\quad\text{in }\Sigma_{\mu,e},
\]
and $w=0$ on $\{x\cdot e=\mu\}$. Since \(e\) is the
	inward unit normal to \(\Sigma_{\mu,e}\), the Hopf boundary lemma
	\cite[Chapter 3]{GilbargTrudinger} gives
	$
	\partial_e w(x_0-re)>0.
	$
	Differentiating \(w=u-u_{\mu,e}\) at the fixed point
	\(x_0-re\) of the reflection, we deduce
	\[
	\partial_e w(x_0-re)
	=2\partial_e u(x_0-re)
	=-2\mathcal{U}'(r).
	\]
	Therefore \(\mathcal{U}'(r)<0\).
\end{proof}

Proposition~\ref{prop:center} completes the proof of the symmetry property stated in Theorem~\ref{thm:symmetry} for positive solutions of the Choquard equation~\eqref{eq:choquard}, for all $ N\geq 3 $ , $ \alpha\in(0,N) $ , and $ p $ satisfying the full range~\eqref{eq:p-range}.
\section{The half-potential Schur complement reduction}
\label{sec:half-potential}
In the remainder of this paper, \(u\) denotes a fixed positive solution
of \eqref{eq:choquard}.  By Theorem~\ref{thm:symmetry}, we may translate
the centre of symmetry to the origin and assume that \(u\) is radial
with \(u'(r)<0\) for all \(r>0\).  Unless stated otherwise, every
function space,  form, and operator is restricted to radial
functions. 

Throughout this section, \((\alpha,p)\) lies in the variational range
\eqref{eq:p-range} with \(0<\alpha<N\).
The half-potential factorization, the two measures, the Markov
operators and the Schur-complement identities are all valid in this
range; the Fredholm and spectral conclusions, however, require in
addition that \(p\geq2\) or \(\alpha\geq2\).
We use the subscript \(r\) to denote the radial subspace of a function
space, for instance, \(L^2_r(\dd\rho)\subset L^2(\RN,\dd\rho)\),
\(H^1_r(\RN)\subset H^1(\RN)\), and
\(C^\infty_{c,r}(\RN)\subset C^\infty_c(\RN)\) denote the subspaces consisting of radial functions. When no ambiguity arises, we omit the domain $\RN$ from the notation.

\subsection{The two-field quadratic form and the half-potential splitting}
\label{subsec:half-splitting}

By the definition \eqref{eq:A-def} of \(A=I_\alpha*u^p\), equation \eqref{eq:choquard} is equivalent to the following nonlinear system for the pair \((u,A)\),
\begin{equation}\label{eq:uA-system}
	-\Delta u+u=Au^{p-1},
	\qquad
	A=I_\alpha*u^p.
\end{equation}
Linearizing \eqref{eq:uA-system} yields the linearized two-field system
\begin{equation}\label{eq:linearized-two-field}
	\bigl(-\Delta+1-(p-1)Au^{p-2}\bigr)\phi=u^{p-1}\eta,
	\qquad
	\eta=I_\alpha*(p,u^{p-1}\phi),
\end{equation}
which is coupled through \(A=I_\alpha*u^p\) and coincides precisely with the linearization~\eqref{linearization} of  \eqref{eq:choquard}. 

The system \eqref{eq:linearized-two-field} is most naturally read off from the quadratic form
\begin{equation}\label{eq:two-field-energy}
	\mathcal C[\phi,\eta]
	=
	\int_{\RN}
	\bigl(|\nabla\phi|^2+\phi^2-(p-1)Au^{p-2}\phi^2\bigr)\dd x
	-2\int_{\RN}u^{p-1}\phi\,\eta\dd x
	+\frac1p\int_{\RN}\eta\,I_\alpha^{-1}\eta\dd x ,
\end{equation}
where \(I_\alpha^{-1}\) is understood in the usual Fourier/Riesz sense.  As written, \eqref{eq:two-field-energy} is a quadratic form in the single field \((\phi,\eta)\), and it is in that sense that \(\mathcal C\) will be used below, with \(\mathcal C[\phi,\eta]\) denoting the value of \(\mathcal C\)
at \((\phi,\eta)\). The associated symmetric bilinear form, obtained by polarization and
used when pairing \(\mathcal C\) with a test field (for instance, in
deriving the Euler--Lagrange system \eqref{eq:linearized-two-field}),
is given by
\begin{align*}
&\mathcal C\bigl((\phi,\eta),(\psi,\xi)\bigr)\\
&=
\int_{\RN}\bigl(\nabla\phi\cdot\nabla\psi+\phi\psi-(p-1)Au^{p-2}\phi\psi\bigr)\dd x
-\int_{\RN}u^{p-1}(\phi\xi+\psi\eta)\,\dd x
+\frac1p\int_{\RN}\xi\,I_\alpha^{-1}\eta\,\dd x .
\end{align*}
With this normalization \(\mathcal C(\cdot,\cdot)\) is exactly the bilinear form of the linearization operator of \eqref{eq:linearized-two-field} on \(L^2(\RN)\times L^2(\RN)\), with no constant left over: \eqref{eq:linearized-two-field} is precisely the statement that \(\mathcal C((\phi,\eta),(\psi,\xi))=0\) for every test field \((\psi,\xi)\).  Throughout the paper we use round brackets for the symmetric bilinear form and square brackets for the associated quadratic form: thus \(t(\cdot,\cdot)\) is a form and \(t[\cdot]=t(\cdot,\cdot)\) its quadratic form.

Since \eqref{eq:linearized-two-field} is precisely the Euler--Lagrange system of the quadratic form \(\mathcal C[\cdot]\), the spectrum of the operator in \eqref{eq:linearized-two-field}, in particular its kernel and its negative directions (that is, its nullity and its Morse index), is governed by the sign of \(\mathcal C[\cdot]\). Moreover, since the second equation of \eqref{eq:linearized-two-field} eliminates \(\eta\), the restriction of \(\mathcal C[\cdot]\) to that constraint is exactly the scalar linearization of \eqref{eq:choquard} at \(u\).

To treat the last term we use \(I_\alpha=I_{\alpha/2}*I_{\alpha/2}\) and set
\[
w:=I_{\alpha/2}^{-1}\eta,
\qquad\text{so that}\qquad
\eta=I_{\alpha/2}*w.
\]  Since \(I_{\alpha/2}\) is self-adjoint,
\[
\int_{\RN}(I_{\alpha/2}*w)\,I_\alpha^{-1}(I_{\alpha/2}*w)\dd x
= \int_{\RN}w\,I_{\alpha/2}*(I_\alpha^{-1}(I_{\alpha/2}*w))\dd x=
\int_{\RN}w^2\dd x .
\]
Thus the two-field form becomes the half-potential form
\begin{equation}\label{eq:E-phi-w}
	\widetilde{\mathcal C}[\phi,w]
	=
	\int_{\RN}
	\bigl(|\nabla\phi|^2+\phi^2-(p-1)Au^{p-2}\phi^2\bigr)\dd x
	-2\int_{\RN}u^{p-1}\phi\,(I_{\alpha/2}*w)\dd x
	+\frac1p\int_{\RN}w^2\dd x .
\end{equation}
The corresponding Euler--Lagrange system, read off from the polarization of \eqref{eq:E-phi-w} as above, is
\begin{equation*}\label{eq:half-step-system}
	\bigl(-\Delta+1-(p-1)Au^{p-2}\bigr)\phi=u^{p-1}I_{\alpha/2}*w,
	\qquad
	w=I_{\alpha/2}*(p\,u^{p-1}\phi).
\end{equation*}
Here \(w\) is the perturbation of the half-potential; the unperturbed half-potential is defined by
\begin{equation}\label{eq:W-def}
	W:=I_{\alpha/2}*u^p,
	\qquad\text{so that}\qquad
	A=I_{\alpha/2}*W .
\end{equation}

The rest of this section recasts this two-component system.  Its two fields will be normalized by \(u\) and by \(W\) respectively (Subsection~\ref{subsec:ground-state}).  This substitution produces a
measure for each field together with the operator that couples them, and
the normalized form is then assembled into a single nonnegative joint
form (Subsection~\ref{subsec:joint}).

\subsection{Normalization: weighted measures and Markov operators}
\label{subsec:ground-state}

We now carry out the two normalizations, \(\phi=uf\) and \(w=Wg\), substituting them in the half-potential form \eqref{eq:E-phi-w}.  Each of the three terms of \eqref{eq:E-phi-w} produces one of the structures on which the rest of this section is built.

For \(f,h\in C^\infty_{c,r}(\RN)\), expanding \(\nabla(uf)\) and \(\nabla(uh)\), integrating the mixed terms by parts, and using \(-\Delta u+u=Au^{p-1}\), we obtain
\begin{equation}\label{eq:ground-state-transform}
	\int_{\RN}u^2\nabla f\cdot\nabla h\,\dd x
	+\int_{\RN}Au^pfh\,\dd x
	=
	\langle uf,uh\rangle_{H^1_r}.
\end{equation}
Together with \(Au^{p-2}\phi^2=Au^pf^2\), this converts the first part of \eqref{eq:E-phi-w} into
\begin{equation}\label{eq:local-part-normalized}
	\int_{\RN}\bigl(|\nabla(uf)|^2+u^2f^2-(p-1)Au^{p-2}u^2f^2\bigr)\dd x
	=
	\langle uf,uf\rangle_{H^1_r}+(1-p)\int_{\RN}f^2\dd\mu,
\end{equation}
where \(\dd\mu(x):=A(x)u(x)^p\dd x\).

The last term of \eqref{eq:E-phi-w} is handled by the second normalization alone: with \(w=Wg\),
\[
\frac1p\int_{\RN}w^2\dd x
=
\frac1p\int_{\RN}W^2g^2\dd x
=
\frac1p\int_{\RN}g^2\dd\nu,
\qquad
\dd\nu(y):=W(y)^2\dd y .
\]
The two normalizations therefore single out the two measures
\begin{equation}\label{eq:joint-measures}
	\dd\mu(x)=A(x)u(x)^p\dd x,
	\qquad
	\dd\nu(y)=W(y)^2\dd y,
\end{equation}
and every space, norm and operator below is an object of the two weighted spaces \(L^2_r(\dd\mu)\) and \(L^2_r(\dd\nu)\).

It remains to read off the coupling.  Since \(u^{p-1}\phi=u^pf\), the middle term of \eqref{eq:E-phi-w} is \(-\int_{\RN}u^pf\,(I_{\alpha/2}*(Wg))\dd x\), a double integral of the product \(f(x)g(y)\); expressed over the product space, it becomes an integral of the product \(f(x)g(y)\) against the positive measure
\begin{equation}\label{eq:joint-coupling-measure}
	\dd\pi(x,y)
	=
	I_{\alpha/2}(x-y)u(x)^pW(y)\dd x\dd y
\end{equation}
on \(\RN\times\RN\).  By \eqref{eq:W-def}, integrating \eqref{eq:joint-coupling-measure} first in \(y\)  yields \(u(x)^p(I_{\alpha/2}*W)(x)\dd x=\dd\mu(x)\), and first in \(x\) gives \(W(y)(I_{\alpha/2}*u^p)(y)\dd y=\dd\nu(y)\): the two marginals of \(\pi\) are precisely the measures of \eqref{eq:joint-measures}, as Lemma~\ref{lem:half-potential-marginals} below establishes.  The coupling thus determines a pairing of \(L^2_r(\dd\mu)\) with \(L^2_r(\dd\nu)\), and it defines operators \(R\) and \(R^*\), the two conditional expectations of \(\pi\), through
\[
-\iint_{\RN\times\RN}f(x)g(y)\dd\pi(x,y)
=
-\langle f,R^*g\rangle_{L^2_r(\dd\mu)}
=
-\langle Rf,g\rangle_{L^2_r(\dd\nu)}.
\]
Explicitly, they take the form
\begin{equation}\label{eq:Markov-operators}
	Rf(y)=\frac{I_{\alpha/2}*(u^pf)(y)}{W(y)},
	\qquad
	R^*g(x)=\frac{I_{\alpha/2}*(Wg)(x)}{A(x)}.
\end{equation}
Both formulas are first verified on the radial test fields \(C^\infty_{c,r}(\RN)\), which are dense in \(L^2_r(\dd\mu)\) and \(L^2_r(\dd\nu)\), respectively, and on which the two convolutions converge absolutely.  The contraction estimates of Lemma~\ref{lem:Markov-contractions} below then extend them uniquely to bounded operators
\[
R:L^2_r(\dd\mu)\longrightarrow L^2_r(\dd\nu),
\qquad
R^*:L^2_r(\dd\nu)\longrightarrow L^2_r(\dd\mu).
\]
In particular \(R^*\) is the Hilbert-space adjoint of \(R\).  Altogether, \(\widehat{\mathcal C}[f,g]:=\widetilde{\mathcal C}[uf,Wg]\) assumes the form
\begin{equation}\label{eq:E-f-g}
	\widehat{\mathcal C}[f,g]
	=
	\langle uf,uf\rangle_{H^1_r}
	+
	(1-p)\|f\|_{L^2_r(\dd\mu)}^2
	-
	2\langle Rf,g\rangle_{L^2_r(\dd\nu)}
	+
	\frac1p\|g\|_{L^2_r(\dd\nu)}^2 .
\end{equation}

The marginals of \eqref{eq:joint-coupling-measure} and the total masses of the two measures are characterized as follows.
\begin{lemma}\label{lem:half-potential-marginals}
	The measures \(\mu\) and \(\nu\) are finite and have the same total mass.  They are the two marginals of \(\pi\):
	\begin{equation}\label{eq:joint-marginals}
		\pi(\cB_1\times\RN)=\mu(\cB_1),
		\qquad
		\pi(\RN\times \cB_2)=\nu(\cB_2)
	\end{equation}
	for all Borel sets \(\cB_1,\cB_2\subset\RN\).  Moreover,
	\begin{equation}\label{eq:equal-joint-masses}
		\mu(\RN)=\nu(\RN)
		=
		\int_{\RN}(I_\alpha*u^p)u^p\dd x<\infty.
	\end{equation}
\end{lemma}

\begin{proof}
	Using \eqref{eq:Riesz-semigroup} and applying Tonelli's theorem, we obtain
	\[
	\begin{aligned}
		\int_{\RN}W(y)^2\dd y
		&=\iiint_{\RN\times\RN\times\RN}
		I_{\alpha/2}(y-x)I_{\alpha/2}(y-z)
		u(x)^pu(z)^p\dd x\dd z\dd y\\
		&=\iint_{\RN\times\RN}I_\alpha(x-z)u(x)^pu(z)^p
		\dd x\dd z
		=\int_{\RN}Au^p\dd x.
	\end{aligned}
	\]
	The last quantity is finite by the Hardy--Littlewood--Sobolev inequality and the strict variational range. Thus \(\mu\) and \(\nu\) are finite and \eqref{eq:equal-joint-masses} holds.
	
	Integrating \eqref{eq:joint-coupling-measure} first in \(y\) and invoking \eqref{eq:Riesz-semigroup} yields
	\[
	\int_{\RN}I_{\alpha/2}(x-y)W(y)\dd y
	=(I_{\alpha/2}*W)(x)=(I_{\alpha/2}*I_{\alpha/2}*u^p)(x)=A(x),
	\]
	and hence the first identity in \eqref{eq:joint-marginals}.  Integrating first in \(x\) produces  \(W(y)^2\dd y\) and confirms the second.
\end{proof}

Next, we summarize the elementary properties of the operators \eqref{eq:Markov-operators}, which ensure that they are Markov operators.
\begin{lemma}\label{lem:Markov-contractions}
	The operators \(R:L^2_r(\dd\mu)\to L^2_r(\dd\nu)\) and \(R^*:L^2_r(\dd\nu)\to L^2_r(\dd\mu)\) of \eqref{eq:Markov-operators} satisfy
	\[
	R1=R^*1=1,\qquad
	Rf\geq0\;\; \text{ if }\; f\geq0,\qquad
	R^*g\geq0\;\; \text{ if }\; g\geq0,
	\]
	and
	\[
	\|Rf\|_{L^2_r(\dd\nu)}\leq\|f\|_{L^2_r(\dd\mu)},
	\qquad
	\|R^*g\|_{L^2_r(\dd\mu)}\leq\|g\|_{L^2_r(\dd\nu)}.
	\]
	Moreover, \(R^*\) is the Hilbert-space adjoint of \(R\).  We define the operator
	\[
	P:=R^*R:L^2_r(\dd\mu)\longrightarrow L^2_r(\dd\mu).
	\]
	Then \(P\) is positivity preserving and self-adjoint, and it satisfies \(P1=1\) and \(\|P\|\leq1\).  It admits the explicit expression
	\begin{equation}\label{eq:P-Markov}
		Pf=\frac{I_\alpha*(u^pf)}{A}.
	\end{equation}
\end{lemma}

\begin{proof}
	We first prove that \(R\) preserves constants and positivity.  By the definition of \(W\),
	\[
	R1(y)
	=
	\frac{1}{W(y)}
	\int_{\RN}I_{\alpha/2}(x-y)u(x)^p\,\dd x
	=1.
	\]
	Since \(I_{\alpha/2}>0\), \(u>0\), and \(W>0\), the definition of \(R\) also entails
	\[
	f\geq0 \Longrightarrow Rf\geq0.
	\]
	In particular, \(0\leq f\leq1\) implies \(0\leq Rf\leq1\).
	
	We next establish the \(L^2\)-contraction estimate, first on \(C^\infty_{c,r}(\RN)\).  By the Cauchy--Schwarz inequality,
	\[
	\begin{aligned}
		|Rf(y)|^2
		&=
		\frac{1}{W(y)^2}
		\left|
		\int_{\RN}
		I_{\alpha/2}(x-y)u(x)^pf(x)\,\dd x
		\right|^2\\
		&\leq
		\frac{1}{W(y)^2}
		\left(
		\int_{\RN}I_{\alpha/2}(x-y)u(x)^p\,\dd x
		\right)
		\left(
		\int_{\RN}
		I_{\alpha/2}(x-y)u(x)^p|f(x)|^2\,\dd x
		\right)\\
		&=
		\frac{I_{\alpha/2}*(u^p|f|^2)(y)}{W(y)}
		=
		R(|f|^2)(y).
	\end{aligned}
	\]
	Multiplying by \(W(y)^2\) and integrating in \(y\), Tonelli's theorem and \eqref{eq:Riesz-semigroup} give
	\[
	\begin{aligned}
		\|Rf\|_{L^2_r(\dd\nu)}^2
		&\leq
		\int_{\RN}
		I_{\alpha/2}*(u^p|f|^2)(y)\,W(y)\,\dd y\\
		&=
		\int_{\RN}u(x)^p|f(x)|^2
		\left(
		\int_{\RN}I_{\alpha/2}(x-y)W(y)\,\dd y
		\right)\dd x\\
		&=
		\int_{\RN}A(x)u(x)^p|f(x)|^2\,\dd x
		=
		\|f\|_{L^2_r(\dd\mu)}^2.
	\end{aligned}
	\]
	Since \(C^\infty_{c,r}(\RN)\) is dense in \(L^2_r(\dd\mu)\), this estimate extends by density: \(R:L^2_r(\dd\mu)\to L^2_r(\dd\nu)\) is a contraction.
	
	We next verify that the second operator in \eqref{eq:Markov-operators} is the Hilbert-space adjoint of \(R\).  For \(f,g\in C^\infty_{c,r}(\RN)\), Fubini's theorem yields
	\[
	\begin{aligned}
		\langle Rf,g\rangle_{L^2_r(\dd\nu)}
		&=
		\iint_{\RN\times\RN}
		I_{\alpha/2}(x-y)u(x)^pf(x)W(y)g(y)
		\,\dd x\dd y\\
		&=
		\int_{\RN}f(x)
		\frac{I_{\alpha/2}*(Wg)(x)}{A(x)}
		A(x)u(x)^p\,\dd x
		=
		\langle f,R^*g\rangle_{L^2_r(\dd\mu)}.
	\end{aligned}
	\]
	Both sides are continuous in \((f,g)\in L^2_r(\dd\mu)\times L^2_r(\dd\nu)\), so by density the identity holds for all \(f\in L^2_r(\dd\mu)\) and \(g\in L^2_r(\dd\nu)\); thus \(R^*\) is the adjoint of \(R\).  The remaining properties follow similarly:
	\[
	R^*1=1,\qquad
	g\geq0\Longrightarrow R^*g\geq0,\qquad
	\|R^*\|=\|R\|\leq1.
	\]
	
	By the definition of \(P=R^*R\), and since \(W Rf=I_{\alpha/2}*(u^pf)\) by \eqref{eq:Markov-operators},
	\[
	\begin{aligned}
		Pf
		&=R^*(Rf)
		=\frac{I_{\alpha/2}*(W Rf)}{A}
		=
		\frac{
			I_{\alpha/2}*I_{\alpha/2}*(u^pf)
		}{A}
		=
		\frac{I_\alpha*(u^pf)}{A},
	\end{aligned}
	\]
	where we used \eqref{eq:Riesz-semigroup}; this is \eqref{eq:P-Markov}, the identity extending from \(f\in C^\infty_{c,r}(\RN)\) to all \(f\in L^2_r(\dd\mu)\) by density.  The positivity preservation, self-adjointness, \(P1=1\), and \(\|P\|\leq1\) follow immediately from the corresponding properties of \(R\) and \(R^*\).
\end{proof}

The normalization has a second consequence: the substitution \(\phi=uf\) is not only a rewriting of the form, but an isometry.  Let \(\mathcal D\) be the completion of \(C^\infty_{c,r}(\RN)\) for the norm
\begin{equation}\label{eq:calD-norm}
	\|f\|_{\mathcal D}^2
	=
	\int_{\RN}u^2|\nabla f|^2
	+
	\int_{\RN}f^2\dd\mu.
\end{equation}

\begin{lemma}\label{lem:ground-state-transform}
	Multiplication by \(u\) is an isometric isomorphism from \(\mathcal D\) onto \(H^1_r(\RN)\).  More precisely, the identity \eqref{eq:ground-state-transform} extends to all \(f,h\in\mathcal D\):
	\[
	\int_{\RN}u^2\nabla f\cdot\nabla h
	+
	\int_{\RN}fh\dd\mu
	=
	\langle uf,uh\rangle_{H^1_r},
	\qquad f,h\in\mathcal D,
	\]
	and
	\begin{equation}\label{eq:calD-realization}
		\mathcal D=
		\bigl\{
		f\in H^1_{\mathrm{loc},r}(\RN):uf\in H^1_r(\RN)
		\bigr\}.
	\end{equation}
\end{lemma}

\begin{proof}
	By \eqref{eq:ground-state-transform}, the map
	\[
	T:C^\infty_{c,r}(\RN)\to H^1_r(\RN),
	\qquad
	Tf=uf,
	\]
	is an isometry for the norm \eqref{eq:calD-norm}.  Since \(C^\infty_{c,r}(\RN)\) is dense in \(\mathcal D\) and \(H^1_r(\RN)\) is complete, \(T\) extends uniquely to an isometry \(T:\mathcal D\to H^1_r(\RN)\); its range is closed, because \(\mathcal D\) is complete.
	
	Since \(u\in L^\infty(\RN)\) and \(Au^p\in L^\infty(\RN)\) by Lemma~\ref{lem:regularity-and-far-field}, there exists a constant \(C\) such that
	\[
	\|g\|_{\mathcal D}\leq C\|g\|_{H^1_r}
	\qquad\text{for every }g\in C^\infty_{c,r}(\RN).
	\]
	
	We first prove that \(T\) is onto.  Let \(\phi\in C^\infty_{c,r}(\RN)\).  Since \(u>0\) is continuous, \(\phi/u\in H^1_r(\RN)\), and by density there exist \(f_n\in C^\infty_{c,r}(\RN)\) with \(f_n\to\phi/u\) in \(H^1_r(\RN)\).  The estimate above shows that \((f_n)\) is Cauchy in \(\mathcal D\), so \(f_n\to f\) in \(\mathcal D\) for some \(f\in\mathcal D\); since \(T\) is continuous and \(Tf_n=uf_n\to\phi\), we infer \(Tf=\phi\).  Hence \(C^\infty_{c,r}(\RN)\subseteq T(\mathcal D)\), and as \(T(\mathcal D)\) is closed while \(C^\infty_{c,r}(\RN)\) is dense in \(H^1_r(\RN)\),
we thus conclude $T(\mathcal D)=H^1_r(\RN)$: multiplication by \(u\) is an isometric isomorphism of \(\mathcal D\) onto \(H^1_r(\RN)\).  Both sides of \eqref{eq:ground-state-transform} are continuous in \((f,h)\) with respect to the norm \eqref{eq:calD-norm} and agree on the dense subset \(C^\infty_{c,r}(\RN)\); hence the identity holds for all \(f,h\in\mathcal D\).
	
	It remains to prove \eqref{eq:calD-realization}.  Fix a compact \(K\subset\RN\).  Since \(1/u\) and \(\nabla(1/u)\) are bounded on \(K\), the product rule gives
	\[
	\|g\|_{H^1(K)}
	\leq C_K\bigl(\|ug\|_{L^2(K)}+\|\nabla(ug)\|_{L^2(K)}\bigr)
	\leq C_K\|ug\|_{H^1_r}
	=C_K\|g\|_{\mathcal D}
	\]
	for every \(g\in C^\infty_{c,r}(\RN)\).  Thus \(g\mapsto g\) is continuous from \((C^\infty_{c,r}(\RN),\|\cdot\|_{\mathcal D})\) to \(H^1(K)\), and therefore extends to a bounded map \(\mathcal D\to H^1(K)\); the extensions obtained for different \(K\) coincide, since they agree on the dense subset \(C^\infty_{c,r}(\RN)\).  Hence every \(f\in\mathcal D\) is represented by a function in \(H^1_{\mathrm{loc},r}(\RN)\) (two representatives agree almost everywhere), and \(uf=Tf\in H^1_r(\RN)\).  Therefore
	\[
	\mathcal D
	\subseteq
	\bigl\{
	f\in H^1_{\mathrm{loc},r}(\RN):uf\in H^1_r(\RN)
	\bigr\}.
	\]
	Conversely, let \(h\in H^1_{\mathrm{loc},r}(\RN)\) satisfy \(uh\in H^1_r(\RN)\).  Since \(T\) is onto, \(uh=Tf\) for some \(f\in\mathcal D\); then \(uf=uh\), and \(u>0\) implies \(f=h\) almost everywhere, so \(h\in\mathcal D\). This completes the proof of \eqref{eq:calD-realization}.
\end{proof}

The form \((f,h)\mapsto\langle uf,uh\rangle_{H^1_r}\) of \eqref{eq:local-part-normalized}, with form domain
\(\mathcal D\), is the local part of the linearization.
Since \(C^\infty_{c,r}(\RN)\subset\mathcal D\) is dense in \(L^2_r(\dd\mu)\), this form is densely defined, symmetric and nonnegative; by Lemma~\ref{lem:ground-state-transform} its form norm is \(\|uf\|_{H^1_r}=\|f\|_{\mathcal D}\), so that it is closed with form domain \(\mathcal D\).

By the representation theorem for closed nonnegative forms \cite[Theorem~2.6 of Chapter~VI]{Kato}, there exists a unique nonnegative self-adjoint operator \(L\) on \(L^2_r(\dd\mu)\) such that
\begin{equation}\label{eq:L-representation}
	\langle uf,uh\rangle_{H^1_r}=\langle Lf,h\rangle_{L^2_r(\dd\mu)}
	\qquad
	\text{for all }f\in D(L),\ h\in\mathcal D.
\end{equation}
Moreover, \(L\ge I\). Indeed, by \eqref{eq:ground-state-transform} and \eqref{eq:calD-norm}, we have  for every \(f\in D(L)\),
 \[\langle Lf,f\rangle_{L^2_r(\dd\mu)}=\langle uf,uf\rangle_{H^1_r}=\|f\|^2_{\mathcal D}\geq\|f\|^2_{L^2_r(\dd\mu)}.\]
More precisely,
\begin{equation}\label{eq:L-domain}
	\begin{aligned}
		D(L)=\bigl\{f\in\mathcal D:\;&
		\text{there exists }F\in L^2_r(\dd\mu)\text{ such that}\
		\langle uf,uh\rangle_{H^1_r}
		=\langle F,h\rangle_{L_r^2(\dd\mu)}
		\text{ for every }h\in\mathcal D
		\bigr\}.
	\end{aligned}
\end{equation}
For \(f\in D(L)\), the function \(F\) in \eqref{eq:L-domain} is unique, and we define \(Lf=F\).  The following differential formula provides its distributional realization:
\[
\langle uf,uh\rangle_{H^1_r}
=
\int_{\RN}Fh\,\dd\mu
=
\int_{\RN}Au^pFh\,\dd x
\]
for every \(h\in C^\infty_{c,r}(\RN)\).  Hence
\[
-\operatorname{div}(u^2\nabla f)+Au^pf=Au^pF
\]
in the distributional sense.  Since \(F=Lf\) and \(Au^p>0\), this can be rewritten as
\[
Lf
=
-\frac{1}{Au^p}\operatorname{div}(u^2\nabla f)+f.
\]

In the operator notation of \eqref{eq:L-representation}, the Euler--Lagrange equations of the normalized quadratic form \eqref{eq:E-f-g}, with variations taken in \(L^2_r(\dd\mu)\) and \(L^2_r(\dd\nu)\), read
\begin{equation}\label{eq:normalized-linearized-system}
	\bigl(L+(1-p)I\bigr)f-R^*g=0,
	\qquad
	g=pRf .
\end{equation}

Since \(L\ge I\), the operator \(L\) is boundedly invertible.
\begin{lemma}\label{lem:L-compact}
	Assume \(p\geq2\).  Then the inclusion map
	$
	\mathcal D\hookrightarrow L^2_r(\dd\mu)
	$
	is compact. Consequently, the operator
	$
	L^{-1}:L^2_r(\dd\mu)\longrightarrow L^2_r(\dd\mu)
	$ is compact.
\end{lemma}

\begin{proof}
  With \(V=Au^{p-2}\), Lemma~\ref{lem:ground-state-transform} gives
	\begin{equation}\label{eq:weighted-identity}
		\int_{\RN}f^2\dd\mu=\int_{\RN}V(uf)^2\dd x,
		\qquad f\in\mathcal D.
	\end{equation}
	Let \((f_j)\) be bounded in \(\mathcal D\) and set \(\phi_j=uf_j\).  By Lemma~\ref{lem:ground-state-transform}, \((\phi_j)\subset H^1_r\) is bounded.  After passing to a subsequence,
	\(\phi_j\rightharpoonup\phi\) in \(H^1_r(\RN)\).  Since \(p\geq2\), Lemma~\ref{lem:regularity-and-far-field} gives  \(V\in L^\infty(\RN)\) and \(V(x)\to0\) as \(|x|\to\infty\).  Therefore multiplication by \(V^{1/2}\) is compact from \(H^1_r(\RN)\) to \(L^2(\RN)\).  Hence
	\[
	\int_{\RN}V|\phi_j-\phi|^2\dd x\to0.
	\]
	By Lemma~\ref{lem:ground-state-transform}, \(\phi/u\in\mathcal D\).  Applying \eqref{eq:weighted-identity} to \(f_j-\phi/u\), we deduce
	\[
	\left\|f_j-{\phi}/{u}\right\|_{L^2_r(\dd\mu)}^2
	=
	\int_{\RN}V|\phi_j-\phi|^2\dd x
	\to0.
	\]
	Thus every bounded sequence in \(\mathcal D\) has a subsequence converging in \(L^2_r(\dd\mu)\), so the inclusion \(\mathcal D\hookrightarrow L^2_r(\dd\mu)\) is compact.
	
	It remains to prove the assertion about the resolvent.  Let \((g_j)\) be bounded in \(L^2_r(\dd\mu)\), and set $f_j=L^{-1}g_j$.
	By the definition of the operator associated with the form,
	\[
	\langle uf_j,uh\rangle_{H^1_r}
	=
	\int_{\RN}g_jh\dd\mu, 
	\qquad\text{ for every }h\in\mathcal D.
	\]
	Taking \(h=f_j\), we find
	\[
	\|f_j\|_{\mathcal D}^2
	=
	\int_{\RN}g_jf_j\dd\mu
	\leq
	\|g_j\|_{L^2_r(\dd\mu)}
	\|f_j\|_{L^2_r(\dd\mu)}
	\leq
	\|g_j\|_{L^2_r(\dd\mu)}
	\|f_j\|_{\mathcal D}.
	\]
	Therefore \((f_j)\) is bounded in \(\mathcal D\).  By the compact inclusion proved above, \((f_j)\) has a subsequence converging in \(L^2_r(\dd\mu)\).  Hence \(L^{-1}\) is compact, and therefore \(L\) has compact resolvent.
\end{proof}

\subsection{The joint quadratic form and the joint operator}
\label{subsec:joint}

We now pass from the normalized form \eqref{eq:E-f-g} to a manifestly nonnegative joint quadratic form.  The two fields are paired in the space
\begin{equation}\label{eq:joint-space}
	\mathcal H=L^2_r(\dd\mu)\oplus L^2_r(\dd\nu),
	\qquad
	\langle(f,g),(h,k)\rangle_{\mathcal H}
	=
	p\langle f,h\rangle_{L^2_r(\dd\mu)}+\langle g,k\rangle_{L^2_r(\dd\nu)}.
\end{equation}

The transformation of the quadratic form amounts to substituting the joint measure \(\pi\) of \eqref{eq:joint-coupling-measure} into the coupling of \eqref{eq:E-f-g}.  Since \(\mu\) and \(\nu\) are the two marginals of \(\pi\), expanding the square gives
\begin{equation}\label{eq:joint-dirichlet-expansion}
	p\iint_{\RN\times\RN}(f(x)-g(y))^2\dd\pi(x,y)
	=
	p\|f\|_{L^2_r(\dd\mu)}^2
	+
	p\|g\|_{L^2_r(\dd\nu)}^2
	-
	2p\langle Rf,g\rangle_{L^2_r(\dd\nu)}.
\end{equation}  Substituting \eqref{eq:joint-dirichlet-expansion} into \eqref{eq:E-f-g} and using \eqref{eq:joint-space} converts the normalized form into
\begin{equation}\label{eq:Cn-transformed}
	p\,\widehat{\mathcal C}[f,g]
	=
	Q[f,g]
	-
	(p-1)\|(f,g)\|_{\mathcal H}^2,
\end{equation}
where the quadratic form \(Q\) is defined by
\begin{equation}\label{eq:joint-form}
	Q[f,g]
	={}p\int_{\RN}u^2|\nabla f|^2
	+p\iint_{\RN\times\RN}(f(x)-g(y))^2\dd\pi(x,y),
\end{equation}
which constitutes the nonnegative part of \(p\,\widehat{\mathcal C}\).
As before, \eqref{eq:joint-form} is a quadratic form in the single field \((f,g)\); the symmetric bilinear form associated with it arises by polarization, and coincides with the bilinear form \(\langle\mathcal L_r(\,\cdot\,),\,\cdot\,\rangle_{\mathcal H}\) of the operator \eqref{eq:joint-operator} introduced below.

By construction \eqref{eq:Cn-transformed} reads
\begin{equation}\label{eq:Q-and-E}
	Q[f,g]
	=
	p\,\widehat{\mathcal C}[f,g]
	+
	(p-1)\|(f,g)\|_{\mathcal H}^2 .
\end{equation}

We now identify the operators associated with these two forms.  By \eqref{eq:joint-dirichlet-expansion}, the form \(p\widehat{\mathcal C}\) decomposes as the sum of the local form \(p\langle uf,uh\rangle_{H^1_r}\), of the two masses, and of the coupling; each summand is represented on \(\mathcal H\): the local part by \(L\) through \eqref{eq:L-representation}, the two masses by the identity operators, and the coupling by the bounded operators of \eqref{eq:Markov-operators}, whose adjointness, established in Lemma~\ref{lem:Markov-contractions}, is matched by the weight \(p\) in \eqref{eq:joint-space}.  The following lemma makes this precise.

\begin{lemma}\label{lem:joint-operator}
	The form \(p\widehat{\mathcal C}\) is represented on \(\mathcal H\) by the normalized linearized block operator
	\begin{equation}\label{eq:linearized-block-operator}
		\mathcal T
		=
		\begin{pmatrix}
			L+(1-p)I & -R^*\\[1mm]
			-pR & I
		\end{pmatrix},
		\qquad
		D(\mathcal T)=D(L)\oplus L^2_r(\dd\nu),
	\end{equation}
	where adjoints are taken with respect to the weighted inner product \eqref{eq:joint-space}, so that the off-diagonal blocks are adjoints of one another.  By \eqref{eq:Q-and-E}, the operator associated with \(Q\) is
	\begin{equation}\label{eq:joint-operator-intro}
		\mathcal L_r=\mathcal T+(p-1)I,
	\end{equation}
	namely
	\begin{equation}\label{eq:joint-operator}
		\mathcal L_r=
		\begin{pmatrix}
			L&-R^*\\[1mm]
			-pR&pI
		\end{pmatrix},
		\qquad
		D(\mathcal L_r)=D(L)\oplus L^2_r(\dd\nu).
	\end{equation}
	Moreover,
	\begin{equation}\label{eq:joint-kernel-recalled}
		\ker\mathcal L_r=\operatorname{span}\{(1,1)\}.
	\end{equation}
\end{lemma}

\begin{proof}
	In view of  \eqref{eq:E-f-g}, the form \(p\widehat{\mathcal C}\) is the sum of the closed form \(p\langle uf,uh\rangle_{H^1_r}\) and  a bounded form, and hence is itself closed. From \eqref{eq:Q-and-E}, it is bounded from below:
	\[
	p\widehat{\mathcal C}[f,g]\geq-(p-1)\|(f,g)\|_{\mathcal H}^2 .
	\]
	By the representation theorem \cite[Theorem~2.6 of Chapter~VI]{Kato}, the form \(p\widehat{\mathcal C}\) is represented by a unique self-adjoint operator on \(\mathcal H\).
	
	We identify this operator.  Substituting \eqref{eq:joint-dirichlet-expansion} into \eqref{eq:E-f-g} and using \eqref{eq:L-representation}, the form \(p\widehat{\mathcal C}\) reads
	\[
	p\widehat{\mathcal C}[f,g]
	=p\,\langle uf,uf\rangle_{H^1_r}
	+p(1-p)\|f\|_{L^2_r(\dd\mu)}^2
	-2p\langle Rf,g\rangle_{L^2_r(\dd\nu)}
	+\|g\|_{L^2_r(\dd\nu)}^2 ,
	\]
	and polarizing term by term yields its bilinear form
	\begin{equation}\label{eq:joint-form-representation}
		\begin{aligned}
			\bigl\langle\mathcal T(f,g),(h,k)\bigr\rangle_{\mathcal H}
			={}&p\,\langle uf,uh\rangle_{H^1_r}+p(1-p)\langle f,h\rangle_{L^2_r(\dd\mu)}
			-p\langle R^*g,h\rangle_{L^2_r(\dd\mu)}\\
			&-p\langle Rf,k\rangle_{L^2_r(\dd\nu)}
			+\langle g,k\rangle_{L^2_r(\dd\nu)},
		\end{aligned}
	\end{equation}
	defined for every \((f,g),(h,k)\in\mathcal D\oplus L^2_r(\dd\nu)\); here the coupling term is polarized with the aid of Lemma~\ref{lem:Markov-contractions}, which shows that \(R^*\) is the Hilbert-space adjoint of \(R\).  If in addition \(f\in D(L)\), then \(\langle uf,uh\rangle_{H^1_r}=\langle Lf,h\rangle_{L^2_r(\dd\mu)}\) by \eqref{eq:L-representation}, and the right-hand side of \eqref{eq:joint-form-representation} reduces to \(\langle\mathcal T(f,g),(h,k)\rangle_{\mathcal H}\) with \(\mathcal T\) as in \eqref{eq:linearized-block-operator}; hence \(D(L)\oplus L^2_r(\dd\nu)\subseteq D(\mathcal T)\).
	
	Conversely, let \((f,g)\in D(\mathcal T)\), so that there exists \((F,G)\in\mathcal H\) such that the right-hand side of \eqref{eq:joint-form-representation} equals
	\[
	p\langle F,h\rangle_{L^2_r(\dd\mu)}+\langle G,k\rangle_{L^2_r(\dd\nu)}
	\]
	for every \((h,k)\in\mathcal D\oplus L^2_r(\dd\nu)\).  Taking \(h=0\) gives \(G=g-pRf\).  Taking \(k=0\) yields
	\[
	\langle uf,uh\rangle_{H^1_r}=\bigl\langle F-(1-p)f+R^*g,\ h\bigr\rangle_{L^2_r(\dd\mu)}
	\qquad\text{for every }h\in\mathcal D,
	\]
	and since \(F-(1-p)f+R^*g\in L^2_r(\dd\mu)\), the definition \eqref{eq:L-domain} of \(L\) implies that \(f\in D(L)\) and \(Lf=F-(1-p)f+R^*g\).  Thus \(D(\mathcal T)=D(L)\oplus L^2_r(\dd\nu)\), and \(\mathcal T(f,g)=\bigl((L+(1-p)I)f-R^*g,\ g-pRf\bigr)\), which establishes \eqref{eq:linearized-block-operator}.
	
	Finally, by \eqref{eq:Q-and-E}, appending the diagonal form \((p-1)\|(f,g)\|_{\mathcal H}^2\)to the closed form \(p\widehat{\mathcal C}\) shifts the representing self-adjoint operator by \((p-1)I\).  Hence the operator associated with \(Q\) is \[\mathcal L_r=\mathcal T+(p-1)I\]
	 with domain \(D(L)\oplus L^2_r(\dd\nu)\), which proves \eqref{eq:joint-operator-intro} and \eqref{eq:joint-operator}.
	
	It remains to prove \eqref{eq:joint-kernel-recalled}.  Suppose that \((f,g)\in\ker\mathcal L_r\).  Since \(\mathcal L_r\) is the nonnegative self-adjoint operator associated with \(Q\),
	\[
	0
	=
	\langle\mathcal L_r(f,g),(f,g)\rangle_{\mathcal H}
	=
	Q[f,g].
	\]
	Both terms of the quadratic form \eqref{eq:joint-form} are nonnegative.  Hence
	\[
	\int_{\RN}u^2|\nabla f|^2=0,
	\qquad
	\iint_{\RN\times\RN}|f(x)-g(y)|^2\dd\pi(x,y)=0.
	\]
	Since \(u>0\) and \(f\in\mathcal D\subset H^1_{\mathrm{loc},r}(\RN)\), the first equality forces \(f=c\) almost everywhere for some constant \(c\).  By the second identity and \eqref{eq:joint-marginals}, we see that \(g=c\) almost everywhere, and therefore $\ker\mathcal L_r\subseteq\operatorname{span}\{(1,1)\}$.
	 Conversely, \eqref{eq:equal-joint-masses} guarantees that constant functions belong to \(L^2_r(\dd\mu)\) and \(L^2_r(\dd\nu)\).  Moreover, \eqref{eq:calD-realization} and \eqref{eq:L-domain} yield \(1\in\mathcal D\), \(1\in D(L)\), and \(L1=1\), while Lemma~\ref{lem:Markov-contractions} gives \(R1=R^*1=1\).  Hence \(\mathcal L_r(1,1)=0\).  Consequently,
$\ker\mathcal L_r=\operatorname{span}\{(1,1)\}$,
	completing  \eqref{eq:joint-kernel-recalled}.
\end{proof}

We now pass from the block operator \(\mathcal L_r\) to a scalar operator by eliminating its second field.  For every \(\lambda<p\), completing the square in the quadratic form of \(\mathcal L_r-\lambda\) yields a decomposition whose residual form \(s_\lambda\) is the Schur complement of \(\mathcal L_r-\lambda\), see e.g., \cite{Tretter2008}; its spectral consequences at the distinguished value \(\lambda=p-1\) are recorded in Proposition~\ref{prop:Schur-complement} below.

\begin{lemma}\label{lem:Schur}
	For every \(\lambda<p\) and \((f,g)\in\mathcal D\oplus L^2_r(\dd\nu)\),
	\begin{equation}\label{eq:Schur-square}
		Q[f,g]-\lambda\|(f,g)\|_{\mathcal H}^2
		=
		p\,s_\lambda[f]
		+
		(p-\lambda)
		\Big\|g-\frac{p}{p-\lambda}Rf\Big\|_{L^2_r(\dd\nu)}^2,
	\end{equation}
	where the residual form \(s_\lambda\) has the bilinear polarization
	\begin{equation}\label{eq:Schur-form}
		s_\lambda(f,h)
		=
		\langle uf,uh\rangle_{H^1_r}
		-
		\lambda\langle f,h\rangle_{L^2_r(\dd\mu)}
		-
		\frac{p}{p-\lambda}\langle Pf,h\rangle_{L^2_r(\dd\mu)},
		\qquad f,h\in\mathcal D,
	\end{equation}
	and is represented by the self-adjoint operator
	\begin{equation}\label{eq:Schur-operator}
		S_\lambda=L-\lambda I-\frac{p}{p-\lambda}P,
		\quad \text{with} \quad
		D(S_\lambda)=D(L).
	\end{equation}
\end{lemma}

\begin{proof}
	Expanding \(Q\) by \eqref{eq:joint-form} and \eqref{eq:joint-dirichlet-expansion}, and invoking the weighted norm of \eqref{eq:joint-space}, we obtain
	\begin{equation*}\label{eq:Schur-expanded}
		Q[f,g]-\lambda\|(f,g)\|_{\mathcal H}^2
		=
		p\,\langle uf,uf\rangle_{H^1_r}
		-p\lambda\|f\|_{L^2_r(\dd\mu)}^2
		-2p\langle Rf,g\rangle_{L^2_r(\dd\nu)}
		+(p-\lambda)\|g\|_{L^2_r(\dd\nu)}^2 .
	\end{equation*}
	Completing the square in the last two terms, we write
	\[
	-2p\langle Rf,g\rangle_{L^2_r(\dd\nu)}+(p-\lambda)\|g\|_{L^2_r(\dd\nu)}^2
	=
	(p-\lambda)\Big\|g-\frac{p}{p-\lambda}Rf\Big\|_{L^2_r(\dd\nu)}^2
	-
	\frac{p^2}{p-\lambda}\|Rf\|_{L^2_r(\dd\nu)}^2.
	\]
	Collecting the residual terms defines
	\[
	s_\lambda[f]
	:=
	\langle uf,uf\rangle_{H^1_r}
	-\lambda\|f\|_{L^2_r(\dd\mu)}^2
	-\frac{p}{p-\lambda}\|Rf\|_{L^2_r(\dd\nu)}^2.
	\]
	This yields \eqref{eq:Schur-square}. Polarizing term by term and recasting the coupling via \(R^*R=P\) (Lemma~\ref{lem:Markov-contractions}), namely
	\[
	\langle Rf,Rh\rangle_{L^2_r(\dd\nu)}=\langle R^*Rf,h\rangle_{L^2_r(\dd\mu)}=\langle Pf,h\rangle_{L^2_r(\dd\mu)},
	\]
	we arrive at \eqref{eq:Schur-form}. 
	
	We next identify the form domain of \(s_\lambda\).  Write \(s_\lambda=t_0-t_1\), where
	\[
	t_0(f,h)=\langle uf,uh\rangle_{H^1_r},
	\qquad
	t_1(f,h)=\lambda\langle f,h\rangle_{L^2_r(\dd\mu)}+\frac{p}{p-\lambda}\langle Pf,h\rangle_{L^2_r(\dd\mu)},
	\] 
By Lemma~\ref{lem:ground-state-transform}, \(t_0\) is a closed form on \(\mathcal D\); the perturbation \(t_1\) is bounded on \(\mathcal D\) since \(P\) is bounded on \(L^2_r(\dd\mu)\).  The form \(s_\lambda\) is therefore densely defined, closed, symmetric and bounded from below on \(\mathcal D\).
	The representation theorem \cite[Theorem~2.6 of Chapter~VI]{Kato} furnishes a unique self-adjoint operator \(S_\lambda\) on \(L^2_r(\dd\mu)\) such that
	\[
	s_\lambda(f,h)=\langle S_\lambda f,h\rangle_{L^2_r(\dd\mu)}
	\qquad
	\text{for all }f\in D(S_\lambda),\ h\in\mathcal D.
	\]
	By \eqref{eq:L-representation}, the local part \(\langle uf,uh\rangle_{H^1_r}\) is represented by \(L\), and the bounded perturbation \(t_1\) by the bounded self-adjoint operator \(\lambda I+\frac{p}{p-\lambda}P\); hence \(S_\lambda=L-\lambda I-\frac{p}{p-\lambda}P\) with \(D(S_\lambda)=D(L)\), which is \eqref{eq:Schur-operator}. Finally, the square in \eqref{eq:Schur-square} vanishes precisely at \(g=\frac{p}{p-\lambda}Rf\), which is the elimination of the second component defining the Schur complement.
\end{proof}
At \(\lambda=p-1\), \eqref{eq:Schur-operator} yields
\[
S_{p-1}=L-(p-1)I-pP,
\]
which is the scalar operator obtained by eliminating \(g\) from the normalized linearized system \eqref{eq:normalized-linearized-system}.  
Thus the original linearized problem is the spectral problem for \(\mathcal L_r\) at the energy \(\lambda=p-1\): the shift by \((p-1)I\) trades the sign-indefinite form \(\widehat{\mathcal C}\) for the manifestly nonnegative joint form \(Q\).  We now relate the spectrum of \(\mathcal L_r-\lambda\) to that of the scalar operator \(S_\lambda\), and then identify \(S_{p-1}\) with the radial restriction of the original linearization \(L_{+,u}\).

\begin{definition}\label{def:index-nullity}
	Let \(T\) be a self-adjoint operator bounded from below, and let \(t\) be its associated closed quadratic form with form domain \(D(t)\).  Following \cite[Section~2]{BerkolaikoKuchment2022}
	we denote the Morse index of $T$ by  \(i_-(T)\) and its nullity by \(i_0(T)\), namely
	\[
	i_-(T)
	=
	\sup\left\{
	\dim E:
	E\subset D(t),\
	t[f]<0\ \text{for every }0\ne f\in E
	\right\},
	\qquad
	i_0(T)=\dim\ker T .
	\]
	If the negative spectrum of \(T\) is discrete, then \(i_-(T)\) is the number of negative eigenvalues of \(T\), counted with multiplicity.	
\end{definition}

\begin{proposition}\label{prop:Schur-complement}
	For every \(\lambda<p\),
	\begin{equation}\label{eq:Schur-kernel}
		\ker(\mathcal L_r-\lambda)
		=
		\left\{
		\left(f,\frac{p}{p-\lambda}Rf\right):
		f\in\ker S_\lambda
		\right\},
	\end{equation}
	and
	\begin{equation}\label{eq:Schur-inertia}
		i_-(\mathcal L_r-\lambda)
		=
		i_-(S_\lambda),
		\qquad
		i_0(\mathcal L_r-\lambda)
		=
		i_0(S_\lambda).
	\end{equation}
	At \(\lambda=p-1\), the form identity
	\begin{equation}\label{eq:scalar-ground-transform}
		s_{p-1}(f,h)
		=\bigl\langle L_{+,u}(uf),uh\bigr\rangle_{H^{-1}_r,H^1_r}
		\qquad(f,h\in\mathcal D)
	\end{equation}
	entails the equivalence
	\begin{equation}\label{eq:S-p-1-equiv}
		S_{p-1}f=\sigma f
		\quad\Longleftrightarrow\quad
		L_{+,u}(uf)=\sigma V(uf),
	\end{equation}
	and hence
	\begin{equation}\label{eq:index-L+u}
		i_-(S_{p-1})=i_-(L_{+,u}),
		\qquad
		i_0(S_{p-1})=i_0(L_{+,u}).
	\end{equation}
\end{proposition}
\begin{proof}
The kernel formula \eqref{eq:Schur-kernel} and the inertia equalities \eqref{eq:Schur-inertia} follow from Lemma~\ref{lem:Schur-inertia} applied to \eqref{eq:joint-operator} with \(A=L-\lambda I\), \(B=-R^*\) and \(C=(p-\lambda)I\): since \(\lambda<p\), \(C\) is positive definite and boundedly invertible, and the Schur complement of \(C\) in \(\mathcal L_r-\lambda\) is
\[
A-BC^{-1}B^*
=
L-\lambda I-\frac{p}{p-\lambda}R^*R
=
S_\lambda
\]
by \eqref{eq:Schur-operator} and \eqref{eq:P-Markov}.

At \(\lambda=p-1\), \eqref{eq:Schur-operator} reduces to
\[
S_{p-1}=L-(p-1)I-pP.
\]
For \(\phi,\psi\in H_r^1(\RN)\), write \(\phi=uf\) and \(\psi=uh\) with \(f,h\in\mathcal D\) (Lemma~\ref{lem:ground-state-transform}).  The identity \(Au^{p-2}\phi\psi=Au^pfh\) and \eqref{eq:P-Markov} then produce \eqref{eq:scalar-ground-transform}.  Since multiplication by \(u\) is an isometric isomorphism from \(\mathcal D\) onto \(H^1_r(\RN)\) (Lemma~\ref{lem:ground-state-transform}), it identifies the form domains, the kernels and the negative subspaces of \(S_{p-1}\) and \(L_{+,u}\).  Because \(\phi=uf\) matches the two quadratic forms exactly, the spectral problem for \(S_{p-1}\) translates, via \(\phi=uf\), into the pencil problem for \(L_{+,u}\) in the weighted space \(L^2_r(V\dd x)\).  Concretely, \(S_{p-1}f=\sigma f\) is equivalent to \(L_{+,u}\phi=\sigma V\phi\) in \(H^{-1}_r(\RN)\), from which \eqref{eq:S-p-1-equiv} and \eqref{eq:index-L+u} follow.
\end{proof}

We conclude this section by collecting the spectral consequences of the Schur complement analysis and the compactness of the resolvent.

\begin{proposition}\label{prop:joint-spectrum}
	Assume \(p\geq 2\) or \(\alpha\geq 2\).  Set
	\[
	\tau:=
	\begin{cases}
		p,&0<\alpha<N,\ p\geq2,\\[1mm]
		1,&2\leq\alpha<N,\ p<2.
	\end{cases}
	\]
	For every \(\lambda<\tau\), both \(S_\lambda\) and \(\mathcal L_r-\lambda\) are Fredholm.  Every spectral point of \(\mathcal L_r\) below \(\tau\) is an isolated eigenvalue of finite multiplicity.  In particular, zero is a simple isolated eigenvalue of \(\mathcal L_r\), and its positive eigenvalues below \(\tau\) can accumulate only at \(\tau\).
\end{proposition}

\begin{proof}
	Fix \(\lambda<\tau\).  In either case \(\lambda<p\), so \(S_\lambda\) is defined by \eqref{eq:Schur-operator}.  We first suppose   \(p\geq2\). By Lemma~\ref{lem:L-compact}, \(L\) has compact resolvent, while \(P\) is bounded and self-adjoint by Lemma~\ref{lem:Markov-contractions}. Note that the form of \(S_\lambda\) differs from that of \(L\) by a bounded symmetric perturbation and a relatively bounded form perturbation of a form with compact resolvent retains compact resolvent \cite[Theorem~3.4 of Chapter~VI]{Kato}. We conclude that \(S_\lambda\) is self-adjoint with compact resolvent.  Its spectrum therefore consists of isolated eigenvalues of finite multiplicity, and in particular \(S_\lambda\) is Fredholm.
	
	Suppose now that \(p<2\), and set \(\phi=uf\).  By \eqref{eq:joint-measures} and \eqref{eq:P-Markov},
	\begin{equation*}\label{eq:P-pairing-subquadratic}
		\begin{aligned}
			\langle Pf,f\rangle_{L^2_r(\dd\mu)}
			&=
			\int_{\RN}(Pf)f\,\dd\mu
			=
			\int_{\RN}
			\frac{I_\alpha*(u^pf)}{A}\,fAu^p\,\dd x
			=
			\int_{\RN}
			\bigl(I_\alpha*(u^{p-1}\phi)\bigr)
			u^{p-1}\phi\,\dd x.
		\end{aligned}
	\end{equation*}
This, together with	\eqref{eq:joint-measures}, \eqref{eq:Schur-form}, and \(V=Au^{p-2}\), implies
		\begin{align}	\label{eq:subquadratic-Schur-form}
			s_\lambda[f]
			=&
			\langle uf,uf\rangle_{H^1_r}
			-
			\lambda\int_{\RN}f^2\,\dd\mu
			-
			\frac{p}{p-\lambda}
			\langle Pf,f\rangle_{L^2_r(\dd\mu)}
			\nonumber\\
			=&
			\int_{\RN}
			\left(|\nabla\phi|^2+(1-\lambda)\phi^2\right)\dd x
			-
			\lambda\int_{\RN}(V-1)\phi^2\,\dd x\nonumber\\
			&-
			\frac{p}{p-\lambda}\int_{\RN}
			\bigl(I_\alpha*(u^{p-1}\phi)\bigr)
			u^{p-1}\phi\,\dd x.
		\end{align}
	Lemma~\ref{lem:regularity-and-far-field} gives \(V(r)\to1\).  Since \(V\) is positive and continuous, the norms of \(L^2(V\,\dd x)\) and \(L^2(\dd x)\) are equivalent; moreover, since \(\lambda<1\), the first term in the preceding formula is coercive on \(H^1_r(\RN)\).  The perturbation
	\[
	\phi\longmapsto\int_{\RN}(V-1)\phi^2\,\dd x
	\]
	is compact on \(H^1_r(\RN)\), by local Rellich compactness and \(V(r)\to1\).  Moreover, Lemma~\ref{lem:regularity-and-far-field} and \eqref{eq:p-range} provide \(u(r)^{p-1}=O(r^{-(p-1)\beta_p})\) with \((p-1)\beta_p>\alpha\); combined with \(\alpha\geq2\), this implies \((p-1)\beta_p>1+{\alpha/2}\), and hence \(u^{p-1}\in L^{2N/(\alpha+2)}(\RN)\).  Multiplication by \(u^{p-1}\) is therefore compact from \(H^1_r(\RN)\) to \(L^{2N/(N+\alpha)}(\RN)\), and the Hardy--Littlewood--Sobolev inequality ensures that the last form in \eqref{eq:subquadratic-Schur-form} is compact as well. In summary, \(s_\lambda\), expressed in the coordinate \(\phi=uf\), decomposes as the sum of a coercive form on \(H^1_r(\RN)\) and a form that is compact on \(H^1_r(\RN)\). By the representation theorem, a coercive closed form corresponds to a self-adjoint operator with bounded inverse, which is in particular Fredholm.  Since \(S_\lambda\) differs from this operator by a compact perturbation, it is itself Fredholm \cite[Theorem~5.26 of Chapter~IV]{Kato}.  The norm equivalence established above transports this conclusion back to \(L^2_r(\dd\mu)\).
	
Since \(\lambda<p\), we apply Lemma~\ref{lem:Schur-inertia}(iii) to \eqref{eq:joint-operator} with \(A=L-\lambda I\), \(B=-R^*\) and \(C=(p-\lambda)I\).  Adjoints are taken in the weighted inner product \eqref{eq:joint-space}, under which the lower-left block \(-pR\) of \eqref{eq:joint-operator} is the adjoint of \(-R^*\).  Since \(S_\lambda\) is Fredholm, the lemma implies that \(\mathcal L_r-\lambda\) is also Fredholm.
	
	Now let \(\lambda<\tau\) be a spectral point of \(\mathcal L_r\).  Since \(\mathcal L_r\) is self-adjoint and \(\mathcal L_r-\lambda\) is Fredholm, and since \(\lambda\in\sigma(\mathcal L_r)\), the standard description of the discrete spectrum of a self-adjoint operator applies: for a self-adjoint operator \(T\), every point of \(\sigma(T)\) at which \(T-\lambda\) is Fredholm is an isolated eigenvalue of finite multiplicity \cite[Section~2.1]{BerkolaikoKuchment2022}.  Hence every spectral point of \(\mathcal L_r\) below \(\tau\) is an isolated eigenvalue of finite multiplicity; in particular, by \eqref{eq:joint-kernel-recalled}, zero is a simple isolated eigenvalue.
	
	If a sequence of distinct eigenvalues in \((0,\tau)\) converged to some \(\lambda_0<\tau\), then the closedness of the spectrum would imply \(\lambda_0\in\sigma(\mathcal L_r)\), contradicting the isolation of \(\lambda_0\).  Thus eigenvalues in \((0,\tau)\) can accumulate only at \(\tau\).
\end{proof}

\section{Monotone zero modes of the Schur complement}
\label{sec:monotone-zero-mode}

The goal of this section is  to establish  Proposition~\ref{prop:monotone-zero-mode}
below: if
\(\sigma(\mathcal L_r)\) has the second eigenvalue $\lambda$ in \((0,p-1]\), then a nonzero radial nonincreasing zero mode of \(S_\lambda\) exists.  By \eqref{eq:Schur-kernel}, this zero mode is equivalently an eigenfunction of \(\mathcal L_r\) with eigenvalue \(\lambda\), and its monotonicity is an ingredient of
the spectral gap argument of Section~\ref{sec:spectral-gap}.  Its
existence follows from the heat-flow limit of \(-S_\lambda\) started at a radial nonincreasing function; the limit inherits radial monotonicity because \(e^{-sS_\lambda}\) preserves this property, a fact established in Proposition~\ref{prop:Schur-monotone} by combining, via the Trotter product formula, the monotonicity-preserving properties \(P\) and of \(e^{-sL}\).

Throughout this section we assume all
functions and operators are radial and restrict \((\alpha,p)\) to the variational range \eqref{eq:p-range} with $\alpha\ge2$; that is, $2\le\alpha<N$ and $\frac{N+\alpha}{N}<p<\frac{N+\alpha}{N-2}$.
This range is narrower than the condition $p\ge2$ or $\alpha\ge2$ of Proposition~\ref{prop:joint-spectrum}: it excludes $\alpha<2$, even when $p\ge2$.
The additional restriction \(\alpha\ge2\) is precisely what the radial monotonicity of \(P\) requires: Lemma~\ref{lem:TP2-MLR} below rests on the log-concavity of the kernel \(k\) (Lemma~\ref{lem:log-concavity}).  This is the only case needed for the spectral gap argument of Section~\ref{sec:spectral-gap}.

\subsection{The coupling operator $P$ preserves radial monotonicity}

To study the coupling operator \(P\), we first express the kernel of the Riesz potential \(I_\alpha\) in polar coordinates.
For \(r,s>0\), set
\[
K_\alpha(r,s)
:=\int_{\mathbb S^{N-1}}I_\alpha(re_1-s\omega)\,\dd\omega.
\]
Writing \(y=s\omega\) in polar coordinates, we obtain, for every radial
\(f\),
\begin{equation}\label{eq:radial-convolution}
	\begin{aligned}
		(I_\alpha*(u^pf))(re_1)
		&=\int_{\RN}I_\alpha(re_1-y)u(y)^pf(y)\,\dd y =\int_0^\infty K_\alpha(r,s)u(s)^pf(s)s^{N-1}\,\dd s.
	\end{aligned}
\end{equation}

Let
\[
\beta=\frac{N-\alpha}{2},
\qquad
\delta=\frac{\alpha-2}{2}.
\]
For \(t\in\mathbb R\), define
\[
k(t)
:=e^{-\beta|t|}
\int_0^1
\theta^{\beta-1}(1-\theta)^\delta
(1-e^{-2|t|}\theta)^\delta\,\dd\theta.
\]
This integral is finite and positive when $1<\alpha<N$.
By Euler's integral representation of the Gauss hypergeometric function \(\cF\) (see \cite[Eq.~9.111]{GradshteynRyzhik}), with
\(B\) denoting the beta function, we have
\[
\int_0^1\theta^{\beta-1}(1-\theta)^\delta(1-z\theta)^\delta\,\dd\theta
=B(\beta,\delta+1)\,
\cF\!\left(-\delta,\beta;\beta+\delta+1;z\right),
\qquad 0<z\leq1,
\]
and hence
\[
k(t)
=\mathrm e^{-\beta|t|}\,B(\beta,\delta+1)\,
\cF\!\left(-\delta,\beta;\beta+\delta+1;\mathrm e^{-2|t|}\right),
\qquad t\in\mathbb R.
\]
In particular \(\cF\equiv1\) when \(\delta=0\), and  
\(\cF(a,b;c;0)=1\) yields that 
\(\mathrm e^{\beta|t|}k(t)\to B(\beta,\delta+1)\) as
\(|t|\to\infty\).

\begin{lemma}\label{lem:log-conjugation}
	Assume \(1<\alpha<N\).  Then 
	\[
	K_\alpha(e^t,e^\sigma)
	=\frac{2^{1-\alpha}}{\Gamma(\alpha/2)^2}e^{-\beta(t+\sigma)}k(t-\sigma),
	\qquad t,\sigma\in\mathbb R.
	\]
\end{lemma}

\begin{proof}
Since
	\[
	|re_1-s\omega|^2=r^2+s^2-2rs\,e_1\cdot\omega,
	\]
and	the integral defining \(K_\alpha(r,s)\) is symmetric in $r$ and $s$, so 
	\(
	K_\alpha(r,s)=K_\alpha(s,r).
	\)
	Therefore it suffices to prove the identity for \(t\ge \sigma\). Set
	\[
	r=\mathrm e^t,\qquad s=\mathrm e^\sigma,\qquad \varsigma =\mathrm e^{-(t-\sigma)}\in(0,1].
	\]
	Then
	\(
	|\mathrm e^t e_1-\mathrm e^\sigma\omega|
	=\mathrm e^t|e_1-\varsigma\omega|,
	\)
	and hence
	\[
	|\mathrm e^t e_1-\mathrm e^\sigma\omega|^{-(N-\alpha)}
	=\mathrm e^{-2\beta t}
	\bigl(1+\varsigma^2-2\varsigma\,e_1\cdot\omega\bigr)^{-\beta}.
	\]
	Since
	\[
	\mathrm e^{-2\beta t}
	=\mathrm e^{-\beta(t+\sigma)}\mathrm e^{-\beta(t-\sigma)},
	\]
	we obtain
	\[
	K_\alpha(\mathrm e^t,\mathrm e^\sigma)
	=c_{N,\alpha}\mathrm e^{-\beta(t+\sigma)}
	\mathrm e^{-\beta(t-\sigma)}
	\int_{\mathbb S^{N-1}}
	\bigl(1+\varsigma^2-2\varsigma\,e_1\cdot\omega\bigr)^{-\beta}
	\,\mathrm d\omega.
	\]
	
We now evaluate the spherical integral. Write \(e_1\cdot\omega=\cos\theta\), \(\theta\in[0,\pi]\), so that
	\[
	\mathrm d\omega=\sin^{N-2}\theta\,\mathrm d\theta\,\mathrm d\Omega_{N-2}.
	\]
	Integrating out the \(N-2\) angular variables, we deduce
	\begin{align*}
		\int_{\mathbb S^{N-1}}
		\bigl(1+\varsigma^2-2\varsigma\,e_1\cdot\omega\bigr)^{-\beta}
		\,\mathrm d\omega
		=&|\mathbb S^{N-2}|
		\int_0^\pi \sin^{N-2}\theta\,
		\bigl(1+\varsigma^2-2\varsigma\cos\theta\bigr)^{-\beta}\,\mathrm d\theta\\
		=&|\mathbb S^{N-2}|\,B\!\left(\frac{N-1}{2},\frac12\right)
		\cF\!\left(\beta,\beta-\frac{N-2}{2};\frac{N}{2};\varsigma^2\right),
	\end{align*}
	where the last identity follows from \cite[Eq.~3.665.2]{GradshteynRyzhik}.
	Note that
	\[
	\beta-\frac{N-2}{2}
	=\frac{N-\alpha}{2}-\frac{N-2}{2}
	=\frac{2-\alpha}{2}
	=-\delta.
	\]
	We thus conclude
	\[
	\int_{\mathbb S^{N-1}}
	\bigl(1+\varsigma^2-2\varsigma\,e_1\cdot\omega\bigr)^{-\beta}
	\,\mathrm d\omega
	=|\mathbb S^{N-2}|\,B\!\left(\frac{N-1}{2},\frac12\right)
	\cF\!\left(\beta,-\delta;\frac{N}{2};\varsigma^2\right).
	\]
	Using  the identity
	\(
	\frac{N}{2}=\beta+\delta+1
	\) and the symmetry \(F(a,b;c;z)=F(b,a;c;z)\) (see e.g., \cite[Eq.~9.100]{GradshteynRyzhik}), we may write
	\[
	\int_{\mathbb S^{N-1}}
	\bigl(1+\varsigma^2-2\varsigma\,e_1\cdot\omega\bigr)^{-\beta}
	\,\mathrm d\omega
	=|\mathbb S^{N-2}|\,B\!\left(\frac{N-1}{2},\frac12\right)
	\cF\!\left(-\delta,\beta;\beta+\delta+1;\varsigma^2\right).
	\]
	
	Since \(t\ge\sigma\), we have \(t-\sigma=|t-\sigma|\), and the form of
	\(k\) recorded above implies
	\[
	\mathrm e^{-\beta(t-\sigma)}
	\cF\!\left(-\delta,\beta;\beta+\delta+1;\mathrm e^{-2(t-\sigma)}\right)
	=\frac{k(t-\sigma)}{B(\beta,\delta+1)} .
	\]
	Substituting the spherical average computed above, we obtain
	\[
	K_\alpha(\mathrm e^t,\mathrm e^\sigma)
	=C_{N,\alpha}\,
	\mathrm e^{-\beta(t+\sigma)}
	k(t-\sigma),
	\]
	where
	\[
	C_{N,\alpha}
	=c_{N,\alpha}|\mathbb S^{N-2}|\,B\!\left(\frac{N-1}{2},\frac12\right)
	\frac{\Gamma(\beta+\delta+1)}{\Gamma(\beta)\Gamma(\delta+1)}
	=\frac{2^{1-\alpha}}{\Gamma(\alpha/2)^2}.
	\]
	By symmetry, the same identity holds when \(t<\sigma\). This completes the proof.
\end{proof}

\begin{lemma} \label{lem:log-concavity}
 For \(\alpha=2\),  \(k(t)=\frac2{N-2}
	e^{-(N-2)|t|/2}\); for $\alpha\in(2,N)$,  \(\log k(t)\) is strictly concave on
	\(\R\).
\end{lemma}

\begin{proof}
	Positivity follows directly from the definition.  If \(\alpha=2\),
	then \(\delta=0\), and
	\[
	k(t)=\frac1\beta e^{-\beta|t|}
	=\frac2{N-2}e^{-(N-2)|t|/2}.
	\]
	Thus the assertion holds in this case.
	
	Assume now that \(2<\alpha<N\), so that \(\beta,\delta>0\), and
	define
	\[
	\varphi(s)=
	\begin{cases}
		e^{-\beta s}(1-e^{-2s})^\delta,&s>0,\\
		0,&s\leq0.
	\end{cases}
	\]
	The function \(\varphi\) is integrable with support \((0,\infty)\).  Moreover, for \(s>0\),
	\[
	\frac{\dd^2}{\dd s^2}\log\varphi(s)
	=-\frac{4\delta e^{2s}}{(e^{2s}-1)^2}<0.
	\]
	Hence \(\varphi\) is log-concave on \(\mathbb R\).
	For \(t\geq0\), the change of variables
	\(\theta=e^{-2s}\) gives
	\[
	\begin{aligned}
		2\int_{\mathbb R}\varphi(s)\varphi(s+t)\,\dd s
		&=2e^{-\beta t}\int_0^\infty
		e^{-2\beta s}(1-e^{-2s})^\delta
		(1-e^{-2(s+t)})^\delta\,\dd s\\
		&=e^{-\beta t}\int_0^1
		\theta^{\beta-1}(1-\theta)^\delta
		(1-e^{-2t}\theta)^\delta\,\dd\theta =k(t).
	\end{aligned}
	\]
	Replacing \(t\) by \(-t\) and substituting \(s\mapsto s-t\) in the
	shifted integral, we get
	\[
	2\int_{\mathbb R}\varphi(s)\varphi(s-t)\,\dd s
	=2\int_{\mathbb R}\varphi(s+t)\varphi(s)\,\dd s
	=2\int_{\mathbb R}\varphi(s)\varphi(s+t)\,\dd s ,
	\]
	so the integral is even in \(t\).  Hence
	\[
	k(t)=2\int_{\mathbb R}\varphi(s)\varphi(s+t)\,\dd s,
	\qquad t\in\mathbb R.
	\]
	Its reflection \(s\mapsto\varphi(-s)\) is also log-concave.
	Moreover, substituting \(s\mapsto-s\), we have
	\[
	\int_{\mathbb R}\varphi(-s)\varphi(t-s)\,\dd s
	=\int_{\mathbb R}\varphi(s)\varphi(t+s)\,\dd s
	=\frac12k(t).
	\]
	Thus \(k/2=\varphi(-\cdot)*\varphi\) is the convolution of two
	log-concave functions, where \((f*g)(t)=\int_{\mathbb R}f(s)g(t-s)\,\dd s\)
	denotes the usual convolution on \(\mathbb R\).
	By Pr\'ekopa's convolution theorem
	\cite[Theorem~7]{Prekopa1973}, \(k\) is log-concave on
	\(\mathbb R\).
	
	It remains to prove the strict concavity.
	The defining integral shows that \(\log k\) is real analytic on
	\(\R\setminus \{0\}\).  If \(\log k\) were not strictly concave, equality in
	the concavity inequality at two distinct points would imply that
	\(\log k\) is affine on the interval between them.
	Assume without loss of generality that this interval contains a subinterval of \((0,\infty)\).
	Analyticity would then imply
	\[
	k(t)=Ce^{at},
	\qquad t>0,
	\]
	for some constants \(C>0\) and \(a\in\mathbb R\).
	On the other hand, by the hypergeometric representation of \(k\)
	proved above, \(\mathrm e^{\beta t}k(t)\to B(\beta,\delta+1)\in(0,\infty)\)
	as \(t\to+\infty\).  Combining this with
	\(k(t)=Ce^{at}\) for \(t>0\), we obtain
	\[
	\mathrm e^{\beta t}k(t)=Ce^{(\beta+a)t},\qquad t>0 .
	\]
	Since the left-hand side has a finite positive limit as
	\(t\to+\infty\), we must have \(a=-\beta\);  the limit
	then equals \(C\), so \(C=B(\beta,\delta+1)\).  Consequently,
	\(\mathrm e^{\beta t}k(t)\equiv B(\beta,\delta+1)\) for every
	\(t>0\).
	
	However, for \(t>0\),
	\[
	e^{\beta t}k(t)
	=\int_0^1\theta^{\beta-1}(1-\theta)^\delta
	(1-e^{-2t}\theta)^\delta\,\dd\theta
	\]
	is strictly increasing in \(t\) since \(\delta>0\) in the present regime
	\(2<\alpha<N\).  This contradicts the
	constancy of \(e^{\beta t}k(t)\) proved above.  Hence \(\log k\)
	is strictly concave.
\end{proof}

\begin{lemma}
	\label{lem:TP2-MLR}
	For \(2\leq\alpha<N\), whenever \(r_1<r_2\) and \(s_1<s_2\), it holds
	\[
	K_\alpha(r_1,s_1)K_\alpha(r_2,s_2)
	\geq
	K_\alpha(r_1,s_2)K_\alpha(r_2,s_1).
	\]
	Consequently, \(P\) maps radial
	nonincreasing functions in \(L^2_r(\dd\mu)\) to radial nonincreasing
	functions.
\end{lemma}

\begin{proof}
	Let \(t_1<t_2\) and \(\sigma_1<\sigma_2\).  Since
	\(\log k\) is concave, its
	derivative is nonincreasing.  Hence
	\[
	\begin{aligned}
		\log k(t_1-\sigma_1)-&\log k(t_1-\sigma_2)
		=\int_0^{\sigma_2-\sigma_1}
		(\log k)'(t_1-\sigma_2+v )\,\dd v\\
		&\geq
		\int_0^{\sigma_2-\sigma_1}
		(\log k)'(t_2-\sigma_2+v)\,\dd v =\log k(t_2-\sigma_1)-\log k(t_2-\sigma_2).
	\end{aligned}
	\]
	Exponentiating gives
	\begin{equation}\label{eq31}
		k(t_1-\sigma_1)k(t_2-\sigma_2)
		\geq
		k(t_1-\sigma_2)k(t_2-\sigma_1).
	\end{equation}
	Now set
	$r_i=e^{t_i}$, $s_i=e^{\sigma_i}$,
	for $i=1,2$.
	Then \(r_1<r_2\) and \(s_1<s_2\).  By
	Lemma~\ref{lem:log-conjugation} and \eqref{eq31},
	\begin{equation}\label{eq33}
		\begin{aligned}
			&K_\alpha(r_1,s_1)K_\alpha(r_2,s_2)
			-K_\alpha(r_1,s_2)K_\alpha(r_2,s_1)\\
			&=C_{N,\alpha}^2
			e^{-\beta(t_1+t_2+\sigma_1+\sigma_2)}
			\Bigl[
			k(t_1-\sigma_1)k(t_2-\sigma_2)
			-k(t_1-\sigma_2)k(t_2-\sigma_1)
			\Bigr]\geq0.
		\end{aligned}
	\end{equation}
	
	We next prove that \(P\) preserves radial monotonicity.  Let
	\(f\in L^2_r(\dd\mu)\cap L^\infty(\RN)\) be a nonincreasing radial function.  By \eqref{eq:A-def}, \eqref{eq:P-Markov} and \eqref{eq:radial-convolution}
	\[
	A(r)
	=\int_0^\infty
	K_\alpha(r,s)u(s)^ps^{N-1}\,\dd s,\qquad
	Pf(r)
	=\frac1{A(r)}
	\int_0^\infty
	K_\alpha(r,s)f(s)u(s)^ps^{N-1}\,\dd s.
	\]
	Hence, for \(0<r_1<r_2\), setting $C_1^{-1}=2A(r_1)A(r_2)$, we obtain
	\[
	\begin{aligned}
		&Pf(r_2)-Pf(r_1)\\
		&=\frac1{A(r_2)}
		\int_0^\infty
		K_\alpha(r_2,s)f(s)u(s)^ps^{N-1}\,\dd s -\frac1{A(r_1)}
		\int_0^\infty
		K_\alpha(r_1,s)f(s)u(s)^ps^{N-1}\,\dd s\\
		&=C_1
		\int_0^\infty\int_0^\infty
		f(s)\bigl[
		K_\alpha(r_2,s)K_\alpha(r_1,t)
		-K_\alpha(r_1,s)K_\alpha(r_2,t)
		\bigr]
		u(s)^pu(t)^ps^{N-1}t^{N-1}\,\dd s\dd t\\
		&=C_1
		\int_0^\infty\int_0^\infty
		\bigl(f(s)-f(t)\bigr)
		\bigl[
		K_\alpha(r_2,s)K_\alpha(r_1,t)
		-K_\alpha(r_1,s)K_\alpha(r_2,t)
		\bigr]
		u(s)^pu(t)^ps^{N-1}t^{N-1}\,\dd s\dd t.
	\end{aligned}
	\]
 The last equality follows by interchanging \(s\) and \(t\) in the
	preceding double integral and then averaging the two expressions.
	When \(s<t\), we have
	$
	f(s)-f(t)\geq0.
	$
	This, together with \eqref{eq33}, implies
	$
	K_\alpha(r_2,s)K_\alpha(r_1,t)
	-K_\alpha(r_1,s)K_\alpha(r_2,t)\leq0.
	$
	When \(s>t\), both inequalities are reversed.  The integrand in
	the last double integral is  nonpositive for all \(s,t>0\), and
	hence
	$
	Pf(r_2)-Pf(r_1)\leq0,
	$
	which shows that \(Pf\) is nonincreasing.
	
	Finally, let \(f\in L^2_r(\dd\mu)\) admit a nonincreasing
	representative and set
	$
	f_\ell=\max\{-\ell,\min\{f,\ell\}\}.
	$
	Then \(f_\ell\to f\) in \(L^2_r(\dd\mu)\).  Since \(P\) is bounded on
	this space,
	\[
	Pf_\ell\to Pf
	\qquad\text{in }L^2_r(\dd\mu).
	\]
	Each \(Pf_\ell\) is nonincreasing by the bounded case established above.  Since \(A(r)u(r)^pr^{N-1}\) is strictly positive on
	\((0,\infty)\), a subsequence converges to \(Pf\) on a set
	\(E\subset(0,\infty)\) of full Lebesgue measure.  For
	\(r_1,r_2\in E\) with \(r_1<r_2\), passing to the limit gives
	\(Pf(r_1)\geq Pf(r_2)\).  Thus \(Pf|_E\) has a nonincreasing
	representative on \((0,\infty)\).
\end{proof}

\subsection{The local semigroup preserves radial monotonicity}

Set $r=e^t$. For radial $f,g\in L^2_r(\dd\mu)$, write $F(t)=f(e^t)$ and $G(t)=g(e^t)$. Define the weights
\[
m(t)=|\mathbb S^{N-1}|e^{Nt}A(e^t)u(e^t)^p,
\qquad
w(t)=|\mathbb S^{N-1}|e^{(N-2)t}u(e^t)^2.
\]
Both belong to $L^1(\mathbb R)$: indeed,
\[
\int_{\mathbb R}m(t)\,\dd t
=\int_{\RN}Au^p\,\dd x
=\mu(\RN)<\infty,
\]
while
\[
\int_{\mathbb R}w(t)\,\dd t
=|\mathbb S^{N-1}|\int_0^\infty r^{N-3}u(r)^2\,\dd r<\infty,
\]
since the integral over $(0,1)$ is finite ($u$ is bounded and $N\geq3$) and the integral over $(1,\infty)$ is bounded by $\|u\|_{L^2(\RN)}^2$. The change of variables
\[
\int_{\RN}fg\,\dd\mu
=\int_{\mathbb R}FG\,m\,\dd t
\]
shows that the map
\begin{equation}\label{uoper}
	(Uf)(t):=f(e^t)
\end{equation}
defines a unitary operator $U\colon L^2_r(\dd\mu)\to L^2(m\,\dd t)$, with $(U^{-1}F)(r)=F(\log r)$. For $f,g\in\mathcal D$,
\begin{equation}\label{eq:log-coordinate-form}
	\int_{\RN}u^2\nabla f\cdot\nabla g\,\dd x
	=\int_{\mathbb R}wF'G'\,\dd t.
\end{equation}

The derivative in \eqref{eq:log-coordinate-form} motivates the flux operator
\[
B_0\colon C_c^\infty(\mathbb R)\subset L^2(w^{-1}\,\dd t)
\longrightarrow L^2(m\,\dd t),
\qquad
B_0Y=\frac{Y'}{m}.
\]
We show that $B_0$ is closable; its closure $B$ and adjoint $B^*$ are the main objects of this subsection.

\begin{lemma}\label{lem:flux-operator}
	The operator $B_0$ is closable. Its closure $B$ satisfies $B^*=B_0^*$ with
	\begin{equation}\label{eq:flux-domain}
		\begin{aligned}
			D(B^*)
			&=\left\{
			F\in L^2(m\,\dd t)\cap H^1_{\mathrm{loc}}(\mathbb R):
			\int_{\mathbb R}w|F'|^2\,\dd t<\infty
			\right\},\\
			B^*F&=-wF'.
		\end{aligned}
	\end{equation}
	Moreover, $U\mathcal D=D(B^*)$, and for $F=Uf$, $G=Ug$,
	\begin{equation}\label{eq:log-flux-form}
		\int_{\RN}u^2\nabla f\cdot\nabla g\,\dd x
		=
		\int_{\mathbb R}wF'G'\,\dd t
		=
		\langle B^*F,B^*G\rangle_{L^2(w^{-1}\,\dd t)}.
	\end{equation}
\end{lemma}

\begin{proof} 
  Suppose $Y_n\to0$ in $L^2(w^{-1}\,\dd t)$ and $B_0Y_n=Y_n'/m\to Z$ in $L^2(m\,\dd t)$. For every $F\in C_c^\infty(\mathbb R)$, integration by parts gives
	\[
	\langle Z,F\rangle_{L^2(m\,\dd t)}
	=\lim_{n\to\infty}\int_{\mathbb R}Y_n'F\,\dd t
	=-\lim_{n\to\infty}\int_{\mathbb R}Y_nF'\,\dd t
	=0,
	\]
	since $wF'\in L^2(w^{-1}\,\dd t)$ and $Y_n\to0$ in that space. Density of $C_c^\infty(\mathbb R)$ in $L^2(m\,\dd t)$ then yields $Z=0$. Hence $B_0$ is closable, and $B^*=B_0^*$ by \cite[Problem~5.24]{Kato}.
	
	 A function $F\in L^2(m\,\dd t)$ belongs to $D(B^*)$ if and only if there exists $H\in L^2(w^{-1}\,\dd t)$ such that $\langle B_0Y,F\rangle_{L^2(m\,\dd t)}=\langle Y,H\rangle_{L^2(w^{-1}\,\dd t)}$ for all $Y\in C_c^\infty(\mathbb R)$. Since $B_0Y=Y'/m$, this reads $\int_{\mathbb R}Y'F\,\dd t=\int_{\mathbb R}YH/w\,\dd t$, hence $F'=-H/w$ distributionally. As $m$ and $w$ are strictly positive and continuous, $F\in H^1_{\mathrm{loc}}(\mathbb R)$, $H=-wF'$, and $\int_{\mathbb R}w|F'|^2\,\dd t<\infty$. The converse follows by reversing the argument, giving \eqref{eq:flux-domain}.
	
	 By \eqref{eq:calD-realization} and \eqref{eq:flux-domain}, $U\mathcal D=D(B^*)$. The identity \eqref{eq:log-flux-form} then follows from \eqref{eq:log-coordinate-form} and the definition of $B^*$.
\end{proof}

Let $L$ be the nonnegative self-adjoint operator from \eqref{eq:L-representation}. Define the conjugated operator via the unitary map $U$ in \eqref{uoper}:
\[
\widetilde L=ULU^{-1},
\qquad
D(\widetilde L)=U\bigl(D(L)\bigr).
\]
The following lemmas collect the properties of $B$ needed to prove that $e^{-s\widetilde L}$ (and hence $e^{-sL}$) preserves radial monotonicity.

\begin{lemma}\label{lem:factorization}
	The conjugated operator admits the factorization
	\begin{equation}\label{eq:flux-factorization}
		\widetilde L=BB^*+I.
	\end{equation}
\end{lemma}

\begin{proof} 
	Let $F\in D(\widetilde L)$. Since $D(\widetilde L)=U(D(L))\subseteq U\mathcal D=D(B^*)$, we have $F\in D(B^*)$. For every $G\in D(B^*)$, set $g=U^{-1}G$. The definitions of $\widetilde L$ and $L$, \eqref{eq:ground-state-transform}, and \eqref{eq:log-flux-form} give
	\[
	\langle\widetilde L F,G\rangle_{L^2(m\,\dd t)}
	=\langle Lf,g\rangle_{L^2_r(\dd\mu)}
	=\langle B^*F,B^*G\rangle_{L^2(w^{-1}\,\dd t)}+\langle F,G\rangle_{L^2(m\,\dd t)}.
	\]
	Since $B$ is closed, $B=(B^*)^*$, so this identifies $B^*F\in D(B)$ with $BB^*F=\widetilde L F-F$. Thus $F\in D(BB^*)$ and $BB^*F=\widetilde L F-F$. Conversely, let $F\in D(BB^*)$. Then $F\in D(B^*)$, and for every $G\in D(B^*)$, 
	$$\langle B^*F,B^*G\rangle_{L^2(w^{-1}\,\dd t)}=\langle BB^*F,G\rangle_{L^2(m\,\dd t)}.$$
	 Setting $f=U^{-1}F$ and $G=Ug$ for $g\in\mathcal D$, the ground-state identity \eqref{eq:ground-state-transform}, \eqref{eq:log-flux-form}, and the unitarity of $U$ give
	\[
	\langle uf,ug\rangle_{H^1_r}
	=\langle B^*F,B^*G\rangle_{L^2(w^{-1}\,\dd t)}+\langle F,G\rangle_{L^2(m\,\dd t)}
	=\langle(BB^*+I)F,G\rangle_{L^2(m\,\dd t)}.
	\]
	By the definition of $L$, this means $f\in D(L)$ and $Lf=U^{-1}(BB^*+I)F$. Therefore $F\in D(\widetilde L)$ and $\widetilde L F=(BB^*+I)F$. The two inclusions together complete \eqref{eq:flux-factorization}.
\end{proof}

\begin{lemma}\label{lem:positivity}
	The semigroup $e^{-sB^*B}$ preserves nonnegative functions:
	\[
	Y\geq0
	\quad\Longrightarrow\quad
	e^{-sB^*B}Y\geq0,
	\qquad s\geq0.
	\]
\end{lemma}

\begin{proof} 
	Since $B$ is densely defined and closed, $B^*B$ is self-adjoint and nonnegative, associated with the closed form $\langle BY,BZ\rangle_{L^2(m\,\dd t)}$ for $Y,Z\in D(B)$ \cite[Proposition~1.31]{Ouhabaz2005}. As $B=\overline{B_0}$, its domain is the completion of $C_c^\infty(\mathbb R)$ under the graph norm
	\[
	\|Y\|_{D(B)}^2
	=
	\int_{\mathbb R}\frac{|Y|^2}{w}\,\dd t
	+
	\int_{\mathbb R}\frac{|Y'|^2}{m}\,\dd t.
	\]
	Since $|Y|$ does not increase this norm---indeed, $\||Y|'\|\leq\|Y'\|$ pointwise a.e.---we have $|Y|\in D(B)$ and $\|B|Y|\|_{L^2(m\,\dd t)}\leq\|BY\|_{L^2(m\,\dd t)}$ whenever $Y\in D(B)$. The first Beurling--Deny criterion \cite[Corollary~2.18]{Ouhabaz2005} then implies that $e^{-sB^*B}$ preserves nonnegativity.
\end{proof}

\begin{lemma}\label{lem:L-monotone}
	If $f\in L^2_r(\dd\mu)$ is nonincreasing, then $e^{-sL}f$ is nonincreasing for every $s\geq0$.
\end{lemma}

\begin{proof} 
	Since $\widetilde L=BB^*+I$, the semigroup satisfies $e^{-s\widetilde L}=e^{-s}e^{-sBB^*}$. Lemma~\ref{lem:product-intertwining} yields the intertwining relation
	\begin{equation*}\label{eq:flux-intertwining}
		B^*e^{-sBB^*}
		=
		e^{-sB^*B}B^*
		\qquad\text{on }D(B^*),\quad s\geq0.
	\end{equation*}
	
	\textit{Step 1: $F\in D(B^*)$.} Let $F\in D(B^*)$ be nonincreasing, so $B^*F=-wF'\geq0$ by \eqref{eq:flux-domain}. The intertwining relation and Lemma~\ref{lem:positivity} give
	\[
	B^*e^{-s\widetilde L}F
	=e^{-s}e^{-sB^*B}B^*F
	\geq0.
	\]
	Since $B^*e^{-s\widetilde L}F=-w(e^{-s\widetilde L}F)'$ and $w>0$, we conclude $(e^{-s\widetilde L}F)'\leq0$ a.e., so $e^{-s\widetilde L}F$ is nonincreasing.
	
	\textit{Step 2: Extension to $L^2(m\,\dd t)$.} Let $F\in L^2(m\,\dd t)$ be nonincreasing. Set $F_\ell=\max\{-\ell,\min\{F,\ell\}\}$ and $F_{\ell,\varepsilon}=\eta_\varepsilon*F_\ell$, where $\eta_\varepsilon$ is a standard nonnegative mollifier. Then $F_{\ell,\varepsilon}$ is smooth, bounded, nonincreasing, with bounded derivative. Since $m,w\in L^1(\mathbb R)$, \eqref{eq:flux-domain} gives $F_{\ell,\varepsilon}\in D(B^*)$, and Step~1 implies $e^{-s\widetilde L}F_{\ell,\varepsilon}$ is nonincreasing.
	
	Since $\widetilde L-I=U(L-I)U^{-1}$ is nonnegative and self-adjoint, the spectral theorem gives $\|e^{-s\widetilde L}\|\leq1$. By dominated convergence ($|F_{\ell,\varepsilon}|\leq\ell$, $m\in L^1$), $F_{\ell,\varepsilon}\to F_\ell$ in $L^2(m\,\dd t)$ as $\varepsilon\downarrow0$, and $F_\ell\to F$ as $\ell\to\infty$. The contraction property then yields
	\[
	e^{-s\widetilde L}F_{\ell,\varepsilon}\to e^{-s\widetilde L}F_\ell\to e^{-s\widetilde L}F
	\qquad\text{in }L^2(m\,\dd t).
	\]
	Passing to a.e.-convergent subsequences, $e^{-s\widetilde L}F$ has a nonincreasing representative.
	
	\textit{Step 3: Transfer to $L$.} Since $\widetilde L=ULU^{-1}$, Lemma~\ref{lem:unitary-calculus} gives $e^{-s\widetilde L}=Ue^{-sL}U^{-1}$. As $t\mapsto e^t$ is increasing, $Uf$ is nonincreasing if and only if $f$ is. Hence $e^{-sL}f$ is nonincreasing whenever $f\in L^2_r(\dd\mu)$ is.
\end{proof}

We conclude with the main result of this subsection.
\begin{proposition}\label{prop:Schur-monotone}
	For every $\lambda<p$ and $s\geq0$, the semigroup $e^{-sS_\lambda}$ preserves the closed cone of radial nonincreasing functions in $L^2_r(\dd\mu)$.
\end{proposition}

\begin{proof}
	By Lemma~\ref{lem:TP2-MLR}, $P$ preserves radial monotonicity: if $f$ is radial nonincreasing, so is $P^j f$ for every $j\geq0$. Since all coefficients in
	\[
	e^{cP}f
	=\sum_{j=0}^\infty\frac{c^j}{j!}P^jf,
	\qquad c\geq0,
	\]
	are nonnegative, every partial sum is radial nonincreasing. The series converges in $L^2_r(\dd\mu)$, so $e^{cP}f$ is radial nonincreasing.
	
	By \eqref{eq:Schur-operator} and the Trotter product formula \cite[\S5.3, Theorem~5.12]{Teschl2014},
	\[
	e^{-sS_\lambda}
	=e^{s\lambda}
	\lim_{n\to\infty}
	\left(
	e^{-sL/n}
	e^{\,spP/[n(p-\lambda)]}
	\right)^n,
	\]
	where the limit is in the strong operator topology. By Lemma~\ref{lem:L-monotone}, $e^{-sL/n}$ preserves radial monotonicity. Since $\lambda<p$, the coefficient $sp/[n(p-\lambda)]$ is nonnegative, so $e^{\,spP/[n(p-\lambda)]}$ shares this property by the first part. Every Trotter factor therefore preserves radial monotonicity, and the scalar $e^{s\lambda}>0$. Passing to a subsequence in the Trotter limit gives the result.
\end{proof}

\subsection{A monotone zero mode of the Schur complement}

\begin{proposition}
	\label{prop:monotone-zero-mode}
	Suppose
	$
	\sigma(\mathcal L_r)\cap(0,p-1]\ne\varnothing.
	$
	Then there exist \(\lambda\in(0,p-1]\) and a nonzero, nonconstant,
	radial nonincreasing function \(f\in\ker S_\lambda\) such that
	\[
	\int_{\mathbb R^N}f\,\dd\mu=0.
	\]
\end{proposition}

\begin{proof}
	By Proposition~\ref{prop:joint-spectrum}, zero is an isolated simple
	eigenvalue of \(\mathcal L_r\), and the positive eigenvalues below the
	threshold \(\tau\) can accumulate only at \(\tau\).  Since
	\(p-1<\tau\), the assumed set has a smallest element; denote it by
	\(\lambda\).
	Since \(\mathcal L_r\geq0\) and \(\lambda\) is its smallest positive spectral
	point, the only negative eigenvalue of
	\(\mathcal L_r-\lambda\) is \(-\lambda\), and it is simple.  Moreover,
	\(\lambda\) is an eigenvalue of \(\mathcal L_r\), so
	$
	\ker(\mathcal L_r-\lambda)\ne\{0\}.
	$
	Therefore, by \eqref{eq:Schur-kernel}, it holds
	$
	i_-(S_\lambda)=1,
	\
	\ker S_\lambda\ne\{0\}.
	$
	Moreover, \(L1=1\), \(P1=1\), and \eqref{eq:Schur-operator} give
	\[
	S_\lambda1
	=-\frac{\lambda(p+1-\lambda)}{p-\lambda}\,1<0.
	\]
	Since \(i_-(S_\lambda)=1\), the negative spectral subspace
	of \(S_\lambda\) is exactly \(\operatorname{span}\{1\}\).
	
	Let \(\Pi_\lambda\) be the orthogonal projection in
	\(L^2_r(\dd\mu)\) onto \(\ker S_\lambda\).
	We claim that there exists a radial nonincreasing function \(h\) such that
	\(\Pi_\lambda h\ne0\).
	Suppose to the contrary that
	\(\Pi_\lambda h=0\) for every radial nonincreasing function \(h\).
	In particular,
	\[
	\Pi_\lambda\mathbf 1_{B_\varrho}=0
	\qquad\text{for every }\varrho>0.
	\]
	For \(0<\varrho_1<\varrho_2\), one has
	\(\mathbf 1_{\{\varrho_1<|x|\leq\varrho_2\}}
	=\mathbf 1_{B_{\varrho_2}}-\mathbf 1_{B_{\varrho_1}}\), and hence
	\(\Pi_\lambda\mathbf 1_{\{\varrho_1<|x|\leq\varrho_2\}}=0\).  Their finite linear combinations are
	the radial step functions, which are dense in \(L^2_r(\dd\mu)\).
	Since \(\Pi_\lambda\) is bounded, this would imply
	\(\Pi_\lambda=0\), contradicting \(\ker S_\lambda\ne\{0\}\).
	Therefore the claim holds.
	
	Set
	\[
	\widetilde h
	=h-\frac{\langle h,1\rangle_{L^2_r(\dd\mu)}}{\|1\|_{L^2_r(\dd\mu)}^2}\,1.
	\]
	Then \(\widetilde h\) is radial nonincreasing and
	\(\widetilde h\perp1\) in \(L^2_r(\dd\mu)\).  Since \(S_\lambda\) is self-adjoint,
	\(\ker S_\lambda\) is orthogonal to its negative spectral subspace.
	Consequently,
	$
	\Pi_\lambda\widetilde h
	=\Pi_\lambda h\ne0.
	$
	Thus \(\widetilde h\) has no component in the negative spectral
	subspace of \(S_\lambda\).  Since
	\(\lambda\leq p-1<\tau\), Proposition~\ref{prop:joint-spectrum} shows
	that \(S_\lambda\) is Fredholm.  Hence zero is isolated in
	\(\sigma(S_\lambda)\), and the positive spectrum of \(S_\lambda\) is
	bounded away from zero.
	Lemma~\ref{lem:kernel-projection} therefore gives
	\[
	e^{-sS_\lambda}\widetilde h
	\to
	\Pi_\lambda\widetilde h
	\quad\text{in }L^2_r(\dd\mu)\quad\text{as \(s\to+\infty\).}
	\]
	By Proposition~\ref{prop:Schur-monotone}, every function on the left
	is radial nonincreasing.  Passing to a subsequence, we conclude that
	$
	f:=\Pi_\lambda\widetilde h
	$
	is radial nonincreasing and nonzero.
	By the orthogonality noted above, \(f\perp1\).  Hence
	\[
	\int_{\RN}f\,\dd\mu=0,
	\]
	and \(f\) is not constant.
\end{proof}

\section{Cumulative identities and the radial spectral gap}
\label{sec:spectral-gap}

Throughout this section, \((\alpha,p)\) satisfy \eqref{eq:range}.
As the remark following Lemma~\ref{lem:combined-prefix} explains, the
cumulative identities are in fact valid for every \(0<\alpha<N\) in the
variational range \eqref{eq:p-range}.

Fix the positive radial solution \(u\) chosen in
Section~\ref{sec:half-potential}.  Recall from \eqref{eq:A-def} and
\eqref{eq:V-def} that \(A=I_\alpha*u^p\) and \(V=Au^{p-2}\), write
\(\dd\mu=Au^p\dd x\), and, for \(r>0\), set
\[
X(r)=-\frac{ru'(r)}{u(r)},
\qquad
Y(r)=-\frac{rA'(r)}{A(r)}.
\]
Both functions are extended radially to \(\RN\), so that
\(X(x)=X(|x|)\) and \(Y(x)=Y(|x|)\).
The radial linearization \(L_{+,u}:H^1_r(\RN)\to H^{-1}_r(\RN)\) is the
scalar operator obtained by eliminating the second field in
\eqref{eq:linearized-two-field}; it is defined weakly by
\[
\bigl\langle L_{+,u}\phi,\psi\bigr\rangle
=\int_{\RN}(\nabla\phi\cdot\nabla\psi+\phi\psi)\,\dd x
-(p-1)\int_{\RN}Au^{p-2}\phi\psi\,\dd x
-p\int_{\RN}\bigl(I_\alpha*(u^{p-1}\phi)\bigr)u^{p-1}\psi\,\dd x,
\]
with formal expression
\(
L_{+,u}\phi=-\Delta\phi+\phi-(p-1)Au^{p-2}\phi
-p\,u^{p-1}\bigl(I_\alpha*(u^{p-1}\phi)\bigr)
\).
As in the introduction, this form is symmetric, closed and bounded from
below on \(H^1_r(\RN)\); it is the form of a self-adjoint operator on
\(L^2_r(\RN)\), which we also denote by \(L_{+,u}\).  Since it acts on
\(H^1_r(\RN)\), all indices and kernels below are understood in the
radial space, with the convention of Definition~\ref{def:index-nullity}.

We shall use the following two tail estimates in the subquadratic case
\(\alpha=2\), \(p<2\).  In the radial variable the decay of
Lemma~\ref{lem:regularity-and-far-field} reads
\(u(r)=(c_{N,\alpha}M)^{1/(2-p)}r^{-\beta_p}(1+o(1))\), where
\(\beta_p=\frac{N-2}{2-p}\) and \(p\beta_p>N+2\).
The continuity of \(V\) and Lemma~\ref{lem:regularity-and-far-field}
imply that \(V\) is bounded.  Applying the interior second-derivative estimate
to \(-\Delta u=(V-1)u\) on unit balls centered at points of large
radius, and then using the radial equation, gives
\begin{equation}\label{eq:u-subquadratic-derivatives}
	|u'(r)|+|u''(r)|=O(r^{-\beta_p}).
\end{equation}

\subsection{Cumulative identities and the strict sign of the weighted pairing}

The goal of this subsection is the strict inequality in
Proposition~\ref{prop:weight-negative} below, which is the first of
the three relations that Proposition~\ref{prop:joint-gap} combines.
It is
obtained in two steps: the cumulative identities of
Lemma~\ref{lem:combined-prefix} express the Riesz kernel through the
prefix integrals of \(H_\lambda\), and the one-dimensional duality of
Lemma~\ref{lem:prefix-duality} transfers their sign to the pairing with
the monotone zero mode \(f\).

\begin{lemma}\label{lem:combined-prefix}
	For every \(\varrho>0\),
	\begin{equation}\label{eq:Riesz-prefix}
		\begin{aligned}
			&\int_{|x|<\varrho}
			\left(Y(x)-\frac{N-\alpha}{2}\right)\dd\mu\\
			&\quad=
			\frac{c_{N,\alpha}(N-\alpha)}{2}
			\iint_{|x|<\varrho<|y|}
			u(x)^pu(y)^p
			\frac{|x|^2-|y|^2}{|x-y|^{N-\alpha+2}}
			\,\dd y\,\dd x<0.
		\end{aligned}
	\end{equation}
	Moreover, for every \(0<\lambda<p\),
	\begin{equation}\label{eq:combined-prefix}
		\begin{aligned}
			&\int_{|x|<\varrho}
			\left[
			X(x)-\frac{N+\alpha}{2p}
			+\frac{Y(x)-(N-\alpha)/2}{p-\lambda}
			\right]\dd\mu\\
			&\quad=
			-\frac{|\mathbb S^{N-1}|}{p}\,\varrho^NA(\varrho)u(\varrho)^p
			+\frac{\lambda}{p(p-\lambda)}
			\int_{|x|<\varrho}
			\left(Y(x)-\frac{N-\alpha}{2}\right)\dd\mu<0.
		\end{aligned}
	\end{equation}
	The left-hand sides of \eqref{eq:Riesz-prefix} and
	\eqref{eq:combined-prefix} tend to zero as
	\(\varrho\to+\infty\).
\end{lemma}

\begin{proof}
	By the regularity and decay properties recorded in
	Lemma~\ref{lem:regularity-and-far-field},
	\(u^p\in L^1(\RN)\cap L^\infty(\RN)\).  In addition,
	\eqref{eq:u-exponential-decay} gives exponential decay when
	\(p\geq2\), whereas \eqref{decayp} of Lemma~\ref{lem:regularity-and-far-field} gives
	\(p\beta_p>N+2\) when \(\alpha=2\) and \(p<2\).  Thus in both cases
	\(\int_{\RN}(1+|x|)u(x)^p\,\dd x<\infty\).
	The differentiation under the convolution was justified in the proof
	of \eqref{eq:A-prime-far-field}; hence \(\nabla A=(\nabla I_\alpha)*u^p\), and
	a direct calculation gives
	\begin{equation}\label{eq:Riesz-kernel}
		\begin{aligned}
			-x\cdot\nabla_xI_\alpha(x-y)
			-\frac{N-\alpha}{2}I_\alpha(x-y)
			&=
			c_{N,\alpha}(N-\alpha)
			\frac{x\cdot(x-y)-|x-y|^2/2}
			{|x-y|^{N-\alpha+2}}\\
			&=
			\frac{c_{N,\alpha}(N-\alpha)}{2}
			\frac{|x|^2-|y|^2}{|x-y|^{N-\alpha+2}}.
		\end{aligned}
	\end{equation}
	The kernel in the last line of \eqref{eq:Riesz-kernel} is absolutely
	integrable against \(u(x)^pu(y)^p\).  Indeed,
	$
	\bigl||x|^2-|y|^2\bigr|
	\leq(|x|+|y|)|x-y|,
	$
	so its singularity at \(x=y\) has order
	\(N-\alpha+1<N\), while the preceding first-moment estimate controls
	infinity.  Hence Fubini's theorem, \eqref{eq:Riesz-kernel}, and radial
	integration give
	\[
	\begin{aligned}
		|\mathbb S^{N-1}|
		\int_0^\varrho r^{N-1}Au^p
		&\left(Y-\frac{N-\alpha}{2}\right)\,\dd r
		=
		\int_{|x|<\varrho}u(x)^p
		\left[
		-x\cdot\nabla A(x)-\frac{N-\alpha}{2}A(x)
		\right]\,\dd x\\
		&=
		\frac{c_{N,\alpha}(N-\alpha)}{2}
		\int_{|x|<\varrho}\int_{\RN}
		u(x)^pu(y)^p
		\frac{|x|^2-|y|^2}{|x-y|^{N-\alpha+2}}
		\,\dd y\,\dd x.
	\end{aligned}
	\]
	The integral over \(B_\varrho\times B_\varrho\) vanishes after
	interchanging \(x\) and \(y\).  On the remaining region
	\(|x|<\varrho<|y|\), the numerator is strictly negative.  This proves
	\eqref{eq:Riesz-prefix}.
	The indicator of the region
	\(\{|x|<\varrho<|y|\}\) tends pointwise to zero as
	\(\varrho\to+\infty\).  The absolute integrability proved above
	therefore allows us to apply dominated convergence, yielding
	\[
	\int_0^\varrho r^{N-1}Au^p
	\left(Y-\frac{N-\alpha}{2}\right)\,\dd r
	\to 0
	\qquad\text{as }\varrho\to+\infty.
	\]
	
	Finally,
	\(Xu^p=-r(u^p)'/p\).  Since the boundary term at \(r=0\) vanishes,
	integration by parts gives
	\begin{equation}\label{eq:X-prefix}
		\begin{aligned}
			\int_0^\varrho r^{N-1}Au^p
			&\left(X-\frac{N+\alpha}{2p}\right)\,\dd r
			=
			-\frac1p\int_0^\varrho r^NA(u^p)'\,\dd r
			-\frac{N+\alpha}{2p}
			\int_0^\varrho r^{N-1}Au^p\,\dd r\\
			&=
			-\frac{\varrho^NA(\varrho)u(\varrho)^p}{p}
			+\frac1p\int_0^\varrho r^{N-1}Au^p
			\left(
			N-Y-\frac{N+\alpha}{2}
			\right)\,\dd r\\
			&=
			-\frac{\varrho^NA(\varrho)u(\varrho)^p}{p}
			-\frac1p\int_0^\varrho r^{N-1}Au^p
			\left(Y-\frac{N-\alpha}{2}\right)\,\dd r.
		\end{aligned}
	\end{equation}
	Multiplying \eqref{eq:X-prefix} by \(|\mathbb S^{N-1}|\) immediately yields
	\eqref{eq:combined-prefix}.  Its strict negativity follows from
	\eqref{eq:Riesz-prefix}.  Finally,
	\[
	\varrho^NA(\varrho)u(\varrho)^p\to 0,
	\]
	by \eqref{eq:u-exponential-decay} and \eqref{eq:A-prime-far-field} if
	\(p\geq2\), and by \eqref{eq:A-prime-far-field} and Item (5) of
	Lemma~\ref{lem:regularity-and-far-field} if \(\alpha=2\) and \(p<2\).
	Together with the first limit, this proves the second asserted limit.
\end{proof}

\begin{remark}
	Lemma~\ref{lem:combined-prefix} has been stated and proved under
	the parameter range \eqref{eq:range} used in this paper, but its
	assertions hold in the full variational range: for every
	\(0<\alpha<N\) and  \eqref{eq:p-range}. 
\end{remark}

\begin{lemma}\label{lem:prefix-duality}
	Let \(\omega\) be a positive, nonatomic, \(\sigma\)-finite Borel
	measure on \((0,\infty)\).  Suppose that \(H\in L^1(\dd\omega)\),
	$
	\int_0^\infty H\,\dd\omega=0,
	\
	\int_0^\varrho H\,\dd\omega<0
	\ \text{for every }\varrho>0.
	$
	If \(f\) is nonconstant and nonincreasing, and
	\(Hf\in L^1(\dd\omega)\), then
	$$
	\int_0^\infty Hf\,\dd\omega<0.
	$$
\end{lemma}

\begin{proof}
	Set \(F(\varrho)=\int_{(0,\varrho)}H\,\dd\omega\) for \(\varrho>0\), so
	that \(F(\varrho)<0\) for every \(\varrho>0\) and \(F(\infty)=0\).
	Since \(\int_0^\infty H\,\dd\omega=0\), we may suppose that
	\(f(1)=0\): replacing \(f\) by \(f-f(1)\) preserves the hypotheses
	on \(f\) and leaves \(\int_0^\infty Hf\,\dd\omega\) unchanged.
	As \(f\) is nonincreasing, there is a positive \(\sigma\)-finite
	Borel measure \(\nu\) on \((0,\infty)\) with
	\[
	\nu((a,b])=f(a)-f(b),
	\qquad 0<a<b<\infty.
	\]
	Since \(f\) is nonconstant, \(\nu\ne0\).  As \(f(1)=0\),
	\[
	f(x)=\nu((x,1])\ \ (x<1),
	\qquad
	f(x)=-\nu((1,x])\ \ (x>1).
	\]
	Substituting this representation of \(f\) and exchanging the order
	of integration gives
	\[
	\begin{aligned}
		\int_0^\infty Hf\,\dd\omega
		&=\int_{(0,1]}\left(\int_{(0,y)}H\,\dd\omega\right)\dd\nu(y)
		-\int_{(1,\infty)}\left(\int_{(y,\infty)}H\,\dd\omega\right)\dd\nu(y)\\
		&=\int_{(0,1]}F\,\dd\nu+\int_{(1,\infty)}F\,\dd\nu
		=\int_{(0,\infty)}F\,\dd\nu.
	\end{aligned}
	\]
	In this computation, \(\int_0^\infty Hf\,\dd\omega\) is split at
	\(1\), and \(\int_{(y,\infty)}H\,\dd\omega=-F(y)\); each step
	discards only sets of zero \(\omega\)-measure, which is precisely
	what the nonatomicity of \(\omega\) supplies.  The exchange of the integrals in
	the two terms is justified by Tonelli's theorem, since \(|f|\) equals
	\(\nu((x,1])\) on \((0,1)\) and \(\nu((1,x])\) on \((1,\infty)\),
	whereas \(Hf\in L^1(\dd\omega)\).  Since \(F(y)<0\) for every
	\(y>0\) and \(\nu\) is a nonzero positive measure, the last integral
	is strictly negative, and the proof is complete.
\end{proof}

\begin{proposition}\label{prop:weight-negative}
	Let \(\lambda\in(0,p-1]\), and let \(f\in\ker S_\lambda\subset L^2_r(\dd\mu)\)
	be nonconstant and radial nonincreasing and satisfy
	\[
	\int_{\RN}f\,\dd\mu=0.
	\]
	Then
	\[
	\int_{\RN}
	\left(X+\frac{Y}{p-\lambda}\right)f\,\dd\mu<0.
	\]
\end{proposition}

\begin{proof}
	Set
	\[
	\dd\omega=|\mathbb S^{N-1}|\,r^{N-1}A(r)u(r)^p\,\dd r,
	\qquad
	H_\lambda
	=
	X-\frac{N+\alpha}{2p}
	+\frac{Y-(N-\alpha)/2}{p-\lambda}.
	\]
	Since \(A,u>0\) and \(\mu(\RN)<\infty\), \(\omega\) is a finite
	measure with strictly positive density.  For every radial integrable
	function \(g\),
	\[
	\int_{\RN}g\,\dd\mu
	=
	\int_0^\infty g\,\dd\omega.
	\]
	When \(p\geq2\), the regularity at the origin,  \eqref{eq:A-tail} and
	\eqref{eq:A-prime-far-field}, the asymptotic behavior of $V$  
	\eqref{eq:u-exponential-decay} give
	$
	X\in L^2(\dd\omega)\cap L^1(\dd\omega),
	\
	Y\in L^\infty(0,\infty).
	$
	When \(\alpha=2\) and \(p<2\), the same conclusion follows from
	\eqref{eq:A-tail}, \eqref{eq:A-prime-far-field}, \eqref{eq:V-limit}, and
	\eqref{eq:u-subquadratic-derivatives}.  Indeed, in that case
	\(X(r)=O(r)\), \(Y(r)=O(1)\), and the density of \(\omega\) is
	\(O(r^{1-p\beta_p})\) at infinity, with \(p\beta_p>N+2\).
	Hence
	\[
	|H_\lambda|
	\leq C_\lambda(1+|X|+|Y|)
	\in L^2(\dd\omega)\cap L^1(\dd\omega).
	\]
	By Lemma~\ref{lem:combined-prefix},
	\[
	\int_0^\varrho H_\lambda\,\dd\omega<0\quad\text{for every }\varrho>0,
	\qquad
	\int_0^\infty H_\lambda\,\dd\omega=0.
	\]
	The radial integration formula gives
	\[
	f\in L^2(\dd\omega),
	\qquad
	\int_0^\infty f\,\dd\omega=0.
	\]
	Therefore, by the Cauchy--Schwarz inequality,
	$H_\lambda f\in L^1(\dd\omega)$.
	
	All the assumptions of Lemma~\ref{lem:prefix-duality} are therefore
	satisfied, and hence
	\[
	\int_0^\infty H_\lambda f\,\dd\omega<0.
	\]
	Using the definition of \(H_\lambda\) and
	\(\int_0^\infty f\,\dd\omega=0\), we obtain
	\[
	\begin{aligned}
		0
		&>
		\int_0^\infty H_\lambda f\,\dd\omega
		=
		\int_0^\infty
		\left(X+\frac{Y}{p-\lambda}\right)f\dd\omega
		=
		\int_{\RN}
		\left(X+\frac{Y}{p-\lambda}\right)f\,\dd\mu,
	\end{aligned}
	\]
	which is the assertion.
\end{proof}

\subsection{The Newtonian cumulative inequality}
\label{subsec:newtonian}

We now assume that \(\alpha=2\) and
\(\frac{N+2}{N}<p<\min\{2,\frac{N+2}{N-2}\}\), the Newtonian
subquadratic range, in which the monotonicity of \(V=Au^{p-2}\) is
not available (see Remark~\ref{rem:V-not-monotone} below).

The equations \eqref{eq:choquard} and \eqref{eq:A-def} give the local system
\begin{equation}\label{eq:local-system}
	-\Delta u+u=Au^{p-1},
	\qquad
	-\Delta A=u^p.
\end{equation}
Both \(u\) and \(A\) tend to zero at infinity.

Define the constant
\begin{equation*}\label{eq:V-star}
	V_*:=\frac{2p}{N+2-(N-2)p}>1.
\end{equation*}
The denominator here is positive because \(p<\frac{N+2}{N-2}\),
so \(V_*>0\); moreover \(V_*>1\) is equivalent to
\(p>\frac{N+2}{N}\), which holds by assumption.
Testing the first equation in \eqref{eq:local-system} by $u$ and using
the Pohozaev identity gives
\begin{align*}\label{eq:Nehari-id}
	\int_{\RN}(|\nabla u|^2+u^2)\dd x
	=\int_{\RN}Au^p\dd x,   \quad 
	\frac{N-2}{2}\int_{\RN}|\nabla u|^2\dd x
	+\frac N2\int_{\RN}u^2\dd x
	=\frac{N+2}{2p}\int_{\RN}Au^p\dd x.                   
\end{align*}
Eliminating the gradient term and using $Au^p=Vu^2$, we obtain
\begin{equation}\label{eq:V-average}
	\displaystyle\int_0^\infty \left(1-\frac{V(r)}{V_*}\right) u(r)^2r^{N-1}\dd r=0.
\end{equation}
Writing $\bar V$ for the $u^2$-weighted running average of $V$,
\begin{equation}\label{eq:V-bar-def}
	\bar V(r):=
	\frac{\displaystyle\int_0^r s^{N-1}u(s)^2V(s)\dd s}
	{\displaystyle\int_0^r s^{N-1}u(s)^2\dd s},
	\qquad r>0,
\end{equation}
we may read \eqref{eq:V-average} as saying that this running average
tends to the global average $V_*$.  
Namely,
\[
\lim_{r\to\infty}\bar V(r)
=\frac{\displaystyle\int_0^\infty s^{N-1}u(s)^2V(s)\dd s}
{\displaystyle\int_0^\infty s^{N-1}u(s)^2\dd s}
=V_* .
\]

\begin{lemma} \label{lem:local-pohozaev}
	For every $r>0$,
	\begin{equation}\label{eq:local-pohozaev}
		\int_0^r s^{N-1}u(s)^2
		\left(1-\frac{V(s)}{V_*}\right)\dd s
		=-\frac12r^{N-2}u(r)^2E(r),
	\end{equation}
	where \begin{equation}\label{eq:E-def}
		E(r)=rX'
		+\frac{A^2}{p u^2}rY'
		-\frac{p-1}{p}r^2V.
	\end{equation}
\end{lemma}

\begin{proposition}\label{prop:E-positive}
	Assume that $\alpha=2$ and
	$\frac{N+2}{N}<p<\min\{2,\frac{N+2}{N-2}\}$.  Then, for every $r>0$,
	\begin{equation*}\label{eq:E-positive}
		E(r)>0.
	\end{equation*}
	Equivalently, $\bar V(r)>V_*$ for every $r>0$, where $\bar V$ is the
	$u^2$-weighted running average of $V$ defined in
	\eqref{eq:V-bar-def}; explicitly,
	\begin{equation}\label{eq:prefix-negative}
		\int_0^r s^{N-1}u(s)^2
		\left(1-\frac{V(s)}{V_*}\right)\dd s<0.
	\end{equation}
\end{proposition}

The proofs of Lemma~\ref{lem:local-pohozaev} and of
Proposition~\ref{prop:E-positive} are given in
Appendix~\ref{app:newtonian}.

\subsection{Scaling identities and the radial spectral gap}

\begin{lemma}\label{lem:scaling}
	Set
	\[
	\Lambda:=\frac{\alpha+2}{2(p-1)}-X,
	\qquad\text{so that}\qquad
	u\Lambda=\frac{\alpha+2}{2(p-1)}u+x\cdot\nabla u .
	\]
	Then \(u\Lambda\in H^1_r(\RN)\), and
	\begin{equation}\label{eq:L-scaling}
		L_{+,u}(u\Lambda)=-2u,\qquad
		PX=\frac{Y+\alpha}{p}.
	\end{equation}
\end{lemma}

\begin{proof}
	Let \(\gamma=(\alpha+2)/(2(p-1))\).  For \(\tau>0\), set
	$
	u_\tau(x)=\tau^\gamma u(\tau x),
	\  A_\tau=I_\alpha*u_\tau^p.
	$
	The homogeneity in \eqref{eq:Riesz-normalization} and the definition
	\eqref{eq:A-def} give
	\[
	A_\tau(x)=\tau^{p\gamma-\alpha}A(\tau x).
	\]
	Since \(2(p-1)\gamma=\alpha+2\), equation
	\eqref{eq:choquard} yields
	\[
	\begin{aligned}
		-\Delta u_\tau+u_\tau-A_\tau u_\tau^{p-1}
		&=\tau^{\gamma+2}
		\bigl(-\Delta u-Au^{p-1}\bigr)(\tau x)
		+\tau^\gamma u(\tau x) =(1-\tau^2)u_\tau.
	\end{aligned}
	\]
	Since the \(\tau\)-derivative of \(u_\tau\) at \(\tau=1\) is
	\(\gamma u+x\cdot\nabla u=u\Lambda\), differentiating this identity at
	\(\tau=1\), in the weak sense used to define \(L_{+,u}\), gives the
	first identity in \eqref{eq:L-scaling}. The decay estimate
	\eqref{eq:u-exponential-decay} gives
	\(u\Lambda\in H^1_r(\RN)\) when \(p\geq2\). In the Newtonian
	subquadratic case the same conclusion follows from Item (5) of
	Lemma~\ref{lem:regularity-and-far-field} and
	\eqref{eq:u-subquadratic-derivatives}.
	
	For the second identity, homogeneity gives
	\begin{equation}\label{eq:density-dilation}
		I_\alpha*\bigl(u(\tau\,\cdot)^p\bigr)(x)
		=\tau^{-\alpha}A(\tau x).
	\end{equation}
	By \eqref{eq:P-Markov}, the definition of \(X\), differentiation of
	\eqref{eq:density-dilation} at \(\tau=1\), and the definition of \(Y\),
	respectively, we obtain
	\[
	\begin{aligned}
		-pA(x)PX(x)
		&=-pI_\alpha*(Xu^p)(x)\\
		&=I_\alpha*\bigl(y\cdot\nabla(u^p)(y)\bigr)(x) =-\alpha A(x)+x\cdot\nabla A(x) =-(\alpha+Y(x))A(x).
	\end{aligned}
	\]
	Since \(A(x)>0\), dividing by \(-pA(x)\) gives
	the second identity in \eqref{eq:L-scaling} .
\end{proof}

\begin{lemma}\label{lem:Schur-pairing}
	Let \(0<\lambda\leq p-1\), and let \(f\in\ker S_\lambda\) satisfy
	\(\int_{\RN}f\,\dd\mu=0\).  Then
	\begin{equation}\label{eq:exact-pairing}
		-(p-1-\lambda)\int_{\RN}f
		\left(X+\frac{Y}{p-\lambda}\right)\,\dd\mu
		=2\int_{\RN}u^2f\,\dd x.
	\end{equation}
\end{lemma}

\begin{proof}
	By Lemma~\ref{lem:scaling} and \eqref{eq:calD-realization},
	\(\Lambda\in\mathcal D\).  By \eqref{eq:scalar-ground-transform}, applied to the
	pair \((\Lambda,f)\), and by \eqref{eq:L-scaling},
	\begin{equation}\label{eq51}
		\begin{aligned}
			s_{p-1}(\Lambda,f)
			&=\bigl\langle L_{+,u}(u\Lambda),uf
			\bigr\rangle_{H^{-1}_r,H^1_r} =-2\int_{\RN}u^2f\,\dd x.
		\end{aligned}
	\end{equation}
	Taking \(h=\Lambda\) in \eqref{eq:Schur-form}, and comparing
	the cases \(\lambda\) and \(p-1\), we obtain
	\begin{equation}\label{eq:Schur-shift}
		\begin{aligned}
			0=s_\lambda(f,\Lambda)
			={}&s_{p-1}(f,\Lambda)
			+(p-1-\lambda)\langle f,\Lambda\rangle_{L^2_r(\dd\mu)} +\frac{p(p-1-\lambda)}{p-\lambda}
			\langle Pf,\Lambda\rangle_{L^2_r(\dd\mu)},
		\end{aligned}
	\end{equation}
	where the first equality follows from \(f\in\ker S_\lambda\).
	Moreover, using the definition of \(\Lambda\) and
	\(\int_{\RN}f\,\dd\mu=0\), we obtain
	\[
	\langle f,\Lambda\rangle_{L^2_r(\dd\mu)}
	=-\int_{\RN}fX\,\dd\mu.
	\]
	Since \(P\) is self-adjoint, \(P1=1\), and
	\(PX=(Y+\alpha)/p\),
	\[
	\begin{aligned}
		\langle Pf,\Lambda\rangle_{L^2_r(\dd\mu)}
		&=\langle f,P\Lambda\rangle_{L^2_r(\dd\mu)}=-\langle f,PX\rangle_{L^2_r(\dd\mu)} =-\frac1p\int_{\RN}fY\,\dd\mu.
	\end{aligned}
	\]
	Substituting these identities and
	\eqref{eq51}
	into \eqref{eq:Schur-shift} gives
	\[
	\begin{aligned}
		0
		={}&-2\int_{\RN}u^2f\,\dd x
		-(p-1-\lambda)\int_{\RN}fX\,\dd\mu -\frac{p-1-\lambda}{p-\lambda}
		\int_{\RN}fY\,\dd\mu.
	\end{aligned}
	\]
	Rearranging yields \eqref{eq:exact-pairing}.
\end{proof}

\begin{lemma}\label{lem:Chebyshev}
	If \(f\in\mathcal D\) is radially nonincreasing, is not constant almost
	everywhere, and satisfies
	\(\int_{\RN}f\,\dd\mu=0\), then
	\begin{equation}\label{eq:u2f}
		\int_{\RN}u^2f\,\dd x<0.
	\end{equation}
\end{lemma}

\begin{proof}
	Since
	\(V=Au^{p-2}\) by \eqref{eq:V-def} and
	\(\dd\mu=Au^p\,\dd x\) by \eqref{eq:joint-measures},  
	\begin{equation}\label{eq54}
		\int_{\RN}Vfu^2\,\dd x
		=\int_{\RN}f\,\dd\mu=0.
	\end{equation}
	
	Suppose first that \(p\geq2\).  We prove the sign of
	\(\int_{\RN}u^2f\,\dd x\) by the strict Chebyshev argument. 
	By Lemma~\ref{lem:ground-state-transform},
	\(uf\in L^2(\RN)\), while
	Lemma~\ref{lem:half-potential-marginals} gives
	\(\mu(\RN)<\infty\).  Hence, by the Cauchy--Schwarz inequality,
	\[
	\begin{aligned}
		\int_{\RN}u^2|f|\,\dd x
		&\leq \|u\|_2\|uf\|_2<\infty,\qquad
		\int_{\RN}Vu^2|f|\,\dd x
		=\int_{\RN}|f|\,\dd\mu
		\leq \mu(\RN)^{1/2}\|f\|_{L^2_r(\dd\mu)}<\infty,
	\end{aligned}
	\]
	where the second equality follows from
	\eqref{eq:V-def} and \eqref{eq:joint-measures}.
	
	By \eqref{eq:range}, \eqref{eq:V-def}, Theorem~\ref{thm:symmetry},
	and Item (4) of Lemma ~\eqref{lem:regularity-and-far-field}, \(V=Au^{p-2}\) is strictly
	radially decreasing.  Hence
	$
	(V(x)-V(y))(f(x)-f(y))\geq0
	\ \text{for all }x,y\in\RN,
	$
	because \(V\) and \(f\) are both radially nonincreasing.
	This inequality is strict on a set of positive measure.  Indeed, since
	the radial representative of \(f\) is not constant, there exist
	\(0<r_1<r_2\) such that \(f(r_1)>f(r_2)\).  Thus, whenever
	$
	|x|<r_1<r_2<|y|,
	$
	we have
	\[
	V(x)>V(y),
	\qquad
	f(x)\geq f(r_1)>f(r_2)\geq f(y).
	\]
	Since \(u>0\), the preceding integrability estimates and Fubini's
	theorem therefore give
	\[
	\begin{aligned}
		0
		<&\frac12\iint_{\RN\times\RN}
		(V(x)-V(y))(f(x)-f(y))
		u(x)^2u(y)^2\,\dd x\dd y\\
		={}&
		\left(\int_{\RN}Vfu^2\,\dd x\right)
		\left(\int_{\RN}u^2\,\dd x\right) -
		\left(\int_{\RN}Vu^2\,\dd x\right)
		\left(\int_{\RN}fu^2\,\dd x\right)\\
		={}&-
		\left(\int_{\RN}Vu^2\,\dd x\right)
		\left(\int_{\RN}fu^2\,\dd x\right).
	\end{aligned}
	\]
	Finally, by
	$
	\int_{\RN}Vu^2\,\dd x
	>0,
	$
	we get \eqref{eq:u2f} in this case.
	
	It remains to consider \(\alpha=2\) and \(p<2\).  Let
	\[
	\dd\omega_0=|\mathbb S^{N-1}|\,r^{N-1}u(r)^2\,\dd r,
	\qquad
	H_0(r)=1-\frac{V(r)}{V_*}.
	\]
	By \eqref{eq:V-average} and \eqref{eq:prefix-negative},
	\[
	\int_0^\infty H_0\,\dd\omega_0=0,
	\qquad
	\int_0^\varrho H_0\,\dd\omega_0<0
	\quad\text{for every }\varrho>0.
	\]
	The measure \(\omega_0\) is finite and has a strictly positive density.
	Moreover, \(V\in L^\infty(0,\infty)\) by
	Lemma~\ref{lem:regularity-and-far-field}, so that the factor
	\(H_0=1-V/V_*\) is bounded, say \(|H_0|\le C\). Then
	\[
	\int_0^\infty|H_0f|\,\dd\omega_0
	\le C\int_{\RN}u^2|f|\,\dd x<\infty,
	\]
	that is, \(H_0f\in L^1(\dd\omega_0)\).  Lemma~\ref{lem:prefix-duality}
	therefore yields
	\[
	\int_0^\infty H_0f\,\dd\omega_0<0.
	\]
	Since \(H_0=1-V/V_*\) and \(\int_{\RN}Vfu^2\,\dd x=0\) by \eqref{eq54},
	\[
	\int_0^\infty H_0f\,\dd\omega_0
	=\int_{\RN}fu^2\,\dd x,
	\]
	which again proves \eqref{eq:u2f}.
\end{proof}

\begin{remark}\label{rem:V-not-monotone}
The Chebyshev argument used for \(p\ge2\) requires \(V=Au^{p-2}\) to be nonincreasing.  In the range \(p<2\), no such monotonicity can be established: \(V\) is the product of the strictly decreasing factor \(A\) and the strictly increasing factor \(u^{p-2}\), and the far-field behavior \(V\to1\) in~\eqref{eq:V-limit} does not restore it. This is the only place in the proof where \(\alpha=2\)
	is used.  The factorization, the radial total positivity, the
	cumulative identities and the compactness estimates of
	Sections~\ref{sec:half-potential}--\ref{sec:nehari} are available for
	every \(2\le\alpha<N\); only the strict sign of
	Lemma~\ref{lem:Chebyshev} in the subquadratic range requires the
	comparison function \(H_0=1-V/V_*\) and the Newtonian cumulative
	inequality of Subsection~\ref{subsec:newtonian}, both of which rest
	on \(-\Delta A=u^p\).  For \(\alpha>2\) and \(p<2\) the endpoint
	\(\lambda=p-1\) of Proposition~\ref{prop:joint-gap} therefore remains
	open: assuming \(i_{-,H^1_r}(L_{+,u})=1\) one still excludes
	\((0,p-1)\) by counting negative directions, but not \(\lambda=p-1\),
	which is precisely the case in which the strict sign is needed.
\end{remark}

\begin{proposition}\label{prop:joint-gap}
	$
	\sigma(\mathcal L_r)\cap(0,p-1]=\varnothing.
	$
\end{proposition}

\begin{proof}
	Suppose that
	$
	\sigma(\mathcal L_r)\cap(0,p-1]\ne\varnothing.
	$
	By Proposition~\ref{prop:monotone-zero-mode}, there exist
	\(\lambda\in(0,p-1]\) and a nonzero, nonconstant, radially
	nonincreasing \(f\in\ker S_\lambda\) such that
	\[
	\int_{\RN}f\,\dd\mu=0.
	\]
	Moreover, \eqref{eq:Schur-operator} gives
	\(f\in D(S_\lambda)=D(L)\subset\mathcal D\), so that
	\eqref{eq:exact-pairing}, Proposition~\ref{prop:weight-negative},
	and \eqref{eq:u2f} all apply to \(f\).  Since
	\(p-1-\lambda\geq0\), these three identities give
	\[
	\begin{aligned}
		0
		&\leq -(p-1-\lambda)
		\int_{\RN}
		\left(X+\frac{Y}{p-\lambda}\right)f\,\dd\mu =2\int_{\RN}u^2f\,\dd x<0,
	\end{aligned}
	\]
	a contradiction.  Therefore
	$
	\sigma(\mathcal L_r)\cap(0,p-1]=\varnothing.
	$
\end{proof}

\begin{Theorem}\label{thm:index-one}
	Every positive radial solution in the range \eqref{eq:range} satisfies
	\[
	i_{-,H^1_r}(L_{+,u})=1,
	\qquad
	i_{0,H^1_r}(L_{+,u})=0.
	\]
\end{Theorem}

\begin{proof}
	By Lemma~\ref{lem:joint-operator},
	\(\sigma(\mathcal L_r)\subset[0,\infty)\), while
	\eqref{eq:joint-kernel-recalled} gives
	$
	\ker\mathcal L_r=\operatorname{span}\{(1,1)\}.
	$
	Thus \(0\) is a simple eigenvalue of \(\mathcal L_r\).  Moreover,
	Proposition~\ref{prop:joint-gap} gives
	$
	\sigma(\mathcal L_r)\cap(0,p-1]=\varnothing.
	$
	Since
	\[
	\sigma\bigl(\mathcal L_r-(p-1)\bigr)
	=\sigma(\mathcal L_r)-(p-1),
	\]
	the negative spectrum of \(\mathcal L_r-(p-1)\) consists only of
	the simple eigenvalue \(-(p-1)\), whose eigenspace is
	\(\operatorname{span}\{(1,1)\}\), and
	\(0\notin\sigma(\mathcal L_r-(p-1))\), because
	\(p-1\notin\sigma(\mathcal L_r)\).  Therefore
	$
	i_-\bigl(\mathcal L_r-(p-1)\bigr)=1,
	\
	i_0\bigl(\mathcal L_r-(p-1)\bigr)=0.
	$
	Proposition~\ref{prop:Schur-complement} with \(\lambda=p-1\) now
	yields
	\[
	i_-(S_{p-1})=1,
	\qquad
	i_0(S_{p-1})=0.
	\]
	
	By \eqref{eq:scalar-ground-transform} and
	Lemma~\ref{lem:ground-state-transform}, multiplication by \(u\) is
	an isometric isomorphism from \(\mathcal D\) onto \(H^1_r(\RN)\)
	that carries the form \(s_{p-1}\) onto the form of \(L_{+,u}\); it
	therefore matches form domains, kernels and negative subspaces, so
	that
	\[
	i_-(L_{+,u})=i_-(S_{p-1})=1,
	\qquad
	i_0(L_{+,u})=i_0(S_{p-1})=0.
	\]
\end{proof}

\section{The Nehari argument and global uniqueness}
\label{sec:nehari}

We fix one exponent $p$ satisfying \eqref{eq:range}.  On
$H^1_r(\RN)$ consider
\begin{equation*}\label{eq:energy}
	\mathcal J(v)=\frac12\|v\|_{H^1_r}^2
	-\frac1{2p}\int_{\RN}(I_\alpha*|v|^p)|v|^p
\end{equation*}
and its radial Nehari manifold
\begin{equation*}\label{eq:Nehari-manifold}
	\cN
	=
	\left\{
	v\in H^1_r(\RN)\setminus\{0\}:
	\|v\|_{H^1_r}^2
	=
	\int_{\RN}(I_\alpha*|v|^p)|v|^p
	\right\}.
\end{equation*}
The Hardy--Littlewood--Sobolev and Sobolev inequalities imply that the
$H^1_r$ norm is bounded away from zero on $\cN$.  Differentiating the
constraint in the radial direction $v$ gives
$
2(1-p)\|v\|_{H^1_r}^2<0.
$
Thus $\cN$ is a complete $C^1$ codimension-one manifold and a natural
constraint.
Moreover,
\begin{equation}\label{eq:energy-on-Nehari}
	\mathcal J(v)=\frac{p-1}{2p}\|v\|_{H^1_r}^2,
	\qquad v\in\cN.
\end{equation}

For completeness, $\mathcal J|_{\cN}$ satisfies the Palais--Smale
condition for this fixed $p$.  A constrained Palais--Smale sequence is
bounded by \eqref{eq:energy-on-Nehari}.  Pairing its multiplier equation
with the sequence shows that the Lagrange multipliers tend to zero.
The exponent $2Np/(N+\alpha)$ lies strictly between $2$ and $2N/(N-2)$ by
\eqref{eq:range}.  Therefore the compact radial embedding
$
H^1_r(\RN)\hookrightarrow L^{2Np/(N+\alpha)}(\RN),
$
together with the Hardy--Littlewood--Sobolev inequality, makes the
nonlinear derivative compact on bounded radial sets.  The sequence
converges strongly in $H^1_r$, and its limit is nonzero because of the
Nehari lower bound.

\subsection{A positive critical point is a Nehari well}

\begin{lemma}\label{lem:Nehari-well}
	For every positive radial critical point $u$, there exist
	$\rho_u,\epsilon_u>0$ such that
	\begin{equation*}\label{eq:Nehari-well}
		\mathcal J(v)\geq\mathcal J(u)
		+\epsilon_u\|v-u\|_{H^1_r}^2
	\end{equation*}
	whenever $v\in\cN$, $v\geq0$, and
	$\|v-u\|_{H^1_r}<\rho_u$.
\end{lemma}

\begin{proof}
	At a positive critical point $u$, differentiation of the Nehari
	constraint gives
	\[
	D\!\left[\langle\mathcal J'(v),v\rangle\right]_{v=u}[\phi]
	=
	\langle L_{+,u}u,\phi\rangle.
	\]
	Consequently,
	\begin{equation}\label{eq:Nehari-tangent}
		T_u\cN
		=
		\left\{\phi\in H^1_r(\RN):
		\langle L_{+,u}u,\phi\rangle=0\right\}.
	\end{equation}
	Moreover,
	$
	\langle L_{+,u}u,u\rangle
	=2(1-p)\|u\|_{H^1_r}^2<0.
	$
	We apply Lemma~\ref{lem:one-index-tangent-gap} to the quadratic form  \eqref{eq:full-linearized-form} below
	$
	\mathcal Q(\phi,\psi)
	:=
	\langle L_{+,u}\phi,\psi\rangle
	$
	on $H^1_r(\RN)$, 
	If $p<2$, then necessarily $\alpha=2$, and
	\[
	\begin{aligned}
		\mathcal Q[\phi]
		={}&
		\int_{\RN}\bigl(|\nabla\phi|^2+(2-p)\phi^2\bigr)
		-(p-1)\int_{\RN}(V-1)\phi^2
		-p\int_{\RN}
		\bigl(I_\alpha*(u^{p-1}\phi)\bigr)u^{p-1}\phi.
	\end{aligned}
	\]
	Here the first term is coercive.  By
	Lemma~\ref{lem:regularity-and-far-field}, multiplication by $V-1$ is
	compact from $H^1_r(\RN)$ to $L^2(\RN)$.  Moreover,
	Lemma~\ref{lem:regularity-and-far-field} and \eqref{eq:range} give
	$u^{p-1}\in L^{N/2}(\RN)$, so multiplication by $u^{p-1}$ is compact
	from $H^1_r(\RN)$ to $L^{2N/(N+2)}(\RN)$.
	
	If $p\geq2$, then
	\[
	\begin{aligned}
		\mathcal Q[\phi]
		={}&
		\|\phi\|_{H^1_r}^2-(p-1)\int_{\RN}V\phi^2 -p\int_{\RN}
		\bigl(I_\alpha*(u^{p-1}\phi)\bigr)u^{p-1}\phi.
	\end{aligned}
	\]
	In this case the coercive part is $\|\phi\|_{H^1_r}^2$, and
	Lemma~\ref{lem:regularity-and-far-field} shows that multiplication by $V$ is compact from
	$H^1_r(\RN)$ to $L^2(\RN)$.  Moreover,
	$u^{p-1}\in L^{2N/(\alpha+2)}(\RN)$, so multiplication by
	$u^{p-1}$ is compact from $H^1_r(\RN)$ to
	$L^{2N/(N+\alpha)}(\RN)$.
	
	In both cases these compactness statements follow from local Rellich
	compactness and the vanishing of the corresponding coefficient norm
	on the tail.  The Hardy--Littlewood--Sobolev inequality then shows that
	the nonlocal term is compact as well.
	
	By Theorem~\ref{thm:index-one}, the quadratic form $\mathcal Q$ has Morse
	index one and trivial nullspace.  The function $u$ is a negative direction
	for $\mathcal Q$, since
	$
	\mathcal Q[u]
	=
	\langle L_{+,u}u,u\rangle
	=
	2(1-p)\|u\|_{H^1_r}^2<0.
	$
	Moreover, \eqref{eq:Nehari-tangent} shows that the condition
	$\mathcal Q(u,\phi)=0$ is equivalent to $\phi\in T_u\cN$.  Thus all
	the hypotheses of Lemma~\ref{lem:one-index-tangent-gap} are satisfied
	with
	$
	X=H^1_r(\RN),
	\
	z=u.
	$
	Consequently, there exists $\gamma_u>0$ such that
	\begin{equation}\label{eq:tangent-gap}
		\langle L_{+,u}\phi,\phi\rangle
		\geq
		\gamma_u\|\phi\|_{H^1_r}^2
		\qquad\text{for every }\phi\in T_u\cN.
	\end{equation}
	
	The following construction is valid throughout the range
	\eqref{eq:range}, because the Nehari constraint is $C^1$.  We set
	$
	\Phi(\xi,s)
	=
	\left\langle
	\mathcal J'(u+\xi+su),u+\xi+su
	\right\rangle
	$ for $\xi\in T_u\cN$, $s\in\mathbb R$.
	Since $\Phi(0,0)=0$ and
	$
	\partial_s\Phi(0,0)
	=
	\langle L_{+,u}u,u\rangle<0,
	$
	the implicit function theorem gives a neighborhood
	$\mathcal U_u\subset T_u\cN$ of $0$ and a unique $C^1$ function
	$\kappa:\mathcal U_u\longrightarrow\mathbb R$ such that
	$
	\Phi(\xi,\kappa(\xi))=0
	\ \text{for every }\xi\in\mathcal U_u.
	$
	Equivalently,
	\[
	u+\xi+\kappa(\xi)u\in\cN
	\qquad\text{for every }\xi\in\mathcal U_u.
	\]
	Since
	$
	H^1_r(\RN)
	=
	T_u\cN\oplus\operatorname{span}\{u\},
	$
	every $v\in\cN$ sufficiently close to $u$ has the unique
	representation
	\begin{equation}\label{eq115}
		v=u+\xi+\kappa(\xi)u,
		\qquad
		\xi\in\mathcal U_u.
	\end{equation}
	Moreover, $\kappa(0)=0$.  Differentiating
	$\Phi(\xi,\kappa(\xi))=0$ at $\xi=0$ gives
	$
	D_\xi\Phi(0,0)
	+\partial_s\Phi(0,0)D\kappa(0)=0.
	$
	Since $D_\xi\Phi(0,0)$ vanishes on $T_u\cN$, while
	$\partial_s\Phi(0,0)=\langle L_{+,u}u,u\rangle\ne0$, one obtains
	$D\kappa(0)=0$.
	So we have
	$\kappa(\xi)=o(\|\xi\|_{H^1_r})$.
	
	Suppose first that $p\geq2$.  Then $\mathcal J$ is $C^2$ near $u$.
	Since
	$\langle L_{+,u}u,\xi\rangle=0$, Taylor's formula,
	\eqref{eq:tangent-gap}, \eqref{eq115}, and
	$\kappa(\xi)=o(\|\xi\|_{H^1_r})$ imply that, for small $\xi$,
	\begin{equation*}\label{eq116}
		\mathcal J(v)-\mathcal J(u)
		=
		\frac12\langle L_{+,u}\xi,\xi\rangle
		+o(\|\xi\|_{H^1_r}^2)\geq\frac{\gamma_u}{4}\|\xi\|_{H^1_r}^2\geq\frac{\gamma_u}{16}\|v-u\|_{H^1_r}^2.
	\end{equation*}
	Assume now that $\alpha=2$ and $p<2$.  A uniform second-order Taylor expansion in
	$H^1_r$ is unavailable because $u^{p-2}$ is singular where $u$ is small.
	We replace it by an exact identity and a uniform estimate on the tail.
	Let $v\geq0$.  Since
	\[
	\int_{\RN}
	\bigl(I_2*v^p\bigr)v^p-\int_{\RN}Au^p
	=
	2\int_{\RN}A\bigl(v^p-u^p\bigr)+\int_{\RN}\bigl(I_2*(v^p-u^p)\bigr)\bigl(v^p-u^p\bigr),
	\]
	and the equation for $u$ gives
	$\langle u,v-u\rangle_{H^1_r}=\int_{\RN}Au^{p-1}(v-u)$, we obtain
	\begin{align}\label{eq117}
		\mathcal J(v)-\mathcal J(u)
		&=\frac12\|v-u\|_{H^1_r}^2-\frac1p\int_{\RN}A
		\{v^p-u^p-pu^{p-1}(v-u)\} \nonumber\\
		&\quad -\frac1{2p}\int_{\RN}
		\bigl(I_2*(v^p-u^p)\bigr)(v^p-u^p).
	\end{align}
	Taylor's formula with integral remainder gives
	\[
	\begin{aligned}
		&v^p-u^p-pu^{p-1}(v-u) =
		p(p-1)(v-u)^2
		\int_0^1(1-\tau)
		\bigl((1-\tau)u+\tau v\bigr)^{p-2}\dd\tau.
	\end{aligned}
	\]
	Since $p-2<0$, $u>0$, and $v\geq0$, one has
	$
	\bigl((1-\tau)u+\tau v\bigr)^{p-2}
	\leq
	\bigl((1-\tau)u\bigr)^{p-2}.
	$
	Consequently,
	\[
	\begin{aligned}
		0
		&\leq v^p-u^p-pu^{p-1}(v-u) \leq
		p(p-1)u^{p-2}(v-u)^2
		\int_0^1(1-\tau)^{p-1}\dd\tau =(p-1)u^{p-2}(v-u)^2.
	\end{aligned}
	\]
	Using this and $Au^{p-2}=V$, we get
	\[
	\begin{aligned}
		&\frac12(v-u)^2
		-\frac A p\{v^p-u^p-pu^{p-1}(v-u)\} \geq
		\left(\frac12-\frac{p-1}{p}V(r)\right)(v-u)^2.
	\end{aligned}
	\]
	By 	Lemma~\ref{lem:regularity-and-far-field}, the coefficient on the right
	tends to $(2-p)/(2p)>0$ as $r\to\infty$.  Hence there exists $\varrho_*>0$
	such that for every $r\geq\varrho_*$  
	\begin{equation}\label{eq:tail-Bregman}
		\begin{aligned}
			&\frac12\bigl(v(r)-u(r)\bigr)^2
			-\frac{A(r)}p
			\left\{
			v(r)^p-u(r)^p-pu(r)^{p-1}\bigl(v(r)-u(r)\bigr)
			\right\} \geq
			\frac{2-p}{4p}\bigl(v(r)-u(r)\bigr)^2.
		\end{aligned}
	\end{equation}

	Suppose, to the contrary, that the quadratic lower bound fails.  Then
	there are $v_j\in\cN$, $v_j\geq0$, $v_j\ne u$, such that
	$v_j\to u$ in $H^1_r$ and
	\begin{equation}\label{eq:no-quadratic-well}
		\limsup_{j\to\infty}
		\frac{\mathcal J(v_j)-\mathcal J(u)}{\varepsilon_j^2}\leq0,\quad\text{where }\varepsilon_j=\|v_j-u\|_{H^1_r}.
	\end{equation}
	Set $
	\xi_j=\frac{v_j-u}{\varepsilon_j},
	$
	after passing to a subsequence,
	$\xi_j\rightharpoonup\xi$ in $H^1_r(\RN)$ and
	$\xi_j\to\xi$ in $L^s_{\rm loc}(\RN)$ for every
	$2\leq s<2N/(N-2)$.  Since
	$\langle\mathcal J'(v_j),v_j\rangle
	=\langle\mathcal J'(u),u\rangle=0$, differentiability of the Nehari
	constraint at $u$ gives
	\[
	0
	=
	\frac{
		\langle\mathcal J'(v_j),v_j\rangle
		-\langle\mathcal J'(u),u\rangle}{\varepsilon_j}
	=
	\langle L_{+,u}u,\xi_j\rangle+o(1).
	\]
	Passing to the limit yields
	\begin{equation}\label{eq:weak-tangent}
		\mathcal Q(u,\xi)=\langle L_{+,u}u,\xi\rangle=0.
	\end{equation}
	
	Since
	$v_j=u+\varepsilon_j\xi_j$, we have
	\[
	\frac{v_j^p-u^p}{\varepsilon_j}
	\to
	p u^{p-1}\xi
	\quad\text{in }L^{2N/(N+2)}(\RN).
	\]
	Here we use the first-order remainder estimate for $t\mapsto t^p$
	and the compactness of multiplication by $u^{p-1}$ established
	above.
	By this and the Hardy--Littlewood--Sobolev inequality, we get
	\begin{equation}\label{eq:nonlocal-quadratic-limit}
		\frac1{2p\varepsilon_j^2}
		\int_{\RN}\bigl(I_2*(v_j^p-u^p)\bigr)(v_j^p-u^p)
		\to
		\frac p2\int_{\RN}
		\bigl(I_2*(u^{p-1}\xi)\bigr)u^{p-1}\xi.
	\end{equation}
	We next identify the local quadratic limit.  Fix $\varrho>0$.  Since
	$u\geq u(\varrho)>0$ on $B_\varrho$, the estimate obtained above gives
	\[
	0\leq
	\frac{A}{p\varepsilon_j^2}
	\bigl\{v_j^p-u^p-pu^{p-1}(v_j-u)\bigr\}
	\leq C_\varrho\xi_j^2
	\qquad\text{on }B_\varrho.
	\]
	$\xi_j\to\xi$ in $L^2(B_\varrho)$, so
	$\xi_j^2\to\xi^2$ in $L^1(B_\varrho)$ and the preceding upper bounds are
	uniformly integrable.  Vitali's convergence theorem therefore gives
	\begin{equation}\label{eq:local-quadratic-limit}
		\frac1{p\varepsilon_j^2}\int_{B_\varrho}A
		\bigl\{v_j^p-u^p-pu^{p-1}(v_j-u)\bigr\}
		\to
		\frac{p-1}{2}\int_{B_\varrho}V\xi^2.
	\end{equation}
	We also set
	\[
	\delta_\nabla
	=\lim_{j\to\infty}\|\nabla\xi_j\|_2^2-\|\nabla\xi\|_2^2,
	\qquad
	\delta_0
	=\lim_{j\to\infty}\|\xi_j\|_2^2-\|\xi\|_2^2.
	\]
	Weak lower semicontinuity gives
	$\delta_\nabla,\delta_0\geq0$.  Since
	$\|\xi_j\|_{H^1_r}=1$,
	\begin{equation*}\label{eq:defect-mass}
		\delta_\nabla+\delta_0=1-\|\xi\|_{H^1_r}^2.
	\end{equation*}
	Fix $\varrho\geq\varrho_*$.  Applying \eqref{eq:tail-Bregman}   with
	$v=v_j$ gives
	\[
	\frac12(v_j-u)^2
	-\frac A p\{v_j^p-u^p-pu^{p-1}(v_j-u)\}
	\geq\frac{2-p}{4p}(v_j-u)^2
	\]
	on $\RN\setminus B_\varrho$.  Moreover, local strong $L^2$ convergence
	and the definition of $\delta_0$ give
	\[
	\int_{\RN\setminus B_\varrho}\xi_j^2
	=
	\|\xi_j\|_2^2-\int_{B_\varrho}\xi_j^2\to
	\|\xi\|_2^2+\delta_0-\int_{B_\varrho}\xi^2
	=
	\delta_0+\int_{\RN\setminus B_\varrho}\xi^2.
	\]
	Using these facts in the exact identity for
	$\mathcal J(v_j)-\mathcal J(u)$ obtained in \eqref{eq117}, together with
	\eqref{eq:nonlocal-quadratic-limit} and
	\eqref{eq:local-quadratic-limit}, yields
	\begin{align*}
		\liminf_{j\to\infty}
		\frac{\mathcal J(v_j)-\mathcal J(u)}{\varepsilon_j^2}
		\geq{}&\frac12(\|\nabla\xi\|_2^2+\delta_\nabla)
		+\frac12\int_{B_\varrho}\xi^2
		-\frac{p-1}{2}\int_{B_\varrho}V\xi^2\\
		&+\frac{2-p}{4p}
		\left(\delta_0+\int_{\RN\setminus B_\varrho}\xi^2\right)
		-\frac p2\int_{\RN}
		\bigl(I_2*(u^{p-1}\xi)\bigr)u^{p-1}\xi.
	\end{align*}
	By
	$V\in L^\infty(\RN)$,  letting $\varrho\to\infty$ therefore gives
	\begin{equation*}\label{eq:Bregman-defect}
		\liminf_{j\to\infty}
		\frac{\mathcal J(v_j)-\mathcal J(u)}{\varepsilon_j^2}
		\geq\frac12\langle L_{+,u}\xi,\xi\rangle
		+\frac12\delta_\nabla+\frac{2-p}{4p}\delta_0.
	\end{equation*}
	By \eqref{eq:weak-tangent} and \eqref{eq:tangent-gap}, the right-hand
	side is bounded below by
	\[
	\min\left\{
	\frac{\gamma_u}{2},\frac12,\frac{2-p}{4p}
	\right\}
	\bigl(\|\xi\|_{H^1_r}^2+\delta_\nabla+\delta_0\bigr)
	=
	\min\left\{
	\frac{\gamma_u}{2},\frac12,\frac{2-p}{4p}
	\right\}>0,
	\]
	contradicting \eqref{eq:no-quadratic-well}.  This contradiction
	completes the proof.
\end{proof}

\subsection{Two positive wells are impossible}
\label{subsection82}
Suppose that $u_0\ne u_1$ are two positive radial critical points.
They can be connected inside the nonnegative part of $\cN$.  Indeed,
for $w_t=(1-t)u_0+tu_1>0$, the path
\begin{equation}\label{eq:positive-Nehari-path}
	t\longmapsto
	\left(
	\|w_t\|_{H^1_r}^2/
	\int_{\RN}(I_\alpha*w_t^p)w_t^p
	\right)^{1/(2p-2)}w_t
\end{equation}
lies in $\cN$ and joins the two solutions.

Let
\begin{equation*}\label{eq:minimax-level}
	\mu=\inf_{\substack{\gamma\in C([0,1],\cN)\\
			\gamma(0)=u_0,\ \gamma(1)=u_1}}
	\max_{t\in[0,1]}\mathcal J(\gamma(t)).
\end{equation*}
The path in \eqref{eq:positive-Nehari-path} shows that this class is
nonempty.  The map $v\longmapsto |v|$ is continuous on
$H^1_r(\RN)$, and
$
\mathcal J(|v|)=\mathcal J(v),
\  |v|\in\cN\ \Longleftrightarrow\ v\in\cN.
$
Thus every admissible path can be replaced by a nonnegative admissible
path with the same energy profile, and the value of $\mu$ is unchanged
if the infimum is restricted to nonnegative paths.

Choose $0<\rho_i<\rho_{u_i}$ so small that
$u_{1-i}\notin\overline B_{\rho_i}(u_i)$, $i=0,1$.  Every
nonnegative path from $u_0$ to $u_1$ crosses both spheres
$\partial B_{\rho_i}(u_i)$.  Lemma~\ref{lem:Nehari-well} therefore gives
\begin{equation}\label{eq:minimax-above}
	\mu\geq\max_{i=0,1}
	\left\{\mathcal J(u_i)+\epsilon_{u_i}\rho_i^2\right\}
	>\max\{\mathcal J(u_0),\mathcal J(u_1)\}.
\end{equation}
Let $\varepsilon_j\downarrow0$.  By the definition of $\mu$, there are
nonnegative paths $\gamma_j\in C([0,1],\cN)$ joining $u_0$ to $u_1$
such that
\[
\max_{t\in[0,1]}\mathcal J(\gamma_j(t))
\leq\mu+\varepsilon_j^2.
\]
Applying the argument in the proof of
\cite[Theorem~2.8]{Willem} to $\mathcal J|_{\cN}$, with
$\epsilon=\varepsilon_j^2$ and $\delta=\varepsilon_j$, gives
$v_j\in\cN$ such that
\begin{equation}\label{eq:localized-minimax}
	|\mathcal J(v_j)-\mu|\leq2\varepsilon_j^2,
	\qquad
	\bigl\|D(\mathcal J|_{\cN})(v_j)\bigr\|_{(T_{v_j}\cN)^*}
	\leq8\varepsilon_j,
	\qquad
	\operatorname{dist}_{H^1_r}
	\bigl(v_j,\gamma_j([0,1])\bigr)\leq2\varepsilon_j.
\end{equation}
The sequence $(v_j)$ is a Palais--Smale sequence for
$\mathcal J|_{\cN}$ at level $\mu$.  After passing to a subsequence,
$v_j\to w$ strongly in $H^1_r(\RN)$.  By
\eqref{eq:localized-minimax}, we may choose $t_j\in[0,1]$ such that
\[
\|v_j-\gamma_j(t_j)\|_{H^1_r}\leq2\varepsilon_j.
\]
Every $\gamma_j(t_j)$ is nonnegative, and the nonnegative cone is closed
in $H^1_r$; hence $w\geq0$.  Moreover, $w\in\cN$, so $w\ne0$.  Since
$\cN$ is a natural constraint, $\mathcal J'(w)=0$.  Elliptic regularity
and the strong maximum principle
give $w>0$.  Also $\mathcal J(w)=\mu$, and
\eqref{eq:minimax-above} yields $w\notin\{u_0,u_1\}$.

Choose
\[
0<\rho<\min\left\{\rho_w,
\frac12\|w-u_0\|_{H^1_r},
\frac12\|w-u_1\|_{H^1_r}\right\}.
\]
Since $\gamma_j(t_j)\to w$, every sufficiently large path $\gamma_j$
enters $B_\rho(w)$, whereas both endpoints lie outside
$\overline B_\rho(w)$.  Thus $\gamma_j$ crosses
$\partial B_\rho(w)$ at a nonnegative point.  Lemma~\ref{lem:Nehari-well}
then gives
\[
\max_{t\in[0,1]}\mathcal J(\gamma_j(t))
\geq\mathcal J(w)+\epsilon_w\rho^2
=\mu+\epsilon_w\rho^2,
\]
contradicting
$\max_t\mathcal J(\gamma_j(t))\leq\mu+\varepsilon_j^2$ for large $j$.
Therefore two distinct positive radial critical points cannot exist.

\begin{proof}[Proof of Theorem~\ref{thm:main}]
	Existence in \eqref{eq:range} is standard.  The symmetry theorem moves
	every positive solution to a centered positive radial solution.
	Theorem~\ref{thm:index-one} applies to every such solution, and
	Lemma~\ref{lem:Nehari-well} makes every one of them a strict
	nonnegative Nehari well.  The preceding two-well contradiction gives
	at most one centered positive radial profile.  Translating back proves
	uniqueness up to translations.
\end{proof}

\section{Angular decomposition and full nondegeneracy}
\label{sec:full-nondegeneracy}

We now prove the full nondegeneracy of Theorem~\ref{thm:main} by removing the
radial restriction from the nondegeneracy statement.
In this section, \(L_{+,u}\) denotes the full-space linearization
\[
L_{+,u}:H^1(\RN)\longrightarrow H^{-1}(\RN),
\]
defined by
\begin{equation}\label{eq:full-linearized-form}
	\begin{aligned}
		\mathcal Q(\phi,\psi)
		 :=\bigl\langle L_{+,u}\phi,\psi\bigr\rangle_{H^{-1},H^1}
		=&\int_{\RN}
		\bigl(\nabla\phi\cdot\nabla\psi+\phi\psi
		-(p-1)V\phi\psi\bigr)\,\dd x\\
		&-p\int_{\RN}
		\bigl(I_\alpha*(u^{p-1}\phi)\bigr)u^{p-1}\psi\,\dd x,
		\qquad \phi,\psi\in H^1(\RN).
	\end{aligned}
\end{equation}
The Hardy--Littlewood--Sobolev and Sobolev inequalities show that
\(\mathcal Q\) is a continuous symmetric form on \(H^1(\RN)\).
Its restriction to \(H^1_r(\RN)\) is the radial linearization
introduced in Section~\ref{sec:half-potential} and used in
Theorem~\ref{thm:index-one}.

Let \(\{Y_{\ell,m}:1\leq m\leq d_\ell\}\) be a real orthonormal
basis of the \(d_\ell\)-dimensional space of spherical harmonics of
degree \(\ell\), so that
\[
-\Delta_{\mathbb S^{N-1}}Y_{\ell,m}
=\kappa_\ell Y_{\ell,m},
\qquad
\kappa_\ell=\ell(\ell+N-2).
\]
Note that for $r,s>0$, $\omega, \theta\in \mathbb S^{N-1}$,
\[
|r\theta-s\omega|=\sqrt{r^2+s^2-2rs\,\theta\cdot\omega},
\]
hence the Funk--Hecke formula \cite[Theorem 9.7.1]{AndrewsAskeyRoy1999}
gives
\begin{equation}\label{eq:angular-Riesz-kernel}
	\int_{\mathbb S^{N-1}}
	I_\alpha(r\theta-s\omega)Y_{\ell,m}(\omega)\,\dd\omega
	=\mathcal K_{\alpha,\ell}(r,s)Y_{\ell,m}(\theta),
\end{equation}
where \(\mathcal K_{\alpha,\ell}(r,s)\) depends only on
\(r\), \(s\) and \(\ell\), and is independent of \(m\).

\begin{lemma}\label{lem:angular-kernel-order}
	For every \(r,s>0\) and every \(\ell\geq0\),
	\begin{equation}\label{eq:angular-kernel-order}
		0<\mathcal K_{\alpha,\ell+1}(r,s)
		<\mathcal K_{\alpha,\ell}(r,s).
	\end{equation}
\end{lemma}

\begin{proof}
	Put
	\[
	\sigma=\frac{N-\alpha}{2},
	\qquad
	\nu=\frac{N-2}{2}.
	\]
	Recall that the Euler Gamma function is defined by
	\[
	\Gamma(\sigma)
	=
	\int_0^\infty y^{\sigma-1}e^{-y}\,\dd y,
	\qquad \sigma>0.
	\]
	For \(x\neq0\), making the change of variables
	\(y=t|x|^2\), we obtain
	\begin{equation}\label{eq128}
		|x|^{-2\sigma}
		=
		\frac1{\Gamma(\sigma)}
		\int_0^\infty
		t^{\sigma-1}e^{-t|x|^2}\,\dd t.
	\end{equation}
	
	Since \(2\sigma=N-\alpha\), it follows from
	\eqref{eq:Riesz-normalization} and
	\eqref{eq:angular-Riesz-kernel} that
	\[
	\begin{aligned}
		\mathcal K_{\alpha,\ell}(r,s)Y_{\ell,m}(\theta)
		&=
		\int_{\mathbb S^{N-1}}
		I_\alpha(r\theta-s\omega)
		Y_{\ell,m}(\omega)\,\dd\omega  \\
		&=
		\frac{c_{N,\alpha}}{\Gamma(\sigma)}
		\int_{\mathbb S^{N-1}}\int_0^\infty
		t^{\sigma-1}e^{-t|r\theta-s\omega|^2}
		Y_{\ell,m}(\omega)\,\dd t\,\dd\omega .
	\end{aligned}
	\]
	
	Since \(Y_{\ell,m}\) is bounded on \(\mathbb S^{N-1}\), we have
	\[
	\big|
	t^{\sigma-1}e^{-t|r\theta-s\omega|^2}
	Y_{\ell,m}(\omega)
	\big|
	\leq
	\|Y_{\ell,m}\|_{L^\infty(\mathbb S^{N-1})}
	t^{\sigma-1}e^{-t|r\theta-s\omega|^2}.
	\]
	The function on the right is integrable on
	\((0,\infty)\times\mathbb S^{N-1}\). Indeed, by \eqref{eq128}
	and Lemma~\ref{lem:log-conjugation}, we have
	\[
	\frac{c_{N,\alpha}}{\Gamma(\sigma)}
	\int_{\mathbb S^{N-1}}\int_0^\infty
	t^{\sigma-1}e^{-t|r\theta-s\omega|^2}
	\,\dd t\,\dd\omega
	=
	\int_{\mathbb S^{N-1}}
	I_\alpha(r\theta-s\omega)\,\dd\omega
	=
	K_\alpha(r,s)<\infty.
	\]
	Therefore, Fubini's theorem and
	\[
	|r\theta-s\omega|^2
	=
	r^2+s^2-2rs\,\theta\cdot\omega
	\]
	yield
	\[
	\begin{aligned}
		\mathcal K_{\alpha,\ell}(r,s)Y_{\ell,m}(\theta)
		={}&
		\frac{c_{N,\alpha}}{\Gamma(\sigma)}
		\int_0^\infty t^{\sigma-1}e^{-t(r^2+s^2)}
		\left(
		\int_{\mathbb S^{N-1}}
		e^{2trs\,\theta\cdot\omega}
		Y_{\ell,m}(\omega)\,\dd\omega
		\right)\dd t.
	\end{aligned}
	\]
	
	We now evaluate the inner integral. For every \(z>0\), we claim that
	\begin{equation}\label{eq:exp-FH}
		\int_{\mathbb S^{N-1}}
		e^{z\theta\cdot\omega}Y_{\ell,m}(\omega)\,\dd\omega
		=
		(2\pi)^{N/2}z^{-\nu}
		I_{\ell+\nu}(z)Y_{\ell,m}(\theta),
	\end{equation}
	where \(I_\eta\) denotes the modified Bessel function of the
	first kind of order \(\eta\).
	By the Funk--Hecke theorem \cite[Theorem 9.7.1]{AndrewsAskeyRoy1999},
	\[
	\int_{\mathbb S^{N-1}}
	e^{z\theta\cdot\omega}Y_{\ell,m}(\omega)\,\dd\omega
	=
	\Lambda_\ell(z)Y_{\ell,m}(\theta),
	\]
	where
	\[
	\Lambda_\ell(z)
	=
	|\mathbb S^{N-2}|
	\int_{-1}^1
	e^{zt}
	\frac{C_\ell^\nu(t)}{C_\ell^\nu(1)}
	(1-t^2)^{\nu-\frac12}\,\dd t,
	\]
	and \(C_\ell^\nu\) is the Gegenbauer polynomial of degree
	\(\ell\) and parameter \(\nu\).
	Now apply Gegenbauer's formula
	\cite[(4.7.7)]{AndrewsAskeyRoy1999},
	\cite[Eq.~7.321]{GradshteynRyzhik}, stated for imaginary argument:
	\[
	\int_{-1}^1
	e^{zt}C_\ell^\nu(t)
	(1-t^2)^{\nu-\frac12}\,\dd t
	=
	\frac{\pi\,2^{1-2\nu}\Gamma(\ell+2\nu)}
	{\ell!\,\Gamma(\nu)}
	\left(\frac z2\right)^{-\nu}
	I_{\ell+\nu}(z).
	\]
	Using
	\[
	|\mathbb S^{N-2}|
	=
	\frac{2\pi^{\nu+\frac12}}{\Gamma(\nu+\frac12)},
	\qquad
	C_\ell^\nu(1)
	=
	\frac{\Gamma(\ell+2\nu)}{\ell!\,\Gamma(2\nu)}
	\]
	(see \cite[Eq.~8.937.4]{GradshteynRyzhik} for the value of \(C_\ell^\nu(1)\)),
	we obtain
	\[
	\begin{aligned}
		\Lambda_\ell(z)
		&=
		\frac{2\pi^{\nu+\frac12}}{\Gamma(\nu+\frac12)}
		\frac{\ell!\,\Gamma(2\nu)}{\Gamma(\ell+2\nu)}
		\frac{\pi\,2^{1-2\nu}\Gamma(\ell+2\nu)}
		{\ell!\,\Gamma(\nu)}
		\left(\frac z2\right)^{-\nu}
		I_{\ell+\nu}(z) \\
		&=
		\frac{2^{2-2\nu}\pi^{\nu+\frac32}\Gamma(2\nu)}
		{\Gamma(\nu+\frac12)\Gamma(\nu)}
		\left(\frac z2\right)^{-\nu}
		I_{\ell+\nu}(z).
	\end{aligned}
	\]
	By Legendre's duplication formula
	\cite[Eq.~8.335]{GradshteynRyzhik},
	\[
	\Gamma(2\nu)
	=
	\frac{2^{2\nu-1}}{\sqrt\pi}
	\Gamma(\nu)\Gamma\!\left(\nu+\frac12\right),
	\]
	the constant simplifies to
	\[
	\frac{2^{2-2\nu}\pi^{\nu+\frac32}\Gamma(2\nu)}
	{\Gamma(\nu+\frac12)\Gamma(\nu)}
	=
	2\pi^{\nu+1}.
	\]
	Therefore
	\[
	\Lambda_\ell(z)
	=
	2\pi^{\nu+1}
	\left(\frac z2\right)^{-\nu}
	I_{\ell+\nu}(z)
	=
	(2\pi)^{\nu+1}z^{-\nu}
	I_{\ell+\nu}(z).
	\]
	Since \(\nu+1=N/2\), this proves \eqref{eq:exp-FH}.
	
	Applying \eqref{eq:exp-FH} with \(z=2trs\) in the preceding
	integral gives
	\begin{equation}\label{eq:angular-kernel-Bessel}
		\mathcal K_{\alpha,\ell}(r,s)
		=
		\frac{c_{N,\alpha}(2\pi)^{N/2}}{\Gamma(\sigma)}
		\int_0^\infty
		t^{\sigma-1}e^{-t(r^2+s^2)}
		(2trs)^{-\nu}
		I_{\ell+\nu}(2trs)\,\dd t.
	\end{equation}
	
	By Soni's inequality \cite{Soni1965},
	\[
	0<I_{\eta+1}(z)<I_\eta(z),
	\qquad z>0,\quad \eta>-\frac12.
	\]
	Since \(\nu=(N-2)/2\geq0\) and \(\ell\geq0\), we have
	\(\ell+\nu>-\frac12\). Hence, for every \(t>0\),
	\[
	0<
	I_{\ell+1+\nu}(2trs)
	<
	I_{\ell+\nu}(2trs).
	\]
	All the remaining factors in
	\eqref{eq:angular-kernel-Bessel} are strictly positive and independent
	of \(\ell\). Consequently,
	\[
	0<\mathcal K_{\alpha,\ell+1}(r,s)
	<\mathcal K_{\alpha,\ell}(r,s),
	\]
	which proves \eqref{eq:angular-kernel-order}.
\end{proof}

For \(\ell\geq0\), set
\[
\mathscr H_\ell
=
\left\{
h\in H^1_{\mathrm{loc}}(0,\infty):
\int_0^\infty
\left(
|h'|^2+h^2+\frac{\kappa_\ell}{r^2}h^2
\right)r^{N-1}\,\dd r<\infty
\right\}.
\]

\begin{proof}[Completing the proof of Theorem~\ref{thm:main}]
	For
	\[
	\phi(r\theta)
	=\sum_{\ell=0}^\infty\sum_{m=1}^{d_\ell}
	h_{\ell,m}(r)Y_{\ell,m}(\theta),
	\]
	the orthonormality of the spherical harmonics and
	\(-\Delta_{\mathbb S^{N-1}}Y_{\ell,m}
	=\kappa_\ell Y_{\ell,m}\) give
	\[
	\begin{aligned}
		\int_{\mathbb S^{N-1}}
		Y_{\ell,m}Y_{\ell',m'}\,\dd\theta
		&=\delta_{\ell\ell'}\delta_{mm'},\qquad
		\int_{\mathbb S^{N-1}}
		\nabla_{\mathbb S^{N-1}}Y_{\ell,m}\cdot
		\nabla_{\mathbb S^{N-1}}Y_{\ell',m'}\,\dd\theta
		=\kappa_\ell\delta_{\ell\ell'}\delta_{mm'}.
	\end{aligned}
	\]
	Since
	$
	|\nabla\phi|^2
	=|\partial_r\phi|^2
	+\frac1{r^2}|\nabla_{\mathbb S^{N-1}}\phi|^2,
	$
	and \(V\) is radial, the two orthogonality relations above give,
	for the local part of \(\mathcal Q\),
	\[
	\int_{\RN}\Bigl(|\nabla\phi|^2+\phi^2-(p-1)V\phi^2\Bigr)\dd x
	=
	\sum_{\ell=0}^\infty\sum_{m=1}^{d_\ell}
	\int_0^\infty
	\Bigl[
	|h_{\ell,m}'|^2
	+\Bigl(1+\frac{\kappa_\ell}{r^2}-(p-1)V(r)\Bigr)h_{\ell,m}^2
	\Bigr]r^{N-1}\,\dd r .
	\]
	For the nonlocal part,
	\eqref{eq:angular-Riesz-kernel} and the preceding orthogonality give
	\[
	\begin{aligned}
		\int_{\mathbb S^{N-1}}\int_{\mathbb S^{N-1}}
		I_\alpha(r\theta-s\omega)
		Y_{\ell,m}(\theta)Y_{\ell',m'}(\omega)
		\,\dd\omega\,\dd\theta 
		&=
		\mathcal K_{\alpha,\ell'}(r,s)
		\int_{\mathbb S^{N-1}}
		Y_{\ell,m}(\theta)Y_{\ell',m'}(\theta)\,\dd\theta \\
		&=
		\mathcal K_{\alpha,\ell}(r,s)
		\delta_{\ell\ell'}\delta_{mm'}.
	\end{aligned}
	\]
	Combining the two computations, the two angular indices
	\((\ell,m)\) and \((\ell',m')\) agree in every nonzero term, and
	\begin{equation}\label{eq:partial-wave-decomposition}
		\mathcal Q[\phi]
		=\sum_{\ell=0}^\infty\sum_{m=1}^{d_\ell}
		\mathcal Q_\ell[h_{\ell,m}],
	\end{equation}
	where
	\begin{equation*}\label{eq:partial-wave-form}
		\begin{aligned}
			\mathcal Q_\ell[h]
			={}&\int_0^\infty
			\left[
			|h'|^2+
			\left(1+\frac{\kappa_\ell}{r^2}-(p-1)V(r)\right)h^2
			\right]r^{N-1}\,\dd r\\
			&-p\int_0^\infty\int_0^\infty
			\mathcal K_{\alpha,\ell}(r,s)
			u(r)^{p-1}h(r)u(s)^{p-1}h(s)
			(rs)^{N-1}\,\dd s\,\dd r.
		\end{aligned}
	\end{equation*}
	for \(h\in\mathscr H_\ell\).  The decomposition is first obtained
	for finite spherical-harmonic sums of smooth functions and then
	extended to \(H^1(\RN)\) by density. We denote by
	\(L_{+,u}^{(\ell)}\) the operator associated with
	\(\mathcal Q_\ell\) on
	\(L^2((0,\infty),r^{N-1}\dd r)\).

	The first eigenvalue of
	\(-\Delta_{\mathbb S^{N-1}}\) is
	\(\kappa_0=0\), and its eigenspace is the one-dimensional
	space of constant functions.  We choose its normalized basis element as
	$
	Y_{0,1}(\theta)=|\mathbb S^{N-1}|^{-1/2}.
	$
	Therefore the term with \(\ell=0\) in the spherical-harmonic
	expansion is
	\[
	h_{0,1}(r)Y_{0,1}(\theta)
	=|\mathbb S^{N-1}|^{-1/2}h_{0,1}(r),
	\]
	which is radial.  By
	\eqref{eq:partial-wave-decomposition}, we have
	$
	\mathcal Q\bigl[hY_{0,1}\bigr]=\mathcal Q_0[h].
	$
	Hence Theorem~\ref{thm:index-one} gives
	\begin{equation}\label{eq:zeroth-channel}
		i_-(\mathcal Q_0)=1,
		\qquad
		\ker\mathcal Q_0=\{0\}.
	\end{equation}

	We next consider \(\ell=1\).  Differentiating
	\eqref{eq:choquard} with respect to \(x_j\) gives
	\begin{equation}\label{eq:translation-zero-modes}
		L_{+,u}(\partial_{x_j}u)=0,
		\qquad 1\leq j\leq N.
	\end{equation}
	Since
	\(\partial_{x_j}u=u'(r)x_j/r\), the function
	$
	w(r)=-u'(r)>0
	$
	satisfies \(L_{+,u}^{(1)}w=0\).  Here positivity follows from
	Theorem~\ref{thm:symmetry}.  The regularity at the origin and
	\eqref{eq:u-exponential-decay} give \(w\in\mathscr H_1\) when
	\(p\geq2\); in the Newtonian subquadratic case, the same conclusion
	follows from \eqref{eq:u-subquadratic-derivatives} and
	\(\beta_p>N\), which is a consequence of \eqref{eq:range}.
	
	Let \(h\in C_c^\infty(0,\infty)\).  Testing the equation
	\(L_{+,u}^{(1)}w=0\) with \(h^2/w\), and subtracting the resulting
	identity from \(\mathcal Q_1[h]\), gives
	\begin{equation}\label{eq:first-channel-ground-state}
		\begin{aligned}
			\mathcal Q_1[h]
			={}&\int_0^\infty
			w(r)^2\left|\left(\frac{h}{w}\right)'(r)\right|^2
			r^{N-1}\,\dd r\\
			&+\frac p2\int_0^\infty\int_0^\infty
			\mathcal K_{\alpha,1}(r,s)
			u(r)^{p-1}w(r)u(s)^{p-1}w(s)
			\left(\frac{h(r)}{w(r)}-\frac{h(s)}{w(s)}\right)^2
			(rs)^{N-1}\,\dd s\,\dd r.
		\end{aligned}
	\end{equation}
	Both terms on the right are nonnegative by
	Lemma~\ref{lem:angular-kernel-order}.  Since
	\(C_c^\infty(0,\infty)\) is dense in \(\mathscr H_1\),
	approximation in the \(\mathscr H_1\)-norm and continuity of
	\(\mathcal Q_1\) show that \(\mathcal Q_1\geq0\) on
	\(\mathscr H_1\).  If \(h\in\ker L_{+,u}^{(1)}\), then
	\(\mathcal Q_1[h]=0\).  Applying
	\eqref{eq:first-channel-ground-state} to an approximating sequence,
	we obtain \((h/w)'=0\) on every compact subinterval of
	\((0,\infty)\).  Hence \(h\) is a constant multiple of \(w\).
	Conversely, \(L_{+,u}^{(1)}w=0\).  Therefore
	\begin{equation}\label{eq:first-channel-kernel}
		\mathcal Q_1[h]\geq0
		\quad\text{for every }h\in\mathscr H_1,
		\qquad
		\mathcal Q_1[h]=0
		\quad\Longleftrightarrow\quad
		h\in\operatorname{span}\{w\}.
	\end{equation}
	
	Finally, let \(\ell\geq2\).  Since
	\(\mathcal K_{\alpha,\ell}>0\), replacing \(h\) by \(|h|\) does not
	increase \(\mathcal Q_\ell\).  Since
	\(\mathscr H_\ell\subset\mathscr H_1\),
	\eqref{eq:angular-kernel-order} and
	\eqref{eq:first-channel-kernel} give, for every
	\(0\ne h\in\mathscr H_\ell\),
	\begin{equation}\label{eq:higher-channel-positive}
		\begin{aligned}
			\mathcal Q_\ell[h]
			\geq\mathcal Q_\ell[|h|] &=\mathcal Q_1[|h|]
			+(\kappa_\ell-\kappa_1)
			\int_0^\infty |h(r)|^2r^{N-3}\,\dd r\\
			&\quad+p\int_0^\infty\int_0^\infty
			\bigl(\mathcal K_{\alpha,1}
			-\mathcal K_{\alpha,\ell}\bigr)(r,s)
			u(r)^{p-1}u(s)^{p-1}|h(r)||h(s)|
			(rs)^{N-1}\,\dd s\,\dd r>0.
		\end{aligned}
	\end{equation}
	
	Combining \eqref{eq:partial-wave-decomposition},
	\eqref{eq:zeroth-channel},
	\eqref{eq:first-channel-kernel}, and
	\eqref{eq:higher-channel-positive}, and using that
	\(\mathcal Q\) diagonalizes in the spherical-harmonic
	decomposition so that indices and nullities add over \(\ell\),
	we find that \(L_{+,u}\) has Morse index one on \(H^1(\RN)\).
	Moreover, the orthogonal spherical-harmonic projections commute with
	\(L_{+,u}\).  Hence every full-space zero mode lies in the first
	channel and has radial coefficient proportional to \(u'\).
	Since the first spherical harmonics are spanned by
	\(x_j/r\), \(1\leq j\leq N\),
	\eqref{eq:translation-zero-modes} proves
	the two identities in \eqref{eq:full-nondegeneracy} of Theorem~\ref{thm:main}.
\end{proof}

\appendix
\renewcommand{\theequation}{\thesection.\arabic{equation}}
\section{Operator-theoretic facts used in the proof}
\label{app:operator-facts}

We collect here the standard operator-theoretic results used above.
They are stated for real Hilbert spaces; the complex versions are
identical.

\begin{lemma}\label{lem:Schur-inertia}
	Let \(\mathcal H=\mathcal H_1\oplus\mathcal H_2\), let \(A\) be self-adjoint
	with domain \(D(A)\subset\mathcal H_1\), let \(C\in\mathcal B(\mathcal H_2)\)
	be self-adjoint and invertible, and let
	\(B\in\mathcal B(\mathcal H_2,\mathcal H_1)\).  Then
	\[
	H=\begin{pmatrix}A&B\\ B^*&C\end{pmatrix},
	\qquad
	D(H)=D(A)\oplus\mathcal H_2,
	\]
	is self-adjoint, so is its Schur complement \(S=A-BC^{-1}B^*\), and
	\begin{itemize}
		\item[(i)] \(i_-(H)=i_-(S)+i_-(C)\) and \(i_0(H)=i_0(S)+i_0(C)\);
		in particular \(i_0(H)=i_0(S)\), and \(i_-(H)=i_-(S)\) if \(C\geq0\);
		\item[(ii)] \(\ker H=\{(u,-C^{-1}B^*u):u\in\ker S\}\), so that
		\(\dim\ker H=\dim\ker S\);
		\item[(iii)] \(H\) is Fredholm if and only if \(S\) is, and then the two
		have the same Fredholm index.
	\end{itemize}
\end{lemma}

\begin{proof}
	Put
	\[
	L=\begin{pmatrix}I&-BC^{-1}\\ 0&I\end{pmatrix},
	\qquad
	L^{-1}=\begin{pmatrix}I&BC^{-1}\\ 0&I\end{pmatrix},
	\qquad
	L^*=\begin{pmatrix}I&0\\ -C^{-1}B^*&I\end{pmatrix}.
	\]
	Since \(B\) and \(C^{-1}\) are bounded, \(L\) is bounded with bounded
	inverse, and a direct computation gives
	\(LHL^*=\operatorname{diag}(S,C)\).  As \(L^*\) maps
	\(D(S)\oplus\mathcal H_2=D(A)\oplus\mathcal H_2\) bijectively onto \(D(H)\),
	the two quadratic forms correspond under a bounded invertible map; by the
	operator form of Sylvester's law of inertia
	\cite[Lemma~2.3]{BerkolaikoKuchment2022} the two operators have the same
	negative and zero subspace dimensions, which is (i).  Item (ii) follows by
	solving the system \(Au+Bv=0\), \(B^*u+Cv=0\) (the same relation is obtained,
	under more general hypotheses, in \cite[Proposition~2.6]{Gerhat2024}).  For (iii) we use
	\(H=L^{-1}\operatorname{diag}(S,C)(L^{-1})^*\): since \(L^{-1}\) is boundedly
	invertible and \(C\) is boundedly invertible, hence Fredholm of index zero,
	\(H\) is Fredholm if and only if \(S\) is, with the same index.  This is the
	Schur complement factorization of block operator matrices
	\cite[Section~1.6]{Tretter2008}; the inertia formula goes back to
	\cite{Haynsworth1968}, and the general transfer of semi-Fredholmness between
	an operator matrix and its Schur complement is
	\cite[Theorems~2.7 and~2.8]{Gerhat2024}.
\end{proof}

\begin{lemma}
	\label{lem:product-intertwining}
	Let \(B:D(B)\subset\mathcal H_1\to\mathcal H_2\) be densely defined and closed.  Then, for every \(s\geq0\),
	\[
	B^*e^{-sBB^*}=e^{-sB^*B}B^*
	\qquad\text{on }D(B^*).
	\]
\end{lemma}

\begin{proof}
	By \cite[Theorem~4.10 and Corollary~4.11]{Teschl2014},
	\(B=U|B|\) with \(|B|=(B^*B)^{1/2}\), \(B^*=|B|U^*\), and \(U\)
	maps \((\ker B)^\perp\) unitarily onto
	\(\overline{\operatorname{Ran}B}=(\ker B^*)^\perp\).  Hence
	\[
	BB^*=U|B|^2U^*=U(B^*B)U^*,
	\qquad\text{that is,}\qquad
	U(B^*B)\subseteq(BB^*)U .
	\]
	Since \(U\) is bounded, \cite[Problem~4.6]{Teschl2014} gives
	\(U\Phi(B^*B)\subseteq\Phi(BB^*)U\) for every Borel function
	\(\Phi\); for bounded \(\Phi\) this is an equality of bounded
	operators, \(U\Phi(B^*B)=\Phi(BB^*)U\).  Put
	\(Q=UU^*\); then \(Q^\perp=\ker B^*\) and \(BB^*=0\) on
	\(Q^\perp\), so that
	\[
	\Phi(BB^*)=U\Phi(B^*B)U^*+\Phi(0)(I-Q).
	\]
	Since \(U^*(I-Q)=U^*-U^*UU^*=0\), \(|B|U^*U=|B|\) (as \(U^*U\)
	is the projection onto \((\ker B)^\perp\) and \(|B|\) vanishes on
	\(\ker B\)), and \(|B|=(B^*B)^{1/2}\) commutes with
	\(\Phi(B^*B)\),
	\[
	\begin{aligned}
		B^*\Phi(BB^*)
		&=|B|U^*\bigl(U\Phi(B^*B)U^*+\Phi(0)(I-Q)\bigr)\\
		&=|B|U^*U\Phi(B^*B)U^*
		=|B|\Phi(B^*B)U^*
		=\Phi(B^*B)B^* .
	\end{aligned}
	\]
	All identities are understood on \(D(B^*)\).  Indeed, for
	\(v\in D(B^*)\) one has \(U^*v\in D(|B|)\) by the polar
	decomposition, while \((I-Q)v\in\ker B^*\) is annihilated by
	\(B^*\), so that every composition above is defined.
	The choice \(\Phi(\lambda)=e^{-s\lambda}\), \(s\geq0\), is the
	claim; the choice \(\Phi(\lambda)=(1+t\lambda)^{-1}\), \(t>0\),
	recovers the resolvent intertwining \cite[(8.50)]{Teschl2014}.
\end{proof}

\begin{lemma}
	\label{lem:unitary-calculus}
	Let \(A\) be self-adjoint on \(\mathcal H_1\), and let
	\(U:\mathcal H_1\to\mathcal H_2\) be unitary.  Then
	\[
	\widetilde A=UAU^{-1},
	\qquad
	D(\widetilde A)=UD(A),
	\]
	is self-adjoint, and the spectral calculus satisfies
	\[
	\Phi(\widetilde A)=U\Phi(A)U^{-1}
	\]
	for every bounded Borel function \(\Phi\).  In particular, if
	\(A\) is bounded from below, then \(e^{-sA}\) and
	\(e^{-s\widetilde A}\) are bounded operators for every \(s\geq0\)
	and
	\[
	e^{-s\widetilde A}=Ue^{-sA}U^{-1},
	\qquad s\geq0 .
	\]
\end{lemma}

This is the unitary invariance of the Borel calculus: \(U\) is bounded
and \(UA\subseteq\widetilde AU\) with \(\widetilde A=UAU^{-1}\), so
\cite[Problem~4.6]{Teschl2014} gives
\(U\Phi(A)\subseteq\Phi(\widetilde A)U\) for every bounded Borel
function \(\Phi\), that is, \(\Phi(\widetilde A)=U\Phi(A)U^{-1}\).
The self-adjointness of \(\widetilde A\) and the relation
\(D(\widetilde A)=UD(A)\) are \cite[(3.34)]{Teschl2014}.

\begin{lemma}
	\label{lem:kernel-projection}
	Let \(A\) be self-adjoint and bounded from below, suppose that zero is
	an isolated spectral point, and denote by \(\Pi_0\) the orthogonal
	projection onto \(\ker A\).  If
	\(h\perp\operatorname{Ran}\mathbf 1_{(-\infty,0)}(A)\), then
	\[
	e^{-sA}h\to \Pi_0h
	\qquad\text{in }\mathcal H
	\quad\text{as }s\to+\infty.
	\]
\end{lemma}

This is immediate from the spectral theorem,
\cite[Theorem~3.6]{Teschl2014}, applied to the bounded Borel function
\(e^{-s\lambda}\).

\begin{lemma}\label{lem:one-index-tangent-gap}
	Let $X$ be a Hilbert space and let $\mathcal Q$ be a continuous
	symmetric bilinear form on $X$.  Write
	$\mathcal Q[\xi]=\mathcal Q(\xi,\xi)$.  Assume that $\mathcal Q$ is
	a compact perturbation of a coercive form, has Morse index one, and
	that
	\[
	\mathcal Q(\xi,\eta)=0\quad\text{for every }\eta\in X
	\quad\Longrightarrow\quad\xi=0.
	\]
	If $z\in X$ satisfies $\mathcal Q[z]<0$, then there exists
	$c_z>0$ such that
	$
	\mathcal Q[\xi]\geq c_z\|\xi\|_X^2
	$
	for every $\xi\in X$ satisfying $\mathcal Q(z,\xi)=0$.
\end{lemma}

\begin{proof}
	Set
	\[
	X_z=\{\xi\in X:\mathcal Q(z,\xi)=0\}.
	\]
	Since $\mathcal Q[z]<0$, one has
	$X=\operatorname{span}\{z\}\oplus X_z$.  The form $\mathcal Q$ is
	nonnegative on $X_z$: otherwise, $z$ and a negative vector in $X_z$
	would span a two-dimensional negative subspace.
	
	Choose a coercive form $\mathcal Q_0$ and a form
	$\mathcal Q_{\rm c}$ represented by a compact self-adjoint operator
	on $X$ such that
	$\mathcal Q=\mathcal Q_0+\mathcal Q_{\rm c}$.  If the asserted
	estimate were false, there would be $\xi_j\in X_z$ such that
	\[
	\|\xi_j\|_X=1,
	\qquad
	\mathcal Q[\xi_j]\to0.
	\]
	After passing to a subsequence, $\xi_j\rightharpoonup\xi$ in $X$.
	The subspace $X_z$ is weakly closed, so $\xi\in X_z$.  Compactness
	and weak lower semicontinuity give
	\[
	\mathcal Q_{\rm c}[\xi_j]\to
	\mathcal Q_{\rm c}[\xi],
	\qquad
	\mathcal Q[\xi]\leq
	\liminf_{j\to\infty}\mathcal Q[\xi_j]=0.
	\]
	Thus $\mathcal Q[\xi]=0$.  Since $\mathcal Q$ is nonnegative on
	$X_z$, $\xi$ is a minimizer of $\mathcal Q$ on $X_z$, and the
	first-order condition at a minimizer gives
	$\mathcal Q(\xi,\eta)=0$ for every $\eta\in X_z$.  Also
	$\mathcal Q(\xi,z)=0$.  The direct-sum decomposition then implies
	$\mathcal Q(\xi,\eta)=0$ for every $\eta\in X$, and hence $\xi=0$.
	It follows that $\mathcal Q_{\rm c}[\xi_j]\to0$ and therefore
	$\mathcal Q_0[\xi_j]\to0$, contradicting the coercivity of
	$\mathcal Q_0$ and $\|\xi_j\|_X=1$.
\end{proof}

\section{Technical proofs for the Newtonian cumulative inequality}
\label{app:newtonian}

We prove here the two results stated in Subsection~\ref{subsec:newtonian}.

\begin{proof}[Proof of Lemma~\ref{lem:local-pohozaev}]
	Let
	\(r>0\). 
	Multiplying
	\(-\Delta u+u=Au^{p-1}\) by \(x\cdot\nabla u\), integrating over
	\(B_r\), and using \(u^{p-1}\nabla u=\frac1p\nabla(u^p)\) together with
	\(\operatorname{div}(Axu^p)
	=(x\cdot\nabla A)u^p+NAu^p+A(x\cdot\nabla u^p)\), we obtain
	\begin{equation}\label{eq:ball-pohozaev-u}
		\begin{aligned}
			-\frac{N-2}{2}\int_{B_r}|\nabla u|^2
			-\frac N2\int_{B_r}u^2
			&+\frac12|\mathbb S^{N-1}|r^N(u^2-u'^2)\\
			=
			&\frac1p|\mathbb S^{N-1}|r^NAu^p
			-\frac1p\int_{B_r}u^p(x\cdot\nabla A)
			-\frac Np\int_{B_r}Au^p,
		\end{aligned}
	\end{equation}
	while multiplying \(-\Delta A=u^p\) by \(x\cdot\nabla A\) and
	integrating over \(B_r\) gives
	\begin{equation}\label{eq:ball-pohozaev-A}
		-\frac{N-2}{2}\int_{B_r}|\nabla A|^2
		-\frac12|\mathbb S^{N-1}|r^NA'^2
		=\int_{B_r}u^p(x\cdot\nabla A).
	\end{equation}
	
	Inserting \eqref{eq:ball-pohozaev-A} into
	\eqref{eq:ball-pohozaev-u},
	\[
	\begin{aligned}
		&-\frac{N-2}{2}\int_{B_r}|\nabla u|^2
		-\frac N2\int_{B_r}u^2
		+\frac12|\mathbb S^{N-1}|r^N(u^2-u'^2)\\
		&\qquad=
		\frac1p|\mathbb S^{N-1}|r^NAu^p
		+\frac{N-2}{2p}\int_{B_r}|\nabla A|^2
		+\frac1{2p}|\mathbb S^{N-1}|r^NA'^2
		-\frac Np\int_{B_r}Au^p .
	\end{aligned}
	\]
	
	Testing \(-\Delta A=u^p\) by \(A\), and \(-\Delta u+u=Au^{p-1}\) by
	\(u\), gives
	\[
	\int_{B_r}|\nabla A|^2=\int_{B_r}Au^p
	+|\mathbb S^{N-1}|r^{N-1}AA',
	\qquad
	\int_{B_r}|\nabla u|^2=\int_{B_r}Au^p-\int_{B_r}u^2
	+|\mathbb S^{N-1}|r^{N-1}uu' .
	\]
	Substituting these two identities into the preceding relation and
	collecting the volume terms, with
	\(-\frac{N-2}2+\frac{N+2}{2p}
	=\frac{N+2-(N-2)p}{2p}=1/V_*\),
	\[
	\frac1{V_*}\int_{B_r}Au^p-\int_{B_r}u^2+|\mathbb S^{N-1}|\mathcal P(r)=0,
	\]
	where the boundary term is
	\[
	\begin{aligned}
		\mathcal P(r)={}&
		-\frac{N-2}{2}r^{N-1}uu'
		+\frac12r^N(u^2-u'^2)
		-\frac1pr^NAu^p\\
		&-\frac{N-2}{2p}r^{N-1}AA'
		-\frac1{2p}r^NA'^2 .
	\end{aligned}
	\]
	Since \(Au^p=Vu^2\), this says
	\(\int_{B_r}u^2\bigl(1-\frac V{V_*}\bigr)\dd x=|\mathbb S^{N-1}|\mathcal P(r)\);
	dividing by \(|\mathbb S^{N-1}|\) and writing the integrals radially,
	\begin{equation}\label{eq:ball-cumulative}
		\int_0^r s^{N-1}u^2\Bigl(1-\frac{V(s)}{V_*}\Bigr)\dd s
		=\mathcal P(r).
	\end{equation}
	
	Differentiating \(-\frac{ru'}u\) and
	\(-\frac{rA'}A\) and using the equations for \(u\) and \(A\), we obtain the
	elementary identities
	\begin{align} 
		r\Bigl(-\frac{ru'}u\Bigr)'
		&=\Bigl(\frac{ru'}u\Bigr)^2+(N-2)\frac{ru'}u+r^2(V-1),\label{eq:log-identities1}\\
		r\Bigl(-\frac{rA'}A\Bigr)'
		&=\Bigl(\frac{rA'}A\Bigr)^2+(N-2)\frac{rA'}A+\frac{r^2u^p}A.\label{eq:log-identities2}
	\end{align}
	Substituting \eqref{eq:log-identities1} and \eqref{eq:log-identities2}  into the definition
	\eqref{eq:E-def} of \(E\) gives
	\[
	E(r)=(N-2)\frac{ru'}u-r^2+r^2\Bigl(\frac{u'}u\Bigr)^2
	+\frac2pr^2V+\frac{N-2}p\frac{rAA'}{u^2}+\frac1p\frac{r^2A'^2}{u^2} = -2r^{2-N}u^{-2}\mathcal P(r).
	\]
	Hence,
	\eqref{eq:ball-cumulative} is exactly \eqref{eq:local-pohozaev}.
\end{proof}

\begin{proof}[Proof of Proposition~\ref{prop:E-positive}]
	Set
	\begin{equation*}\label{eq:H-def}
		H=\frac{A-u^{2-p}}u=(V-1)u^{1-p}.
	\end{equation*}
	Then
	\begin{equation}\label{eq:V-H}
		V-1=u^{p-1}H,
		\qquad
		-\Delta u=u(V-1)=u^pH.
	\end{equation}
	We shall prove that the $u^2$-weighted running average $\bar V$ of
	$V$, defined in \eqref{eq:V-bar-def}, is strictly decreasing and
	stays strictly above its limit $V_*$; by
	Lemma~\ref{lem:local-pohozaev} this is equivalent to $E(r)>0$.
	
	\medskip
	\noindent
	\emph{Step 1. An algebraic barrier.}
	
	For $s>0$, define $\theta(s)>0$ by
	\begin{equation}\label{eq:barrier-theta-equation}
		\theta(s)\left(\theta(s)+(2-p)s^{1-p}\right)=1 ;
	\end{equation}
	explicitly,
	\[
	\theta(s)=\frac12\left(
	\sqrt{(2-p)^2s^{2-2p}+4}
	-(2-p)s^{1-p}\right)>0.
	\]
	Since $(2-p)s^{1-p}>0$, \eqref{eq:barrier-theta-equation} gives
	$1>\theta(s)^2$, that is,
	$
	0<\theta(s)<1.
	$
	The defining equation also gives, for every $s>0$,
	\begin{equation}\label{eq:barrier-first-derivative}
		(2-p)s^{1-p}=\frac{1-\theta^2}{\theta},
		\qquad
		2\theta+(2-p)s^{1-p}
		=\frac{1+\theta^2}{\theta},
	\end{equation}
	where $\theta=\theta(s)$.  Differentiating
	\eqref{eq:barrier-theta-equation} with respect to $s$ and using
	\eqref{eq:barrier-first-derivative}, we obtain
	\begin{equation}\label{eq:barrier-first-derivative-gap}
		s\theta'(s)
		=(p-1)\frac{\theta(s)\left(1-\theta(s)^2\right)}
		{1+\theta(s)^2}>0 ;
	\end{equation}
	in particular $\theta$ is strictly increasing.  We now define the
	barrier $B$ by
	\[
	B(s)=s^{2-p}+s\theta(s)
	=\frac{s^{2-p}}2
	\left(p+\sqrt{(2-p)^2+4s^{2p-2}}\right).
	\]  We shall use
	\begin{equation}\label{eq:B-properties}
		B''(s)<0,
		\qquad
		B'(s)>\frac{s}{B(s)-s^{2-p}}=\frac{1}{\theta(s)}.
	\end{equation}
	Since $B(s)-s^{2-p}=s\theta(s)$ and, by
	\eqref{eq:barrier-first-derivative},
	$
	B'(s)=(2-p)s^{1-p}+\theta+s\theta'
	=\frac1\theta+s\theta',
	$
	\eqref{eq:barrier-first-derivative-gap} gives both the second
	inequality in \eqref{eq:B-properties} and
	\begin{equation}\label{eq:barrier-second-derivative}
		B'(s)
		=\frac1\theta
		+(p-1)\frac{\theta(1-\theta^2)}{1+\theta^2}.
	\end{equation}
	Differentiating \eqref{eq:barrier-second-derivative} therefore gives
	\[
	B''(s)
	=\theta'\left[
	-\frac1{\theta^2}
	+(p-1)\frac{1-4\theta^2-\theta^4}
	{(1+\theta^2)^2}\right]
	\qquad\text{with }\theta=\theta(s).
	\]
	Since $0<\theta<1$,
	$
	-1/\theta^2<-1.
	$
	Moreover,
	$
	(1+\theta^2)^2-(1-4\theta^2-\theta^4)
	=6\theta^2+2\theta^4>0,
	$
	so
	$
	(1-4\theta^2-\theta^4)/(1+\theta^2)^2<1.
	$
	Consequently,
	\[
	-\frac1{\theta^2}
	+(p-1)\frac{1-4\theta^2-\theta^4}
	{(1+\theta^2)^2}
	<-1+(p-1)=p-2<0.
	\]
	Together with \eqref{eq:barrier-first-derivative-gap}, this proves
	$B''(s)<0$ and hence \eqref{eq:B-properties}.
	
	\medskip
	\noindent
	\emph{Step 2. Comparison principle.}
	
	We claim that
	$
	A(r)<B(u(r))$ for $
	r\geq0.
	$
	Since $u,A\to0$ as $r\to+\infty$ and $B(s)\to0$ as $s\downarrow0$, one has
	$
	B(u(r))-A(r)\to0
	\ \text{as }r\to\infty.
	$
	If $B(u)-A$ were not strictly positive, it would attain a
	nonpositive global minimum at some $x_0\in\RN$.  At $x_0$, we have
	$B(u)\leq A$, and hence
	\[
	H=\frac{A-u^{2-p}}u
	\geq\frac{B(u)-u^{2-p}}u
	=\theta(u)>0.
	\]
	The second inequality in \eqref{eq:B-properties} therefore gives
	\begin{equation}\label{eq:barrier-laplacian}
		B'(u)H\geq B'(u)\theta(u)
		>\frac1{\theta(u)}\theta(u)=1.
	\end{equation}
	Moreover, the chain rule, $-\Delta u=u^pH$, and $-\Delta A=u^p$
	give
	\begin{equation}\label{eq:barrier-contact}
		\begin{aligned}
			-\Delta\bigl(B(u)-A\bigr)
			&=B'(u)(-\Delta u)-B''(u)|\nabla u|^2+\Delta A
			=u^p\bigl(B'(u)H-1\bigr)
			-B''(u)|\nabla u|^2.
		\end{aligned}
	\end{equation}
	By \eqref{eq:barrier-laplacian} and $B''<0$,
	\eqref{eq:barrier-contact} yields
	$
	-\Delta(B(u)-A)(x_0)>0.
	$
	On the other hand, since $x_0$ is a minimum point of $B(u)-A$, it holds
	$
	\Delta(B(u)-A)(x_0)\geq0,
	$
	a contradiction.  Therefore
	\[
	A(r)<B(u(r)),
	\quad r\geq0.
	\]
	It follows that
	\[
	H(r)
	=\frac{A(r)-u(r)^{2-p}}{u(r)}
	<
	\frac{B(u(r))-u(r)^{2-p}}{u(r)}
	=\theta(u(r)).
	\]
	Substituting $s=u(r)$ into \eqref{eq:barrier-theta-equation} gives
	$
	\theta(u(r))
	\left(\theta(u(r))+(2-p)u(r)^{1-p}\right)=1.
	$
	Therefore, wherever $H(r)>0$, we have
	$
	0<H(r)<\theta(u(r)).
	$
	Since $(2-p)u(r)^{1-p}>0$, it follows that
	\begin{equation}\label{eq:H-barrier}
		\begin{aligned}
			H(r)\left(H(r)+(2-p)u(r)^{1-p}\right)
			&<
			\theta(u(r))
			\left(\theta(u(r))+(2-p)u(r)^{1-p}\right)=1.
		\end{aligned}
	\end{equation}
	
	\medskip
	\noindent
	\emph{Step 3. The flux ratio.}
	
	For $r>0$, define the flux ratio
	\begin{equation}\label{eq:R-def}
		\mathcal F(r)=\frac{u'(r)}{A'(r)}.
	\end{equation}
	Radial integration of $-\Delta u=u^pH$ gives
	\[
	-r^{N-1}u'(r)
	=\int_0^r s^{N-1}u(s)^pH(s)\dd s.
	\]
	Combining this identity with \eqref{eq:A-flux}, and using
	$u'(r)<0$ and $A'(r)<0$, we obtain
	\begin{equation}\label{eq:R-average}
		\mathcal F(r)=
		\frac{\displaystyle\int_0^r
			s^{N-1}u(s)^pH(s)\dd s}
		{\displaystyle\int_0^r s^{N-1}u(s)^p\dd s}>0.
	\end{equation}
	Thus $\mathcal F(r)$ is the $u^p$-weighted average of $H$ over $(0,r)$.
	Differentiating \eqref{eq:R-average}, we have
	\begin{equation}\label{eq:R-derivative}
		\mathcal F'
		=\frac{r^{N-1}u^p}
		{\displaystyle\int_0^r s^{N-1}u^p\dd s}
		(H-\mathcal F).
	\end{equation}
	On the other hand,
	\[
	A=u^{2-p}+uH.
	\]
	Differentiating this identity and using $A'=u'/\mathcal F$, which follows from
	\eqref{eq:R-def}, we obtain
	\[
	\frac{u'}{\mathcal F}
	=(2-p)u^{1-p}u'+u'H+uH'.
	\]
	Solving for $H'$ gives
	\begin{equation}\label{eq:H-derivative}
		H'
		=-\frac{u'}u
		\left(
		H+(2-p)u^{1-p}-\frac1{\mathcal F}
		\right).
	\end{equation}
	
	We claim that
	\begin{equation}\label{eq:R-separation}
		\mathcal F(r)>H(r),
		\qquad r>0.
	\end{equation}
	First, \eqref{eq:R-average} and the continuity of $H$ give
	$
	\mathcal F(r)\to H(0)
	\ \text{as }r\downarrow0.
	$
	Since $\mathcal F(r)>0$ given in \eqref{eq:R-average}, we have $H(0)\geq0$. Since $u'(0)=A'(0)=0$ and $u(0)>0$, differentiating
	$
	H=A/u-u^{1-p}
	$
	gives
	$
	H'
	=\frac{A'}u-\frac{Au'}{u^2}-(1-p)u^{-p}u',
	$
	and hence $H'(0)=0$.
	Moreover, the equations at the origin give
	\[
	-Nu''(0)=u(0)^pH(0),
	\qquad
	-NA''(0)=u(0)^p.
	\]
	Since $H=A/u-u^{1-p}$, differentiating twice at the origin yields
	\[
	H''(0)
	=\frac{A''(0)}{u(0)}
	-\frac{A(0)u''(0)}{u(0)^2}
	-(1-p)u(0)^{-p}u''(0).
	\]
	Substituting
	$
	A''(0)=-u(0)^p/N
	$
	and
	$
	u''(0)=-u(0)^pH(0)/N
	$
	into the preceding formula for $H''(0)$ gives
	\begin{align*}
		H''(0)
		&=-\frac{u(0)^{p-1}}N
		+\frac{A(0)u(0)^{p-2}H(0)}N
		+\frac{(1-p)H(0)}N.
	\end{align*}
	Moreover,
	$
	A(0)=u(0)^{2-p}+u(0)H(0)
	$
	implies
	$
	A(0)u(0)^{p-2}=1+u(0)^{p-1}H(0).
	$
	Therefore
	\begin{equation}\label{eq:H-centre-expansion}
		\begin{aligned}
			H''(0)
			&=-\frac{u(0)^{p-1}}N
			+\frac{H(0)+u(0)^{p-1}H(0)^2}N
			+\frac{(1-p)H(0)}N\\
			&=-\frac{u(0)^{p-1}}N
			+\frac{u(0)^{p-1}H(0)^2+(2-p)H(0)}N\\
			&=-\frac{u(0)^{p-1}}N
			\left[
			1-H(0)\left(H(0)+(2-p)u(0)^{1-p}\right)
			\right].
		\end{aligned}
	\end{equation}
	If $H(0)=0$, then \eqref{eq:H-centre-expansion} gives
	$H''(0)<0$.  Since $H'(0)=0$, this implies that $H(r)<0$ for all
	sufficiently small $r>0$.  Its positive-weight average
	\eqref{eq:R-average} would then also be negative, contradicting
	$\mathcal F(r)>0$.  Therefore $H(0)>0$.
	Applying \eqref{eq:H-barrier} at the origin, we have
	\[
	H(0)\left(H(0)+(2-p)u(0)^{1-p}\right)<1.
	\]
	Equation \eqref{eq:H-centre-expansion} therefore gives $H''(0)<0$.
	Since $H'(0)=0$ and $H''$ is continuous, $H'(r)<0$ for all
	sufficiently small $r>0$.  Hence $H$ is strictly decreasing near the
	origin.  For every sufficiently small $r>0$, \eqref{eq:R-average}
	then gives
	\[
	\mathcal F(r)-H(r)
	=
	\frac{\displaystyle\int_0^r
		s^{N-1}u(s)^p\bigl(H(s)-H(r)\bigr)\dd s}
	{\displaystyle\int_0^r s^{N-1}u(s)^p\dd s}>0.
	\]
	Suppose that \eqref{eq:R-separation} fails, and let $r_0>0$ be the
	first point at which equality occurs.  Then
	\[
	\mathcal F(r_0)=H(r_0)>0.
	\]
	Equation \eqref{eq:R-derivative} gives $\mathcal F'(r_0)=0$.  Since
	$\mathcal F-H>0$ on $(0,r_0)$ and $\mathcal F(r_0)-H(r_0)=0$, we have
	$
	(\mathcal F-H)'(r_0)\leq0.
	$
	It follows from $\mathcal F'(r_0)=0$ that $H'(r_0)\geq0$.  On the other hand,
	\eqref{eq:H-barrier} and $H(r_0)>0$ give
	\[
	H(r_0)+(2-p)u(r_0)^{1-p}-\frac1{H(r_0)}<0.
	\]
	Since $\mathcal F(r_0)=H(r_0)$ and $-u'(r_0)/u(r_0)>0$,
	\eqref{eq:H-derivative} gives $H'(r_0)<0$, a contradiction.  This
	proves \eqref{eq:R-separation}.  Substituting $\mathcal F>H$ into
	\eqref{eq:R-derivative} also gives
	\begin{equation}\label{eq:R-decreasing}
		\mathcal F'(r)<0,
		\qquad r>0.
	\end{equation}
	
	\medskip
	\noindent
	\emph{Step 4. Transfer to the $u^2$-weighted average of $V$.}
	
	We prove that $\bar V(r)>V_*$ for every $r>0$, which is
	\eqref{eq:prefix-negative}.  By Lemma~\ref{lem:local-pohozaev},
	this will imply $E(r)>0$.  Fix $r>0$.  If $0<t\leq r$,
	then \eqref{eq:R-separation} and \eqref{eq:R-decreasing} give
	$
	\mathcal F(t)\geq \mathcal F(r)>H(r).
	$
	Using \eqref{eq:R-average}, we therefore obtain
	\begin{align}
		&\int_0^t s^{N-1}u(s)^p
		\bigl(H(s)-H(r)\bigr)\dd s =
		\left(\int_0^t s^{N-1}u(s)^p\dd s\right)
		\bigl(\mathcal F(t)-H(r)\bigr)>0.
		\label{eq:H-prefix-positive}
	\end{align}
	Therefore, integration by parts gives
	\begin{equation}\label{eq:Stieltjes-positive}
		\begin{aligned}
			&\int_0^r t^{N-1}u(t)^{p+1}
			\bigl(H(t)-H(r)\bigr)\dd t\\
			&\quad=
			\int_0^r u(t)\frac{\dd}{\dd t}
			\left[
			\int_0^t s^{N-1}u(s)^p
			\bigl(H(s)-H(r)\bigr)\dd s
			\right]\dd t\\
			& \quad=
			u(r)\int_0^r s^{N-1}u(s)^p
			\bigl(H(s)-H(r)\bigr)\dd s
			-\int_0^r u'(t)
			\left[
			\int_0^t s^{N-1}u(s)^p
			\bigl(H(s)-H(r)\bigr)\dd s
			\right]\dd t>0.
		\end{aligned}
	\end{equation}
	The strict inequality follows from
	\eqref{eq:H-prefix-positive}, $u(r)>0$, and $u'<0$.
	
	We next prove that the running average exceeds the pointwise value,
	\begin{equation}\label{eq:V-endpoint-comparison}
		\int_0^r s^{N-1}u(s)^2
		\bigl(V(s)-V(r)\bigr)\dd s>0,
	\end{equation}
	which is the same as $\bar V(r)>V(r)$ after division by the positive
	integral $\int_0^r s^{N-1}u(s)^2\dd s$.
	Suppose first that $H(r)\geq0$.  Since
	$
	V-1=u^{p-1}H
	$ given in \eqref{eq:V-H},
	we have
	\begin{align*}
		\int_0^r s^{N-1}u(s)^2
		\bigl(V(s)-V(r)\bigr)\dd s &=
		\int_0^r s^{N-1}u(s)^{p+1}
		\bigl(H(s)-H(r)\bigr)\dd s\\
		&\quad
		+H(r)\int_0^r s^{N-1}u(s)^2
		\bigl(u(s)^{p-1}-u(r)^{p-1}\bigr)\dd s,
	\end{align*}
	the two constants $1$ supplied by \eqref{eq:V-H} cancelling.
	The first term on the right is strictly positive by
	\eqref{eq:Stieltjes-positive}.  The second is nonnegative when
	$H(r)\geq0$.  Thus
	\eqref{eq:V-endpoint-comparison} holds when $H(r)\geq0$.
	
	Suppose now that $H(r)<0$.  Then
	$
	V(r)-1=u(r)^{p-1}H(r)<0.
	$
	The radial equation $-\Delta u=u(V-1)$ can be written as
	$
	-\bigl(s^{N-1}u'(s)\bigr)'
	=s^{N-1}u(s)(V(s)-1).
	$
	Multiplying this identity by $u(s)$ and integrating from $0$ to $r$,
	we obtain
	\[
	\begin{aligned}
		\int_0^r s^{N-1}u(s)^2(V(s)-1)\dd s
		&=-\int_0^r u(s)\bigl(s^{N-1}u'(s)\bigr)'\dd s\\
		&=-\left[s^{N-1}u(s)u'(s)\right]_{0}^{r}
		+\int_0^r s^{N-1}(u'(s))^2\dd s\\
		&=-r^{N-1}u(r)u'(r)
		+\int_0^r s^{N-1}(u'(s))^2\dd s>0.
	\end{aligned}
	\]
	Since $V(r)-1<0$, this also gives \eqref{eq:V-endpoint-comparison}.
	Combining these two cases and dividing  by the positive quantity 
	$\int_0^r s^{N-1}u(s)^2\dd s>0$, we obtain
	\begin{equation}\label{eq:V-running-average}
		\bar V(r)>V(r),\qquad r>0 .
	\end{equation}
	
	Differentiating \eqref{eq:V-bar-def} and using
	\eqref{eq:V-running-average}, we obtain
	\[
	\bar V'(r)=
	\frac{r^{N-1}u(r)^2}
	{\displaystyle\int_0^r s^{N-1}u(s)^2\dd s}
	\bigl(V(r)-\bar V(r)\bigr)<0,
	\]
	so that $\bar V$ is strictly decreasing on $(0,\infty)$.  Therefore, at every $r>0$,
	\[
	\bar V(r)>\lim_{\rho\to\infty}\bar V(\rho)=V_* .
	\]
	After multiplying by the positive denominator, and using \(V_*>0\),
	this is exactly \eqref{eq:prefix-negative}.
	Lemma~\ref{lem:local-pohozaev} now gives $E(r)>0$ for every $r>0$.
\end{proof}

		\noindent{\bf Acknowledgements}  
	This work was partially supported by the National Natural Science Foundation of China (Grant Nos. 12271436, 12371107 and 11901582). The authors used OpenAI language models as auxiliary tools in preparing this manuscript; all mathematical content was independently verified and is the sole responsibility of the authors.
		\vskip.1in
		\noindent{\bf  Data availability statement} Data sharing is not applicable to this article, as no datasets were generated or analyzed during the current study.
		
		\vskip.1in
		\noindent{\bf  Conflict of interest statement}  
		On behalf of all authors, the corresponding author states that there is no conflict of interest.

\bigskip
\addcontentsline{toc}{section}{References}


	\end{document}